\documentclass[12pt]{article}
\usepackage{amsmath}
\usepackage{graphicx,psfrag,epsf}
\usepackage{enumerate}
\usepackage{enumitem}
\usepackage{natbib}
\usepackage{url} 
\usepackage{amssymb}
\usepackage{algpseudocode}
\usepackage{algorithmicx}
\usepackage{algorithm2e}
\usepackage{color}
\usepackage{xcolor}
\usepackage{mathtools}
\usepackage{siunitx}
\usepackage{relsize}
\usepackage{multirow}
\usepackage{tikz}

\usepackage[colorlinks]{hyperref}
\hypersetup{
  colorlinks=true,
  linkcolor=red,
  citecolor=red
}
\AtBeginDocument{%
  \hypersetup{%
    linkcolor=red,%
    citecolor=blue,%
  }%
}%

\renewcommand{\sectionautorefname}{Section}
\renewcommand{\subsectionautorefname}{Section}

\let\orgautoref\autoref
\providecommand{\autorefs}
        {\def\equationautorefname{Equations}%
         \def\figureautorefname{Figures}%
         \def\subfigureautorefname{Figures}%
         \def\sectionautorefname{Sections}%
         \def\subsectionautorefname{Sections}%
         \def\subsubsectionautorefname{Sections}%
         \def\Itemautorefname{items}%
         \def\tableautorefname{Tables}%
         \orgautoref}
         
\renewcommand{\autoref}
        {\def\equationautorefname{Equation}%
         \def\figureautorefname{Figure}%
         \def\subfigureautorefname{Figure}%
         \def\sectionautorefname{Section}%
         \def\subsectionautorefname{Section}%
         \def\subsubsectionautorefname{Section}%
         \def\Itemautorefname{item}%
         \def\tableautorefname{Table}%
         \def\chapterautorefname{Chapter}%
         \orgautoref}

\def\bSig\mathbf{\Sigma}

\newcommand{\widesim}[2][1.5]{
  \mathrel{\overset{#2}{\scalebox{#1}[1]{$\sim$}}}
}
\newcommand{\m}[1]{\boldsymbol{#1}}

\newcommand*{\spr}[3]{\left\langle \,#1\, \middle| \,#2\, \right\rangle_{#3}}
\DeclareMathOperator{\spn}{span}

\newcommand{\blind}{0}

\date{}

\begin{document}

\def\spacingset#1{\renewcommand{\baselinestretch}%
{#1}\small\normalsize} \spacingset{1}


\if0\blind
{
\title{\bf Component--based regularisation of multivariate generalised linear mixed models}
\author{Jocelyn Chauvet
\thanks{IMAG, Univ Montpellier, CNRS, Montpellier, France. \newline 
\href{mailto:jocelyn.chauvet@umontpellier.fr}{jocelyn.chauvet@umontpellier.fr} ; \href{mailto:xavier.bry@umontpellier.fr}{xavier.bry@umontpellier.fr}}
    \\
    Catherine Trottier\thanks{Univ Paul--Val\'ery Montpellier 3, Montpellier, France. \newline 
\href{mailto:catherine.trottier@univ-montp3.fr}{catherine.trottier@univ-montp3.fr}} \footnotemark[1]
    \\
    and \\
    Xavier Bry\footnotemark[1]
    }
  \maketitle
} \fi

\if1\blind
{
  \bigskip
  \bigskip
  \bigskip
  \begin{center}
    {\LARGE\bf Component--based regularisation of multivariate generalised linear mixed models}
\end{center}
  \medskip
} \fi

\bigskip
\begin{abstract}
We address the component--based regularisation of a multivariate Generalised Linear Mixed Model (GLMM) in the framework of grouped data. 
A set $\m{Y}$ of random responses is modelled with a multivariate GLMM, based on a set $\m{X}$ of explanatory variables, a set $\m{A}$ of additional explanatory variables, and random effects to introduce the within--group dependence of observations.
Variables in $\m{X}$ are assumed many and redundant so that regression demands regularisation.  
This is not the case for $\m{A}$, which contains few and selected variables.
Regularisation is performed building an appropriate number of orthogonal components that both contribute to model $\m{Y}$ and capture relevant structural information in $\m{X}$.
To estimate the model, we propose to maximise a criterion specific to the Supervised Component--based Generalised Linear Regression (SCGLR) within an adaptation of Schall's algorithm.
This extension of SCGLR is tested on both simulated and real grouped data, and compared to ridge and LASSO regularisations.
Supplementary material for this article is available online.
\end{abstract}

\noindent%
\vfill

\newpage
\spacingset{1.45} 
\section{Introduction}
\label{sec:intro}

In the framework of regression models on a large number of explanatory variables with redundancies and collinearities, the search for a reduced number of relevant dimensions to model responses has been an ongoing research over the last decades. 
In particular, the case where 
the explanatory variables outnumber the observations tends to be a new standard.
Generalised Linear Models (GLMs) are the most widely used regression models, 
because they are easy to interpret and address a very large scope of applications with a variety of response distributions. 
For instance, Epidemiology, Biology and Social Sciences need to model binary outcomes, count data and survival times. All these fields often have to deal with both grouped data and multivariate responses combining variables of different types (e.g. one binary and another Poisson). 
In this work, we particularly aim at modelling abundances of several tree genera 
on plots of land grouped into forest concessions, using multiple redundant explanatory variables.

\paragraph*{}
As far as dimension--reduction is concerned, two main approaches have been developed. 
The first one is variable--selection, whereas the second one builds components, 
i.e. linear combinations of the explanatory variables, which synthesise the useful part of their information. 
As far as variable--selection is concerned, the most popular method is currently the LASSO, introduced by \citet{tibshirani1996}, which combines the likelihood with a penalty based on the 
$L_1$--norm of the coefficient vector.
LASSO is one of the penalty--based regularisation methods, 
as are also ridge \citep{hoerl1970} and elastic--net \citep{zou2005}. 
This LASSO selection approach has proved efficient to explain the phenomenon of interest when some of the explanatory variables are the ``true'' ones, surrounded by a high number of irrelevant others. Nevertheless, it may be very unstable and helpless when the true explanatory dimensions are latent and indirectly measured through highly correlated proxies. 
This is where the component--based approach turns out to be useful. 
\citet{bry2013} have developed a new methodology 
named Supervised Component--based Generalised Linear Regression (SCGLR), later extended and refined in \citet{bry2014,bry2016,bry2018}.
As in any PLS--type method, the construction of components in SCGLR is guided both by the correlation--structure of variables in  the explanatory space and by the prediction quality of the responses. 
Nevertheless, unlike PLS, SCGLR involves a general and flexible criterion allowing to specify the type of structure components are wanted to align with in the explanatory space (e.g. variable bundles, principal components, other subspaces). 
Moreover, SCGLR searches for  explanatory directions common to multiple responses with probability distributions in the exponential family, each response being entitled to their own distribution.
The current SCGLR method is implemented in the R package \texttt{SCGLR} \citep{packageSCGLR} available at \url{https://CRAN.R-project.org/package=SCGLR} and 
\url{https://github.com/SCnext/SCGLR}.

\paragraph*{}
In the present work, we aim at modelling responses with a repeated or grouped design.
For this purpose, the use of mixed models with random effects is widespread.
Research on variance--component estimation in Generalised Linear Mixed Models (GLMMs) has been very active since the 1980s.
For the most general distribution assumptions in such models, parameter estimation faces the  intractability of the likelihood expressed as an integral with respect to the random effects. Several numerical approximations of the integral have been proposed: Gaussian quadrature \citep{anderson1985} or  adaptive versions of it \citep{pinheiro1995}, Laplace approximation leading to the definition of the penalised quasi--likelihood \citep{breslow1993} or modified versions of it \citep{shun1995}. An alternative to this type of analytic approximation is a stochastic approximation of the integral calculation via MCMC techniques. In this approach, \citet{zeger1991} described an approximate Gibbs sampling for GLMMs, which was extended by \citet{clayton1996} to more general Metropolis--Hastings algorithms. In parallel, \citet{mcCulloch1997} developed the Monte Carlo EM algorithm where the expectation is computed numerically through a Monte Carlo approximation, after generating random effects with a Metropolis--Hastings sampler. 
Mention can also be made of the recent work by \citet{knudson2016}: 
her strategy is to approximate the entire likelihood function using random effects simulated from a parametrised importance sampling distribution depending on the data. Unfortunately, these different approaches are not necessarily suitable for the same types of random effect designs (one--dimensional random effect, embedded random effects, etc).
In the wake of the first type of approximations, we here adopt the ``Joint--Maximisation'' strategy \citep{mcCulloch1997}, as introduced for instance by \citet{schall1991}.  
The model is iteratively linearised conditional on the random effects and  variance components are then estimated using adapted linear mixed models methods. This strategy can be used for any random effect design and is less computationally intensive than Monte Carlo methods.
Moreover, it provides us with a linear setting more suitable for the computation of components. Once the components calculated, model parameters can be estimated using any of the aforementioned strategies (see \autoref{sec:conc}).

\paragraph*{}
Modelling grouped responses through a GLMM with a large number of explanatory variables is the focus of this paper. The need for dimension--reduction and regularisation has to accommodate the presence of random effects in the model, but our main purpose still remains to investigate the explanatory structure and link it to interpretable dimensions.
For Gaussian responses, \citet{eliot2011} proposed to extend the ridge regression to Linear Mixed Models (LMMs). Based on a penalised complete log--likelihood, the adaptation of the Expectation--Maximisation algorithm they suggest includes a new step to find the best shrinkage parameter using a generalised cross--validation scheme at each iteration. 
More recently, \cite{schelldorfer2014} --- and also \citet{groll2014}  ---  proposed an $L_1$--penalised algorithm for fitting a high--dimensional GLMM, using Laplace approximation and an efficient coordinate gradient descent. 
In this work, we combine Schall's iterative model linearisation with regularisation at each step. 
However, we do not use a penalty on the coefficient vector's norm --- as proposed by \citet{zhang2017} within the framework of multivariate count data. We rather propose to combine dimension--reduction and predictor--regularisation using supervised components aligning on the most predictive and interpretable directions in the explanatory space.

\paragraph*{}
The paper is organised as follows. In \autoref{Chauvet::Model_definition}, we formalise the model and set the main notations used throughout the paper.
In \autoref{Chauvet::SCGLR}, we present the key features of SCGLR. 
\autoref{Chauvet::Group_SCGLR} designs an extension of this methodology to mixed models, and particularly to grouped data. 
In \autoref{Chauvet::simulatedDATA}, our extended method ``mixed--SCGLR'' is evaluated on simulations and compared to ridge-- and LASSO--based regularisations. 
Finally, in order to highlight the power of mixed--SCGLR in terms of model interpretation, \autoref{Chauvet::realDATA} presents an application to real data in the Poisson case.

\section{Model definition and notations}
\label{Chauvet::Model_definition}

In the framework of a multivariate GLMM, we consider $q$ response--vectors
$ \boldsymbol{y_1}, \ldots, \boldsymbol{y_q}$ forming matrix 
$\boldsymbol{Y}_{n \times q}$, to be explained by two categories of explanatory variables. 
The first category consists of few weakly correlated variables
$\boldsymbol{A}_{n \times r}=
\big[\, \boldsymbol{a_1} \mid \ldots \mid \boldsymbol{a_r} \,\big]$.
These variables are assumed to be interesting per se and 
their marginal effects need to be precisely quantified. 
The second category consists of abundant and highly correlated variables
$\boldsymbol{X}_{n \times p}=
\big[\, \boldsymbol{x_1} \mid \ldots \mid \boldsymbol{x_p} \,\big]$ 
considered as proxies to latent dimensions which must be found and interpreted.
Since explanatory variables in $\m{A}$ are few, non--redundant and of interest, they are kept as such in the model. 
By contrast, $\m{X}$ may contain several unknown structurally relevant dimensions $K<p$ important to model and predict $\m{Y}$, how many we do not know.
$\m{X}$ is thus to be searched for an appropriate number of orthogonal components that both capture relevant structural information in $\m{X}$ and contribute to model $\m{Y}$.

\paragraph*{}
This work addresses grouped data: 
the $n$ observations form $N$ groups. 
Within each group, observations are not assumed independent.
For each response $\m{y_k}$, a $N$--level random effect 
$\m{\xi_k}$ is used to model the dependence of observations within each group. 
Hence, each $\m{y_k}$ is modelled with a GLMM assuming a conditional distribution from the exponential family.

\paragraph*{Notations and conventions}
\begin{itemize}[parsep=0cm,itemsep=0cm,topsep=0cm]
\item[$\blacktriangleright$]
All variables 
(namely the $\m{a_i}$'s, $\m{x_j}$'s and $\m{y_k}$'s)
will be identified with $n$--vectors. 
\item[$\blacktriangleright$]
We will use bold lowercase letters for vectors (e.g. $\m{u}$) and bold capital letters for matrices (e.g. $\m{M}$).
\item[$\blacktriangleright$]
$\m{M}$ being any matrix, ${\m{M}}^{\mbox{\tiny T}}$ denotes the transpose of $\m{M}$.
\item[$\blacktriangleright$]
$\m{I}_n$ denotes the identity matrix of size $n$.
\item[$\blacktriangleright$]
$\mathbf{1}_m$ denotes the all--ones vector of size $m$.
\item[$\blacktriangleright$]
Let $\m{u}$ and $\m{v}$ be non--zero vectors in $\mathbb{R}^d$ and let $\m{M}$ be a symmetric positive definite matrix of size $d \times d$.
Then $\spr{\m{u}}{\m{v}}{\m{M}} = {\m{u}}^{\mbox{\tiny T}} \m{M} \m{v}$ 
refers to the Euclidean scalar product of $\m{u}$ and $\m{v}$ with respect to metric $\m{M}$. The cosine of the angle between $\m{u}$ and $\m{v}$ with respect to $\m{M}$ is given by
$\cos_{\m{M}} \left( \m{u}, \m{v} \right) = 
\dfrac{\spr{\m{u}}{\m{v}}{\m{M}}}
{\left\lVert \m{u} \right\rVert_{\m{M}} \left\lVert \m{v} \right\rVert_{\m{M}}}$, where $\left\lVert \m{u} \right\rVert_{\m{M}} = \sqrt{\spr{\m{u}}{\m{u}}{\m{M}}}$.
\item[$\blacktriangleright$]
The space spanned by vectors $\m{u_1}, \ldots, \m{u_h}$ is denoted by 
$\spn \left\lbrace \m{u_1}, \ldots, \m{u_h} \right\rbrace$.
$\m{U}$ being any matrix, $\spn \left\lbrace \m{U} \right\rbrace$ refers to the space spanned by the column--vectors of $\m{U}$.
\item[$\blacktriangleright$]
Let $\mathbb{R}^n$ be endowed with metric $\m{W}$ and let $\m{Z}$ be a matrix of size $n \times p$.
Then $\Pi_{\spn \left\lbrace \m{Z} \right\rbrace}^{\m{W}}$ refers to the 
$\m{W}$--orthogonal projector onto $\spn \left\lbrace \m{Z} \right\rbrace$.
Let $\m{b}$ be a vector in $\mathbb{R}^n$.
The cosine of the angle between $\m{b}$ and 
$\spn \left\lbrace \m{Z} \right\rbrace$ with respect to $\m{W}$ is defined by 
$\cos_{\m{W}} \left( \m{b},\, \spn \left\lbrace \m{Z} \right\rbrace \right) 
= \cos_{\m{W}} \left( \m{b},\, 
\Pi_{\spn \left\lbrace \m{Z} \right\rbrace}^{\m{W}} \m{b} \right)$.
\end{itemize}

\section{SCGLR with additional explanatory variables}
\label{Chauvet::SCGLR}

In this section we consider the situation where each 
$\m{y_k}$ is modelled with a GLM (without random effect).
For the sake of simplicity, we focus on the single--component SCGLR ($K=1$).
\autoref{Chauvet::background_GLM} briefly recalls some standards for univariate GLMs.
\autoref{Chauvet::New_linear_predictors} defines the linear predictors considered in the SCGLR methodology, in a multivariate GLM framework with additional explanatory variables.
Finally, \autoref{Chauvet::criterion} introduces the criterion SCGLR maximises to compute the component.

\subsection{Notations and main features of univariate GLMs}
\label{Chauvet::background_GLM}

We refer the reader to \citet{mcCullagh1989} for a thorough overview of GLMs.
This section is only intended to recall the classical iterative scheme performing maximum likelihood (ML) estimation.
Let $\m{X}$ denote the $n \times p$ matrix of explanatory variables 
and $\boldsymbol{\beta}$ the $p$--dimensional parameter vector. 
At iteration $t+1$, the Fisher Scoring Algorithm (FSA) for ML estimation calculates
\begin{equation}
{\m \beta}^{[t+1]}
= \left( {\m X}^{\mbox{\tiny T}} {\m W}^{[t]} {\m X} \right)^{-1} 
{\m X}^{\mbox{\tiny T}} {\m W}^{[t]} {\m z}^{[t]}, 
\label{Chauvet::GLM_FSA}
\end{equation}
where ${\m z}^{[t]}$ and ${\m W}^{[t]}$ respectively denote the classical working variable and the associated weight matrix at iteration $t$.
As pointed out by \citet{nelder1972}, update (\ref{Chauvet::GLM_FSA}) 
may be interpreted as a weighted least squares step in the linearised model $\mathcal{M}^{[t]}$ defined by
\begin{equation}
\mathcal{M}^{[t]}: \;
\left|
\begin{aligned}
&{\m z}^{[t]} = {\m X} {\m \beta} + {\m \zeta}^{[t]} \\
&\mbox{with: }
\mathbb{E} \left( {\m \zeta}^{[t]} \right) = {\m 0} \mbox{ and }
\mathbb{V} \left( {\m \zeta}^{[t]} \right) = {{\m W}^{[t]}}^{-1}.
\end{aligned}
\right.
\label{Chauvet::GLM_model_Mt}
\end{equation}

\subsection{Linear predictors for SCGLR with multiple responses}
\label{Chauvet::New_linear_predictors}

We are now considering a multivariate GLM \citep{fahrmeir1994}.
In this context, SCGLR searches for a 
component common to all the $\m{y_k}$'s. 
This component will be denoted $\m{f}$ and its 
$p$--dimensional loading--vector will be denoted $\m{u}$, so that 
$\m{f} = \m{Xu}$.
The linear predictor associated with response--vector $\m{y_k}$ then writes
\begin{equation} 
\m{\eta_k} = \left( \m{Xu} \right)\gamma_k + \m{A\delta_k}, 
\label{Chauvet::firstpredictor}
\end{equation}
where $\gamma_k$ and $\m{\delta_k}$ are the regression parameters associated respectively with component $\m{f}$ and additional explanatory variables $\m{A}$. 
$\m{f}$ being common to all the $\m{y_k}$'s, predictors are collinear in their $\m{X}$--part.
For identification purposes, we impose 
${\m{u}}^{\mbox{\tiny T}} \m{M}^{-1} \m{u} = 1$, 
where $\m{M}$ may so far be any $p \times p$ symmetric positive definite matrix.
Let us note $y_{k,i}$ the $i$--th observation of the $k$--th response--vector and 
$\m{H} = \left\lbrace \eta_{k,i} \mid 
1 \leqslant k \leqslant q, 1 \leqslant i \leqslant n \right\rbrace$ 
the predictor set.
We assume that the $q$ responses are independent conditional on $\m{f}$, 
and that the $n$ observations are independent. 
The log--density then writes
\begin{equation*}
\ell \left(\m{Y} \lvert \m{H} \right) = \sum_{i=1}^n \sum_{k=1}^q \ell_k 
\left( y_{k,i} \lvert \eta_{k,i} \right),
\end{equation*}
where $\ell_k$ is the log--density of the $k$--th response, conditional on its linear predictor.
As a result, $\m{z_k}$ being the working variable associated with $\m{y_k}$ 
and $\m{W}_{\m{k}}^{-1}$ its variance matrix, the corresponding linearised model derived from the FSA at iteration $t$ is
\begin{equation}
\mathcal{M}_k^{[t]}: \;
\left|
\begin{aligned}
&\m{z}_{\m{k}}^{[t]} = \left( \m{Xu} \right)\gamma_k + \m{A\delta_k} + 
\m{\zeta}_{\m{k}}^{[t]}  \\
&\mbox{with: }
\mathbb{E} \left( \m{\zeta}_{\m{k}}^{[t]} \right) = 0 \mbox{ and }
\mathbb{V} \left( \m{\zeta}_{\m{k}}^{[t]} \right) = 
{\m{W}_{\m{k}}^{[t]}}^{-1}.
\end{aligned}
\right.
\label{Chauvet::GLM_model_Mkt}
\end{equation}

Although linearised models (\ref{Chauvet::GLM_model_Mt}) and (\ref{Chauvet::GLM_model_Mkt}) seem very similar, (\ref{Chauvet::GLM_model_Mkt}) is no longer linear, owing to the product $\m{u}\gamma_k$. 
An alternate version of the FSA must therefore be used:
\begin{itemize}
\item[(i)] 
Given current values of all the $\gamma_k$'s and $\m{\delta_k}$'s, 
a new loading--vector $\m{u}$ is obtained by solving an SCGLR--specific program (see \autoref{Chauvet::criterion} for details).
\item[(ii)]
Given a current value of $\m u$, 
each $\m{z}_{\m{k}}$ is regressed independently on 
$\big[\, \m{Xu} \mid \m{A} \,\big]$ with respect to weight matrix 
$\m{W_k}$, yielding new regression parameters 
$\gamma_k$ and $\m{\delta_k}$.
\end{itemize}

\subsection{Calculating the component maximising an SCGLR--specific criterion}
\label{Chauvet::criterion}

For an easier reading of this part, we omit the $[t]$ index. 
For each $k \in \left\lbrace 1,\ldots,q \right\rbrace$, consider model 
$\mathcal{M}_k$ endowed with weight matrix $\m{W_k}$.
As suggested in \cite{bry2013},
the best loading--vector in the weighted least--squares sense would be the solution of
\begin{equation*}
\underset{\m{u}:\,
\m{u}^{\mbox{\tiny T}}\m{M}^{-1}\m{u}=1}{\min} \;\, 
\sum_{k=1}^q  
\left\lVert 
\m{z_k} - 
\Pi_{\spn \left\lbrace \m{Xu},\m{A} \right\rbrace}^{\m{W_k}} \m{z_k} 
\right\rVert_{\m{W_k}}^2 
\quad  \Longleftrightarrow  \quad
\underset{\m{u}:\,
\m{u}^{\mbox{\tiny T}}\m{M}^{-1}\m{u}=1}{\max}  \;\,
\sum_{k=1}^q  \left\lVert 
\Pi_{\spn \left\lbrace \m{Xu},\m{A} \right\rbrace}^{\m{W_k}} \m{z_k} 
\right\rVert_{\m{W_k}}^2. 
\end{equation*}
The maximisation program also writes \; 
$\underset{\m{u}:\,\m{u}^{\mbox{\tiny T}}\m{M}^{-1}\m{u}=1}{\max} \;\, 
\psi_{\m{A}} \left( \m{u} \right)$, where
\begin{align}
\psi_{\m{A}} \left( \m{u} \right) 
&= 
\sum_{k=1}^q \left\lVert \m{z_k}^{\vphantom{t}} \right\rVert_{\m{W_k}}^2  
\cos_{\m{W_k}}^2 
\Big( \m{z_k},\, \spn \left\lbrace \m{Xu},\m{A} \right\rbrace \Big)
\nonumber \\
&=
\sum_{k=1}^q \left\lVert \m{z_k}^{\vphantom{t}} \right\rVert_{\m{W_k}}^2  
\cos_{\m{W_k}}^2 
\left( \m{z_k}, \,
\Pi_{\spn \left\lbrace \m{Xu},\m{A} \right\rbrace}^{\m{W_k}} \m{z_k} 
\right).
\label{Chauvet::GoF_criterion}
\end{align}
Now, $\psi_{\m{A}}$ is a mere goodness--of--fit (GoF) measure that does not take into account the closeness of component $\m{f}=\m{Xu}$ to interpretable directions 
in $\m{X}$. 
The GoF measure, $\psi_{\m{A}}$, must therefore be combined with a measure $\phi$ of structural relevance (SR).

\paragraph*{}
Assume matrix $\m{X}$ consists of $p$ standardised numeric variables. 
Consider a weight system $\m{\omega} = \left\lbrace \omega_1, \ldots, \omega_p \right\rbrace$
--- e.g. $\omega_j = \frac{1}{p} \; \forall j \in \left\lbrace 1,\ldots,p \right\rbrace$ --- 
reflecting the a priori relative importance of variables.
Also consider a weight matrix $\m{P}$ 
--- e.g. $\m{P}=\frac{1}{n} \m{I}_n$ --- 
reflecting the a priori relative importance of observations. 
We define the most structurally relevant loading--vector as the solution of
\begin{equation*}
\underset{\m{u}:\,\m{u}^{\mbox{\tiny T}}\m{M}^{-1}\m{u}=1}{\max} \;\, 
\phi \left( \m{u} \right),
\end{equation*}
where
\begin{equation}
\phi \left( \m{u} \right)
= \left[    
\sum_{j=1}^p \omega_j 
\left(
\spr{\m{Xu}}{\m{x_j}}{\m{P}}^2 
\right)^l
\right]^{\frac{1}{l}}
= \left[   
\sum_{j=1}^p \omega_j 
\Big( 
\m{u}^{\mbox{\tiny T}} \, 
\m{X}^{\mbox{\tiny T}} \m{P} \m{x_j} \m{x}_{\m{j}}^{\mbox{\tiny T}} \m{P X} \, 
\m{u} 
\Big)^{l}      
\right]^{\frac{1}{l}}, \; l \geqslant 1,
\label{Chauvet::SR_criterion_particulier}
\end{equation}
for the scalar product is commutative. 
Formula (\ref{Chauvet::SR_criterion_particulier}) is in fact a particular case of the SR criterion proposed by \cite{bry2015theme,bry2016}. 
It can be viewed as a generalised average version of the usual dual PCA criterion:
$\sum_{j=1}^p \cos^2_{\m{P}} \left( \m{Xu}, \m{x_j} \right) = 
\sum_{j=1}^p \spr{\m{Xu}}{\m{x_j}}{\m{P}}^2$.
For $\m{M}=\left( \m{X}^{\mbox{\tiny T}} \m{PX} \right)^{-1}$,
(\ref{Chauvet::SR_criterion_particulier}) is called ``Variable--Powered Inertia'' (VPI). 
It should be stressed that for $\m{X}^{\mbox{\tiny T}} \m{PX}$ to be invertible, $\m{X}$ must be a column full rank matrix.
In case of strict collinearities within $\m{X}$,
as it always happens in high--dimensional settings,
we replace $\m{X}$ with the matrix $\m{C}$ of its principal components associated with non--zero eigenvalues. 
The component is then sought as $\m{f}=\m{Cu}$. 
We have $\m{C}=\m{XV}$, where $\m{V}$ is the matrix of corresponding unit-eigenvectors. 
Then, $\m{f}=\m{Cu}=\m{X\widetilde{u}}$ with $\m{\widetilde{u}}=\m{Vu}$.  
\cite{bry2018} show that among all loading--vectors $\m{t}$ such that $\m{Xt}=\m{f}$, 
$\m{\widetilde{u}}$ is that which has the minimum $L_2$--norm.

\begin{figure}[!ht]
\centering
\includegraphics[width=0.8\textwidth]{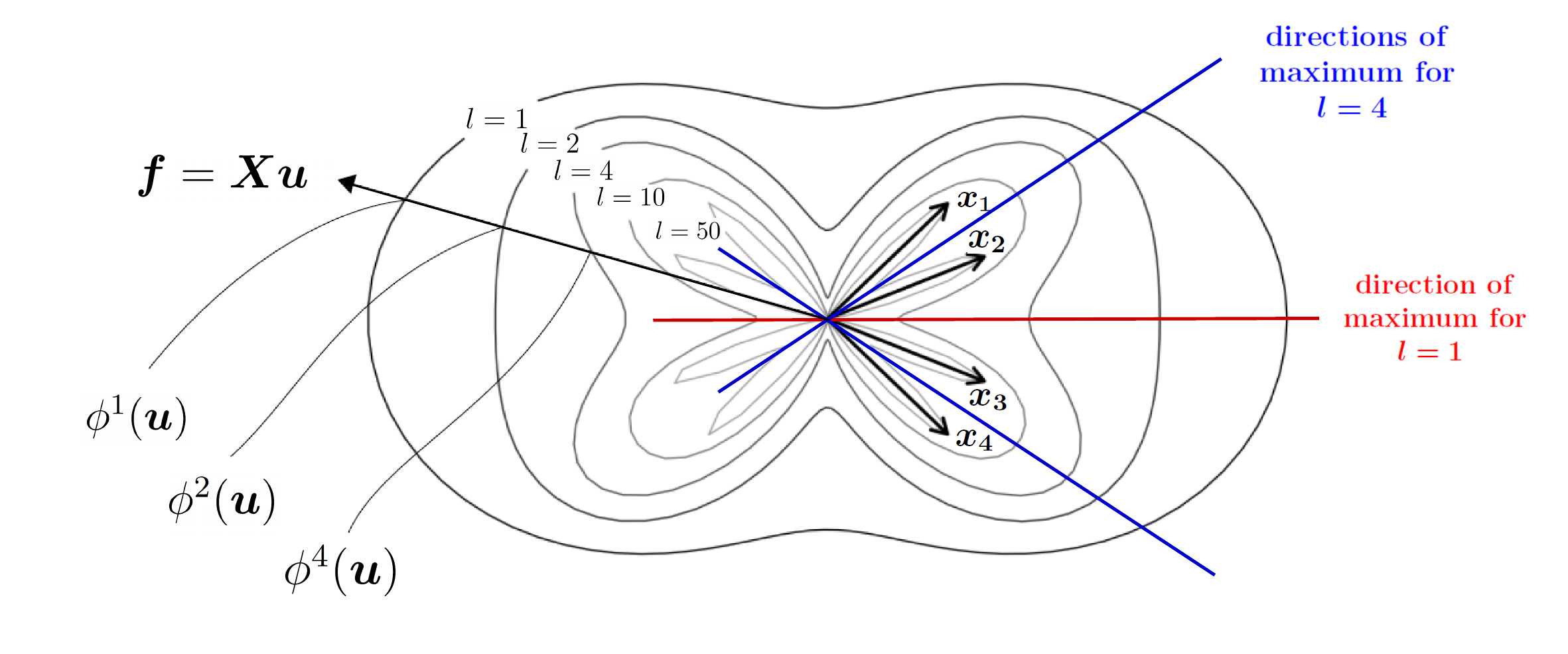}
\caption{
Polar representation of the VPI according to the value of $l$ in the elementary case of four coplanar variables, $\m{x_1}, \m{x_2}, \m{x_3}, \m{x_4}$, with 
$\omega_j=\frac{1}{4} \, \forall j \in \left\lbrace 1,2,3,4 \right\rbrace$.
\textit{
Loading--vector $\m{u}$ is identified with complex number $e^{i\theta}$, where 
$\theta \in \left[ 0, 2\pi \right)$.
Curves $z_l \left(\theta\right) \vcentcolon = \left[\phi\left(e^{i\theta}\right)\right]^l e^{i\theta}$ are graphed for $l \in \left\lbrace 1,2,4,10,50 \right\rbrace$.
The intersection of curve $z_l$ with $\m{f} = \m{Xu}$ has a radius equal to 
$\left[\phi\left(e^{i\theta}\right)\right]^l$.
The red line is the direction of maximum for $l=1$, which is in fact the first principal component. These four variables are then regarded as a unique bundle. 
By contrast, the blue lines represent the two directions of maximum for $l=4$.
The variables are then seen as two bundles containing two variables each. Finally, when $l=50$, each variable is considered a bundle in itself.
}}
\label{Chauvet::locality_bundles}
\end{figure}

\paragraph*{}
Tuning parameter $l$ allows to draw component towards more (greater $l$) or less (smaller $l$) local bundles of correlated variables, as depicted on \autoref{Chauvet::locality_bundles} in the particular instance of four coplanar variables. 
Informally, a bundle is a set of variables correlated ``enough'' to be viewed as proxies to the same latent dimension.
The notion of bundle is flexible, and parameter $l$ tunes the level of within--bundle correlation to be considered: the higher the correlation, the more local the bundle.
Overall, taking $l=1$ draws the components towards global structural directions (namely the principal components) while taking $l$ higher leads to more local ones (ultimately, the variables themselves). 
The goal is to focus on the most interpretable directions.

\paragraph*{}
Finally, let $s \in [0,1]$ be a parameter tuning the importance of the SR relative to the GoF.
SCGLR attempts a trade--off between (\ref{Chauvet::GoF_criterion}) and (\ref{Chauvet::SR_criterion_particulier}) by solving
\begin{equation*}
\underset{\m{u}:\,\m{u}^{\mbox{\tiny T}}\m{M}^{-1}\m{u}=1}{\max} \; 
\left[\phi \left( \m{u} \right)\right]^s \: 
\left[\psi_{\m{A}} \left( \m{u} \right) \right]^{1-s} 
\end{equation*}
or equivalently
\begin{equation}
\underset{\m{u}:\,\m{u}^{\mbox{\tiny T}}\m{M}^{-1}\m{u}=1}{\max} \; 
s \log \left[\phi \left( \m{u} \right)\right] + (1-s) \log \left[\psi_{\m{A}} \left( \m{u} \right)\right].
\label{Chauvet::CompromiseCriterion}  
\end{equation}
More detail can be found in \cite{bry2018}.



\section{Extension to mixed models}
\label{Chauvet::Group_SCGLR}

We now propose to extend SCGLR to mixed models. This extension will be called ``mixed--SCGLR''. 
A particular focus is placed on grouped data, for which the independence assumption of observations is no longer valid. 
The within--group dependence of each response is modelled with a random group--effect. 
Consequently, each $\m{y_k}$ is modelled with a GLMM. 
As in SCGLR, the responses are assumed to be independent conditional on the components.
\autoref{Chauvet::Onecomp} presents the single--component mixed--SCGLR method.
The underlying algorithm is given in \autoref{Chauvet::ALGORITHM}.
Considering only one component is generally not enough to explain the responses making it necessary to search for $K$ explanatory components,
with $1 \leqslant K \leqslant \text{rank} \left( {\m X} \right)$.
The way in which we extract higher rank components is explained in \autoref{Chauvet::HIGHERrank}.

\subsection{First component}
\label{Chauvet::Onecomp}

The random group--effect is assumed different across responses. This leads to $q$ random--effect vectors 
$\m{\xi_1}, \ldots, \m{\xi_q}$, which are assumed independent and normally distributed: 
\begin{equation*}
\forall k \in \left\lbrace 1, \ldots, q\right\rbrace, \;\; 
\m{\xi_k} \, 
\widesim{\text{ind.}} \, 
\mathcal{N}_N \left(\m{0}, \m{D_k} \right),
\end{equation*}
where $N$ denotes the number of groups.
In this paper, variance components models will be considered. 
We assume 
$\m{D_k} = \sigma_k^2 \,\m{I}_N$, 
where $\sigma_k^2$ is the group variance component associated with response $\m{y_k}$. 
Linear predictors involved in mixed--SCGLR are expressed as
\begin{equation}
\forall k \in \left\lbrace 1, \ldots, q\right\rbrace, \;\;
\m{\eta^{\xi}_k} = (\m{Xu})\gamma_k + \m{A\delta_k} + \m{U\xi_k},
\label{Chauvet::random_predictor}
\end{equation}
where $\m{U}$ is the known random effects' design matrix.
Predictor $\m{\eta^{\xi}_k}$ epitomises the way we capture the dependence between outcomes. Indeed, as component $\m{f} = \m{Xu}$ does not depend on $k$, it captures a structural dependence between the various 
$\m{y_k}$'s. 
By contrast, the random effect $\m{\xi_k}$ models the within--group stochastic dependence of outcomes forming response--vector 
$\m{y_k}$.

\paragraph*{}
Recall that the distribution of the data conditional on the random effects is supposed to belong to the exponential family. 
The FSA was adapted by \citet{schall1991} to the GLMM dependence structure.
The key idea is to extend Schall's algorithm to the component--based predictors in (\ref{Chauvet::random_predictor}).

\subsubsection{Linearisation step}
Let $g_k$ denote the link function for response $\m{y_k}$, 
$g'_k$ its first derivative
and $\m{\mu^{\xi}_k}$ the conditional expectation 
(i.e. $\m{\mu^{\xi}_k} \vcentcolon= 
\mathbb{E} \left( \m{y_k} \, \lvert \, \m{\xi_k} \right)$). 
The working variable associated with $y_{k,i}$ is calculated through
\begin{align*}
z_{k,i}^{\xi} 
&= g_k \left( \mu_{k,i}^{\xi} \right) + 
\left( y_{k,i} - \mu_{k,i}^{\xi} \right) 
g'_k \left( \mu_{k,i}^{\xi} \right) \\
&= 
\eta_{k,i}^{\xi} + e_{k,i}, \quad
\text{where} \quad 
e_{k,i} = \left( y_{k,i} - \mu_{k,i}^{\xi} \right) 
g'_k \left( \mu_{k,i}^{\xi} \right).
\end{align*}
In view of the conditional independence assumption,
the conditional variance matrix for $\m{z^{\xi}_k}$ is
\begin{equation*}
\mbox{Var} \left( 
\m{z^{\xi}_k} \,| \, \m{\xi_k} 
\right) 
= 
{\m{W^{\xi}_k}}^{-1} 
= 
\textbf{Diag} \left(
\left[ g'_k\left( \mu_{k,i}^{\xi} \right) \right]^2 
a_{k,i}(\phi_k) \, v_k \left( \mu_{k,i}^{\xi} \right)
\right)_{i=1, \ldots, n},
\label{Chauvet::conditional_variance}
\end{equation*}
where $a_{k,i}$ and $v_k$ are known functions, and $\phi_k$ is the dispersion parameter related to $\m{y_k}$.
At iteration $t$, the conditional linearised model for working vector $\m{z_{k}^{\xi}}$ is then defined by
\begin{equation}
{\mathcal{M}^{\xi}_k}^{[t]}: \;
\left|
\begin{aligned}
&{\m{z^{\xi}_k}}^{[t]} = \left(\m{Xu} \right)\gamma_k + \m{A\delta_k} + 
\m{U\xi_k} + \m{e}_{\m{k}}^{[t]}  \\
&\mbox{with: }
\mathbb{E} \left( \m{e}_{\m{k}}^{[t]} \, \lvert \, \m{\xi_k} \right) = 0 \mbox{ and }
\mathbb{V} \left( \m{e}_{\m{k}}^{[t]} \, \lvert \, \m{\xi_k} \right) = 
{ {\m{W^{\xi}_k}}^{-1} }^{[t]}.
\end{aligned}
\right.
\label{Chauvet::modele_lin_cond}
\end{equation} 
Besides the variance component estimation,
an alternated estimation step has to be developed
(as aforementioned in \autoref{Chauvet::New_linear_predictors})
to deal with the non--linearity of (\ref{Chauvet::modele_lin_cond}).

\subsubsection{Estimation step}

\paragraph{Calculating the component:}
Given current values of all the $\gamma_k$'s, $\m{\delta_k}$'s, 
$\m{\xi_k}$'s and $\sigma_k^2$'s, a new
component $\m{f}=\m{Xu}$ is calculated by solving a (\ref{Chauvet::CompromiseCriterion})--type program.
However, (\ref{Chauvet::GoF_criterion}) has to be adapted to conditional linearised models $\mathcal{M}^{\xi}_k$'s, involving weight matrices 
$\m{W^{\xi}_k}$'s. The appropriate goodness--of--fit measure is
\begin{equation}
\psi_{\m{A}} \left( \m{u} \right) = \sum_{k=1}^q 
\left\lVert \m{z^{\xi}_k} \right\rVert_{\m{W^{\xi}_k}}^2 
\cos_{\m{W^{\xi}_k}}^2 
\left( \m{z^{\xi}_k}, \, \spn \left\lbrace \m{Xu}, \m{A} \right\rbrace \right).
\label{Chauvet::psi_modified}
\end{equation}

\paragraph{Estimating the regression parameters and variance--components:}
Given a current value of component $\m{f}$, we apply Schall's method with the linear predictors given in (\ref{Chauvet::random_predictor}).
New values of parameters $\gamma_k$ and $\m{\delta_k}$ as well as new prediction 
$\m{\xi_k}$ are obtained by solving
the following Henderson system \citep{Henderson1975}:
\begin{equation*}
\begin{pmatrix}
\m{f}^{\mbox{\tiny T}} \hspace{0.2mm} \m{W^{\xi}_k} \, \m{f}  & 
\m{f}^{\mbox{\tiny T}} \hspace{0.2mm} \m{W^{\xi}_k} \, \m{A}  & 
\m{f}^{\mbox{\tiny T}} \hspace{0.2mm} \m{W^{\xi}_k} \, \m{U}  \smallskip \\
\m{A}^{\mbox{\tiny T}} \hspace{0.2mm} \m{W^{\xi}_k} \, \m{f}  & 
\m{A}^{\mbox{\tiny T}} \hspace{0.2mm}  \m{W^{\xi}_k} \, \m{A}  & 
\m{A}^{\mbox{\tiny T}} \hspace{0.2mm}  \m{W^{\xi}_k} \, \m{U}  \smallskip \\
\m{U}^{\mbox{\tiny T}} \hspace{0.2mm}  \m{W^{\xi}_k} \, \m{f}  & 
\m{U}^{\mbox{\tiny T}} \hspace{0.2mm}  \m{W^{\xi}_k} \, \m{A}  & 
\m{U}^{\mbox{\tiny T}} \hspace{0.2mm}  \m{W^{\xi}_k} \, \m{U} + \m{D}_{\m{k}}^{-1}
\end{pmatrix}
\begin{pmatrix}
\gamma_k \smallskip \\ 
\m{\delta_k} \smallskip \\ 
\m{\xi_k}
\end{pmatrix} 
= 
\begin{pmatrix}
\m{f}^{\mbox{\tiny T}} \hspace{0.2mm}  \m{W^{\xi}_k} \, \m{z^{\xi}_k}  \smallskip \\
\m{A}^{\mbox{\tiny T}} \hspace{0.2mm}  \m{W^{\xi}_k} \, \m{z^{\xi}_k}  \smallskip \\
\m{U}^{\mbox{\tiny T}} \hspace{0.2mm}  \m{W^{\xi}_k} \, \m{z^{\xi}_k}
\end{pmatrix}.
\end{equation*}
Finally, as mentioned by \citet{schall1991}, given prediction 
$\widehat{\m{\xi_k}}$ for $\m{\xi_k}$, the update of the ML estimation of variance component $\sigma_k^2$ is
\begin{equation*}
\sigma_{k}^2 \: \longleftarrow \:
\dfrac{\widehat{\m{\xi_k}}^{\mbox{\tiny T}} \; \widehat{\m{\xi_k}}}
{N - \frac{1}{\sigma_{k}^2} \, \text{Trace} 
\left[ \left( 
\m{U}^{\mbox{\tiny T}} \hspace{0.2mm} \m{W^{\xi}_k} \, \m{U} + \m{D}_{\m{k}}^{-1} \right)^{-1}  \right]}.
\end{equation*}

\subsection{The algorithm}
\label{Chauvet::ALGORITHM}

The conditional linearised models considered at iteration $t$ are given by (\ref{Chauvet::modele_lin_cond}).
\autoref{Chauvet::PSEUDOCODE} describes the $(t+1)$--th iteration of the single--component mixed--SCGLR. 
It is repeated until stability of parameters is reached.

\vspace{-3mm}
\begin{center}
\fbox{\begin{minipage}{1\linewidth}
\small
\begin{algorithm}[H]
\begin{description}
\vspace{1mm}
\item[Step 1:] \textbf{Computing the component.} Set 
\vspace{-3mm}
\begin{flalign*}
\hspace{-2em}
\m{u}^{[t+1]} 
= 
\underset{\m{u}:\,\m{u}^{\mbox{\tiny T}}\m{M}^{-1} \m{u}=1}{\mbox{arg max}} \; 
\left[ \phi\left( \m{u} \right)\right]^s \: 
\left[ \psi_{\m{A}}^{[t]} \left( \m{u} \right)\right]^{1-s}, \;
\text{where } \psi_{\m{A}} \left( \m{u} \right) \, \text{is given by }   (\ref{Chauvet::psi_modified}) 
\text{ and } \phi\left( \m{u} \right) \text{ by } (\ref{Chauvet::SR_criterion_particulier}) &&
\end{flalign*}
\vspace{-12mm}
\begin{flalign*}
\hspace{-2em}
\m{f}^{[t+1]} 
= \m{X} \m{u}^{[t+1]} &&
\end{flalign*}
\vspace{-12mm}
\item[Step 2:] \textbf{Henderson systems.}
For each $k \in \left\lbrace 1, \ldots, q \right\rbrace $, solve the system
\begin{equation*}
\hspace{-5em}
\begin{pmatrix}
{\m{f}^{[t+1]}}^{\mbox{\tiny T}} {\m{W^{\xi}_k}}^{[t]} {\m{f}^{[t+1]}} & 
{\m{f}^{[t+1]}}^{\mbox{\tiny T}} {\m{W^{\xi}_k}}^{[t]} \m{A} & 
{\m{f}^{[t+1]}}^{\mbox{\tiny T}} {\m{W^{\xi}_k}}^{[t]} \m{U} 
\smallskip \\
\m{A}^{\mbox{\tiny T}} {\m{W^{\xi}_k}}^{[t]} {\m{f}^{[t]}}    & 
\m{A}^{\mbox{\tiny T}} {\m{W^{\xi}_k}}^{[t]} \m{A}    & 
\m{A}^{\mbox{\tiny T}} {\m{W^{\xi}_k}}^{[t]} \m{U}   \smallskip \\
\m{U}^{\mbox{\tiny T}} {\m{W^{\xi}_k}}^{[t]} {\m{f}^{[t]}}    & 
\m{U}^{\mbox{\tiny T}} {\m{W^{\xi}_k}}^{[t]} \m{A}    & 
\m{U}^{\mbox{\tiny T}} {\m{W^{\xi}_k}}^{[t]} \m{U} + 
{\m{D}_{\m{k}}^{[t]}}^{-1}
\end{pmatrix}
\begin{pmatrix}
\gamma_k \smallskip \\ 
\m{\delta_k} \smallskip \\ 
\m{\xi_k}
\end{pmatrix} 
= 
\begin{pmatrix}
{\m{f}^{[t+1]}}^{\mbox{\tiny T}} {\m{W^{\xi}_k}}^{[t]} \, 
\m{z^{\xi}_k}^{[t]} \smallskip \\
\m{A}^{\mbox{\tiny T}} {\m{W^{\xi}_k}}^{[t]} \, 
\m{z^{\xi}_k}^{[t]}   
\smallskip \\
\m{U}^{\mbox{\tiny T}} {\m{W^{\xi}_k}}^{[t]} \, 
\m{z^{\xi}_k}^{[t]}
\end{pmatrix}
\end{equation*}
Call $\gamma_k^{[t+1]}$, $\m{\delta}_{\m{k}}^{[t+1]}$ and 
$\m{\xi}_{\m{k}}^{[t+1]}$ the solutions.
\item[Step 3:] \textbf{Updating variance--component estimates.}
For each $k \in \left\lbrace 1, \ldots, q \right\rbrace $, compute
\begin{equation*}
\hspace{-5em}
{\sigma_{k}^2}^{[t+1]} =
\dfrac{ {\m{\xi}_{\m{k}}^{[t+1]}}^{\mbox{\tiny T}} \; \m{\xi}_{\m{k}}^{[t+1]} }
{N - \frac{1}{{\sigma_{k}^2}^{[t]}} \,
\text{Trace} \left[ \left( \m{U}^{\mbox{\tiny T}} \, 
{\m{W^{\xi}_k}}^{[t]} \, \m{U} + {\m{D}_{\m{k}}^{[t]}}^{-1} \right)^{-1} \right]}
\quad \text{and} \quad
\m{D}_{\m{k}}^{[t+1]} = {\sigma_{k}^2}^{[t+1]} \m{I}_{N} 
\end{equation*} 
\item[Step 4:] \textbf{Updating working variables and weighting matrices.} \\
For each $k \in \left\lbrace 1, \ldots, q \right\rbrace $, compute
\begin{align*}
{\m{\eta^{\xi}_k}}^{[t+1]} 
&=
\m{f}^{[t+1]} \gamma_k^{[t+1]} + 
\m{A} \m{\delta}_{\m{k}}^{[t+1]} + 
\m{U} \m{\xi}_{\m{k}}^{[t+1]} 
\\
{\mu_{k,i}^{\xi}}^{[t+1]} 
&=
g_{k}^{-1} \left( 
{\eta_{k,i}^{\xi}}^{[t+1]} 
\right), 
\; i=1, \ldots, n \\
{z_{k,i}^{\xi}}^{[t+1]}
&= 
{\eta_{k,i}^{\xi}}^{[t+1]} + 
\left( y_{i}^{k} - {\mu_{k,i}^{\xi}}^{[t+1]} \right) 
g'_{k} \left( {\mu_{k,i}^{\xi}}^{[t+1]} \right), 
\; i=1, \ldots, n \\
{\m{W^{\xi}_k}}^{[t+1]} &= 
\textbf{Diag} \left( 
\left\lbrace 
\left[ g'_{k} \left( {\mu_{k,i}^{\xi}}^{[t+1]} \right) \right]^{2} 
a_{k,i}(\phi_k) \, v_k\left( {\mu_{k,i}^{\xi}}^{[t+1]} \right)
\right\rbrace^{-1}
\right)_{i=1, \ldots, n} 
\end{align*}
\vspace{-12mm}
\item[Incrementing:] $t \longleftarrow t+1$
\end{description}
\caption{Iteration of the single--component mixed--SCGLR}
\label{Chauvet::PSEUDOCODE}
\end{algorithm}
\end{minipage}}
\end{center}

\subsection{Extracting higher rank components}
\label{Chauvet::HIGHERrank}

Let ${\m{F_h}} = \big[\, \m{f_1} \mid \ldots \mid \m{f_h} \,\big]$ be the matrix of the first $h$ components, where $h<K$.
An extra component $\m{f_{h+1}}$ must best complement the existing ones plus 
$\m{A}$, i.e. 
$\m{A_h} = \big[\, \m{F_h} \mid \m{A} \,\big]$.
So $\m{f_{h+1}}$ must be calculated using $\m{A_h}$ as additional explanatory variables.
Moreover, we must impose that $\m{f_{h+1}}$ be orthogonal to $\m{F_h}$, i.e.
$\m{F}_{\m{h}}^{\mbox{\tiny T}} \m{P} \m{f_{h+1}} = \m{0}$.
Component $\m{f_{h+1}} = \m{X} \m{u_{h+1}}$ is thus obtained by solving
\begin{equation}
\left\lbrace
\begin{aligned}
&\text{max} \quad  
s \log \left[\phi\left(\m{u}\right)\right] + 
(1-s) \log \left[\psi_{\m{A_h}}\left(\m{u}\right)\right]   \\
&\text{subject to:}  \quad \m{u}^{\mbox{\tiny T}} \m{M}^{-1} \m{u} = 1 \; \text{and} \; \m{D}_{\m{h}}^{\mbox{\tiny T}} \m{u} = \m{0},
\end{aligned}
\right.
\label{Chauvet::program_plusieurs_comp}
\end{equation}
where 
$\psi_{\m{A_h}}\left(\m{u}\right) = 
\displaystyle{\sum_{k=1}^q} 
\left\lVert \m{z^{\xi}_k} \right\rVert_{\m{W^{\xi}_k}}^2 
\cos_{\m{W^{\xi}_k}}^2 \left( 
\m{z^{\xi}_k}, \, \spn \left\lbrace \m{Xu}, \m{A_h} \right\rbrace 
\right)$
and $\m{D_h} = \m{X}^{\mbox{\tiny T}} \m{P F_h}$.

\paragraph*{}
In the online Supplementary Material, we give a simple tool to maximise, at least locally, any criterion on the unit sphere: the Projected Iterated Normed Gradient (PING) algorithm. 
In particular, PING solves (\ref{Chauvet::program_plusieurs_comp})--type programs, 
which give all components of rank $h>1$.
The rank--one component is computed using the same program with
$\m{A_0} = \m{A}$ and $\m{D_0}=\m{0}$.



\section{Comparative results on simulated data}
\label{Chauvet::simulatedDATA}

Five simulation studies have been implemented to assess our method. The first one 
(discussed in \autorefs{Chauvet::datageneration} -- \ref{Chauvet::model_interpret_simulated_data})
focuses on LMMs. It compares the performances of mixed--SCGLR, LMM--ridge \citep{eliot2011} and GLMM--LASSO \citep{groll2014, schelldorfer2014}.
The second simulation (\autoref{Chauvet::simulated_data_binary_poisson}) extends the first one to binary and Poisson outcomes.
All simulation studies have been performed using R \citep{logicielR}.
To compute LASSO regressions, we have used the R package \texttt{glmmLasso} \citep{packageGLMMLASSO}.
The extension of SCGLR to mixed models is available at 
\url{https://github.com/SCnext/mixedSCGLR}.
Three additional simulations are presented in the online Supplementary Material. 
The first one reproduces the simulation scheme of \autoref{Chauvet::simulated_data_binary_poisson} with binomial and Poisson outcomes. 
The second one assesses the performance of mixed--SCGLR on a different bundle structure and presents results concerning variance component estimates.
The third one deals with high dimensional data.

\subsection{Data generation}
\label{Chauvet::datageneration}

To generate grouped data,
we consider $N=10$ groups and $R=10$ observations per group 
(i.e. a total of $n=100$ observations). 
The random effects' design matrix is then
$\m{U} = \m{I}_N \otimes \mathbf{1}_R$. 
Explanatory variables $\m{X}$ consist of
three independent bundles: 
$\m{X_0}$ (15 variables), $\m{X_1}$ (10 variables) 
and $\m{X_2}$ (5 variables).
Each explanatory variable is normally simulated with mean $0$ and variance $1$. 
Parameter $\tau \in \left\lbrace 0.1, 0.3, 0.5, 0.7, 0.9 \right\rbrace$
tunes the level of redundancy within each bundle: 
the correlation matrix of bundle $\m{X_j}$ is
\begin{equation*}
\text{cor} \left(\m{X_j} \right) = \tau \mathbf{1}_{p_j} 
\mathbf{1}_{p_j}^{\mbox{\tiny T}} + (1-\tau) \m{I}_{p_j},
\end{equation*}
where $p_j$ is the number of variables in $\m{X_j}$.
In order to enable comparison with LASSO and ridge and to focus on regularisation, 
our simulations do not involve additional explanatory variables ($\m{A}=\m{0}$).
Two random responses $\m{Y} = \big[\, \m{y_1} \mid \m{y_2} \,\big]$ are generated as
\begin{equation}
\left\lbrace
\begin{aligned}
\m{y_1} &= \m{X\beta_1} + \m{U\xi_1} + \m{\varepsilon_1} \\
\m{y_2} &= \m{X\beta_2} + \m{U\xi_2} + \m{\varepsilon_2},
\end{aligned}
\right.
\label{Chauvet:modele_simulateddata}
\end{equation}
such that $\m{y_1}$ is predicted only by bundle $\m{X_1}$, $\m{y_2}$ only by bundle $\m{X_2}$, 
and bundle $\m{X_0}$ plays no explanatory role.
Our choice for the fixed--effect parameters is
\begin{align*}
\m{\beta_1} &= (\; 
\underbrace{0, \ldots\ldots\ldots\ldots\ldots\ldots, 0}_{15 \: \text{times}},
\underbrace{0.3, .. , 0.3}_{3 \: \text{times}}, 
\underbrace{0.4, \ldots, 0.4}_{4 \: \text{times}}, 
\underbrace{0.5, .. , 0.5}_{3 \: \text{times}}, 
\underbrace{0, \ldots, 0}_{5 \: \text{times}}
\;)^{\mbox{\tiny T}}, \\
\m{\beta_2} &=  (\; 
\underbrace{0, 
\ldots\ldots\ldots\ldots\ldots\ldots
\ldots\ldots\ldots\ldots\ldots\ldots\ldots, 0}_{25 \: \text{times}}, 
0.3,0.3,0.4,0.5,0.5
\;)^{\mbox{\tiny T}}.
\end{align*}
Finally, for each $k \in \left\lbrace 1,2 \right\rbrace$, 
random effect and noise vectors are simulated respectively from
\begin{equation*}
\m{\xi_k} \widesim{} \mathcal{N}_N
\left(\m{0}, \; \sigma_k^2 \, \m{I}_N \right) \; \text{and} \;\,
\m{\varepsilon_k} \widesim{} \mathcal{N}_n
\left(\m{0}, \; \omega_k^2 \, \m{I}_n \right),
\end{equation*}
where
$\sigma_k^2 = \omega_k^2 = 1$. 
For each value of $\tau$, $B=100$ samples are generated according to Model (\ref{Chauvet:modele_simulateddata}). 

\subsection{Parameter calibration}
\label{Chauvet::section_parameter_calibration}

In order to compare mixed--SCGLR with the ridge and LASSO regressions, we recall the regularisation parameters required by each method.
For both LMM--ridge and GLMM--LASSO methods, a unique shrinkage parameter has to be calibrated: 
$\lambda_{\text{ridge}}$ and $\lambda_{\text{LASSO}}$ respectively.
For mixed--SCGLR, three tuning parameters need to be calibrated: 
the number of components $K$ and the trade--off parameter $s$, 
which are both regularisation parameters,
and the bundle--locality parameter $l$. 
For greater clarity, the simulation focuses on the behaviour of $K$ and $s$. 
As recommended by \citet{bry2013}, we set $l=4$.
In case--studies, $l$ has to be tuned to maximise the interpretability of components.

\paragraph*{}
For both mixed--SCGLR and GLMM--LASSO, 
optimal regularisation parameters
are obtained through a $5$--fold cross--validation, withdrawing $2$ observations from each group every time. 
This could be termed ``leave--two--observations--out per group.''
The data are thus divided into five parts $\mathcal{P}_{1}, \ldots, \mathcal{P}_{5}$,
each $\mathcal{P}_{j}$ containing $20$ observations, $2$ for each of the $10$ groups. 
Let $y_{k,i}^{(b)}$ be the $i$--th observation of the $k$--th response vector in the $b$--th sample. Let also 
$\widehat{y^{(b)}_{k,i (-j)}}$ be the fit for $y_{k,i}^{(b)}$ with part 
$\mathcal{P}_j$ removed.
The cross--validation error in the $b$--th sample, $E^{(b)}$, is defined as
\begin{equation}
E^{(b)} = 
\dfrac{1}{2} \sum_{k=1}^2 E^{(b)}_k, \label{Chauvet::crossvalerror}
\end{equation}
where
\begin{equation*}
E^{(b)}_k = 
\dfrac{1}{5} \sum_{j=1}^5 
\sqrt{
\dfrac{1}{20} \sum_{i \in \mathcal{P}_j} 
\left( y_{k,i}^{(b)} - \widehat{y^{(b)}_{k,i (-j)}} \right)^2
}.
\end{equation*}
In the $b$--th sample, the optimal number of components ${K^{\star}}^{(b)}$, 
the trade--off parameter ${s^{\star}}^{(b)}$, 
and the shrinkage parameter 
$\lambda_{\text{LASSO}}^{\star^{\mathlarger{(b)}}}$ are selected to minimise the cross--validation error (\ref{Chauvet::crossvalerror}).
We then define
\begin{equation*}
s^{\star} = 
\dfrac{1}{B} \sum_{b=1}^B {s^{\star}}^{(b)}, \;\,
K^{\star} = 
\mbox{mode} \left( \left\lbrace 
{K^{\star}}^{(1)}, \ldots, {K^{\star}}^{(B)} \right\rbrace \right) 
\;\, \text{and} \;\,
\lambda_{\text{LASSO}}^{\star} = 
\dfrac{1}{B} \sum_{b=1}^B \lambda_{\text{LASSO}}^{\star^{\mathlarger{(b)}}}.
\end{equation*}
By contrast, \citet{eliot2011} suggest to calibrate the ridge parameter at each step of their EM implementation, using the generalised cross--validation. We thus define
\begin{equation*}
\lambda_{\text{ridge}}^{\star} = 
\dfrac{1}{B} \sum_{b=1}^B \lambda_{\text{ridge}}^{\star^{\mathlarger{(b)}}},
\end{equation*}
where $\lambda_{\text{ridge}}^{\star^{\mathlarger{(b)}}}$ denotes the average of the ridge parameter values obtained over all the iterations of the EM algorithm in the $b$--th sample.

\paragraph*{}
\autoref{Chauvet::res_shrinkage_par} summarises the optimal regularisation parameters selected through cross-- \linebreak validation.
In both ridge and LASSO, the shrinkage parameter value increases with the level of redundancy $\tau$. 
Whereas for mixed-SCGLR, when $\tau$ increases, $K^{\star}$ decreases towards the true number of predictive variable--bundles: the greater the value of $\tau$, the better mixed--SCGLR focuses on the  structures in $\m{X}$ that contribute to model $\m{Y}$.
Moreover, when $\tau$ increases, the trade--off parameter $s^{\star}$ increases, meaning that regularisation requires a greater importance of the structural relevance relative to the goodness--of--fit.

\begin{table}[!ht]
\caption{Optimal regularisation parameter values obtained through cross--validation over $100$ simulations.} 
\label{Chauvet::res_shrinkage_par}
\vspace{3mm}
\begin{tabular}{l p{3.5cm} p{3.25cm} p{3cm} p{2.75cm} }
& \centering \textbf{GLMM--LASSO} & \centering \textbf{LMM--ridge} &  \multicolumn{2}{c}{\centering \textbf{mixed--SCGLR}} \\
& \centering  shrinkage parameter 
$\;\lambda_{\text{LASSO}}^{\star}\;$ & 
\centering  shrinkage parameter 
$\lambda_{\text{ridge}}^{\star}$ & 
\centering  number of components $K^{\star}$ & 
\centering  trade--off parameter $s^{\star}$
\tabularnewline 
\hline \noalign{\smallskip}
$\tau=0.1$ & \centering  65 & \centering 24 & \centering 15 & \centering 0.50 \tabularnewline  
$\tau=0.3$ & \centering  92 & \centering 54 & \centering 5 & \centering 0.58 \tabularnewline 
$\tau=0.5$ & \centering 124 & \centering 73 & \centering 3 & \centering 0.70 \tabularnewline 
$\tau=0.7$ & \centering 163 & \centering 78 & \centering 3 & \centering 0.73 \tabularnewline
$\tau=0.9$ & \centering 175 & \centering 85 & \centering 2 & \centering 0.80
\tabularnewline 
\hline
\end{tabular}
\end{table}

\subsection{Comparison of the estimate accuracies}
\label{Chauvet::gaussian_comparison_estimates}

Once tuning parameters are obtained,
we focus on the fixed--effect estimates' accuracy.
Since the response--vectors $\m{y_1}$ and $\m{y_2}$ are normally distributed and have comparable orders of magnitude, the fixed--effect relative errors are on the same scale.
Then we consider a risk--averse comparison criterion called ``Mean Upper Relative Squared Error'' (MURSE) defined as
\begin{equation*}
\text{MURSE} \left( \m{\beta_1}, \m{\beta_2} \right) = 
\dfrac{1}{B} \sum_{b=1}^B \max \left\lbrace
\dfrac{\left\lVert \widehat{\m{\beta}}_{\m{1}}^{(b)} - \m{\beta_1} \right\rVert^2}{\left\lVert \m{\beta_1} \right\rVert^2}, 
\dfrac{\left\lVert \widehat{\m{\beta}}_{\m{2}}^{(b)} - \m{\beta_2} \right\rVert^2}{\left\lVert \m{\beta_2} \right\rVert^2}
\right\rbrace,
\end{equation*}
where $\widehat{\m{\beta}}_{\m{k}}^{(b)}$ is the estimate of 
$\m{\beta_k}$ associated with sample $b$.
The MURSE values for mixed--SCGLR, LMM--ridge and GLMM--LASSO are presented in \autoref{Chauvet::res_estim_study}.
The LMM results obtained without regularisation are also presented. They were computed using the R package \texttt{lme4} \citep{packagelme4}.
In the latter case, relative errors increase dramatically with $\tau$.
Those of ridge and LASSO increase less drastically (but increase anyway) because these methods suffer from the high correlations among the explanatory variables.
Except for $\tau=0.1$, mixed--SCGLR provides the most accurate fixed effect estimates.
Indeed, if there are no real bundles in $\m{X}$ ($\tau \simeq 0$), searching for structures in $\m{X}$ may lead mixed--SCGLR to be slightly less accurate. 
Conversely, mixed--SCGLR takes advantage of the high correlations among the explanatory variables: the stronger the structures (high $\tau$), the more efficient the method.

\begin{table}[!ht]
\caption{Mean Upper Relative Squared Error (MURSE) values associated with the optimal parameter values. } 
\label{Chauvet::res_estim_study}
\vspace{3mm}
\begin{tabular}{lcccc} 
 & \centering \textbf{LMM} & \multirow{2}{*}{\textbf{GLMM--LASSO}}
 & \multirow{2}{*}{\textbf{LMM--ridge}}
 & \multirow{2}{*}{\textbf{mixed--SCGLR}}
 \tabularnewline 
 & (no regularisation) & & & 
 \tabularnewline
\hline \noalign{\smallskip}
$\tau=0.1$ & \centering \tablenum{0.12} & \centering \tablenum{0.05} & 
\centering \tablenum{0.08} & \centering \tablenum{0.12}
\tabularnewline  
$\tau=0.3$ & \centering \tablenum{0.33} & \centering \tablenum{0.12} & 
\centering \tablenum{0.13} & \centering \tablenum{0.10}
\tabularnewline
$\tau=0.5$ & \centering \tablenum{0.61} & \centering \tablenum{0.20} & 
\centering \tablenum{0.16} & \centering \tablenum{0.07}
\tabularnewline
$\tau=0.7$ & \centering \tablenum{1.32} & \centering \tablenum{0.25} & 
\centering \tablenum{0.20} & \centering \tablenum{0.06} 
\tabularnewline  
$\tau=0.9$ & \centering \tablenum{4.62} & \centering \tablenum{0.26} & 
\centering \tablenum{0.31} & \centering \tablenum{0.05}
\tabularnewline
\hline
\end{tabular}
\end{table}

\subsection{Model interpretation}
\label{Chauvet::model_interpret_simulated_data}

\begin{figure}[!ht]
\centering
\includegraphics[width=.45\linewidth, trim={2.5cm 0 2.25cm 0},clip]
{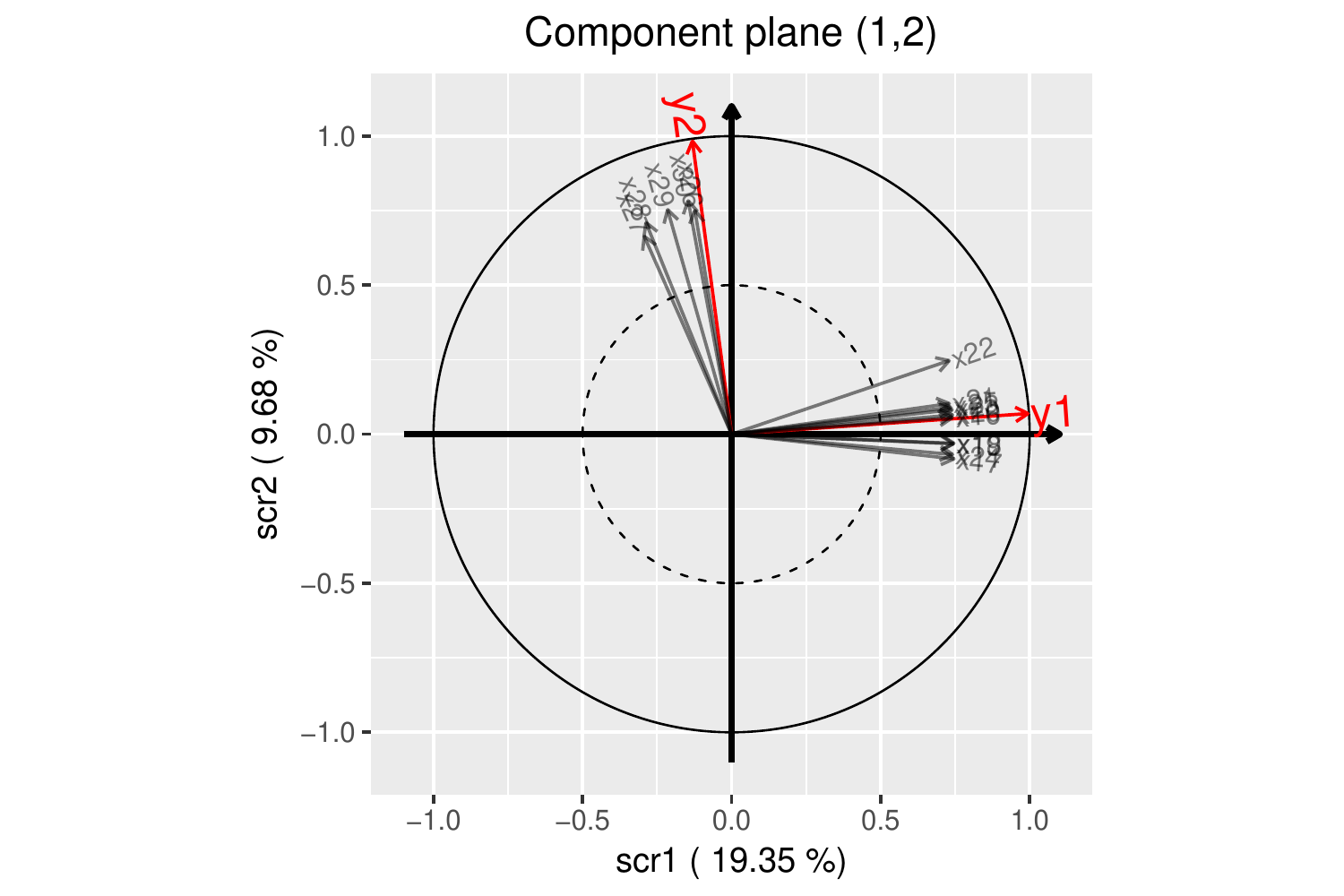}%
\hspace{5mm}
\includegraphics[width=.45\linewidth, trim={2.5cm 0 2.25cm 0},clip]
{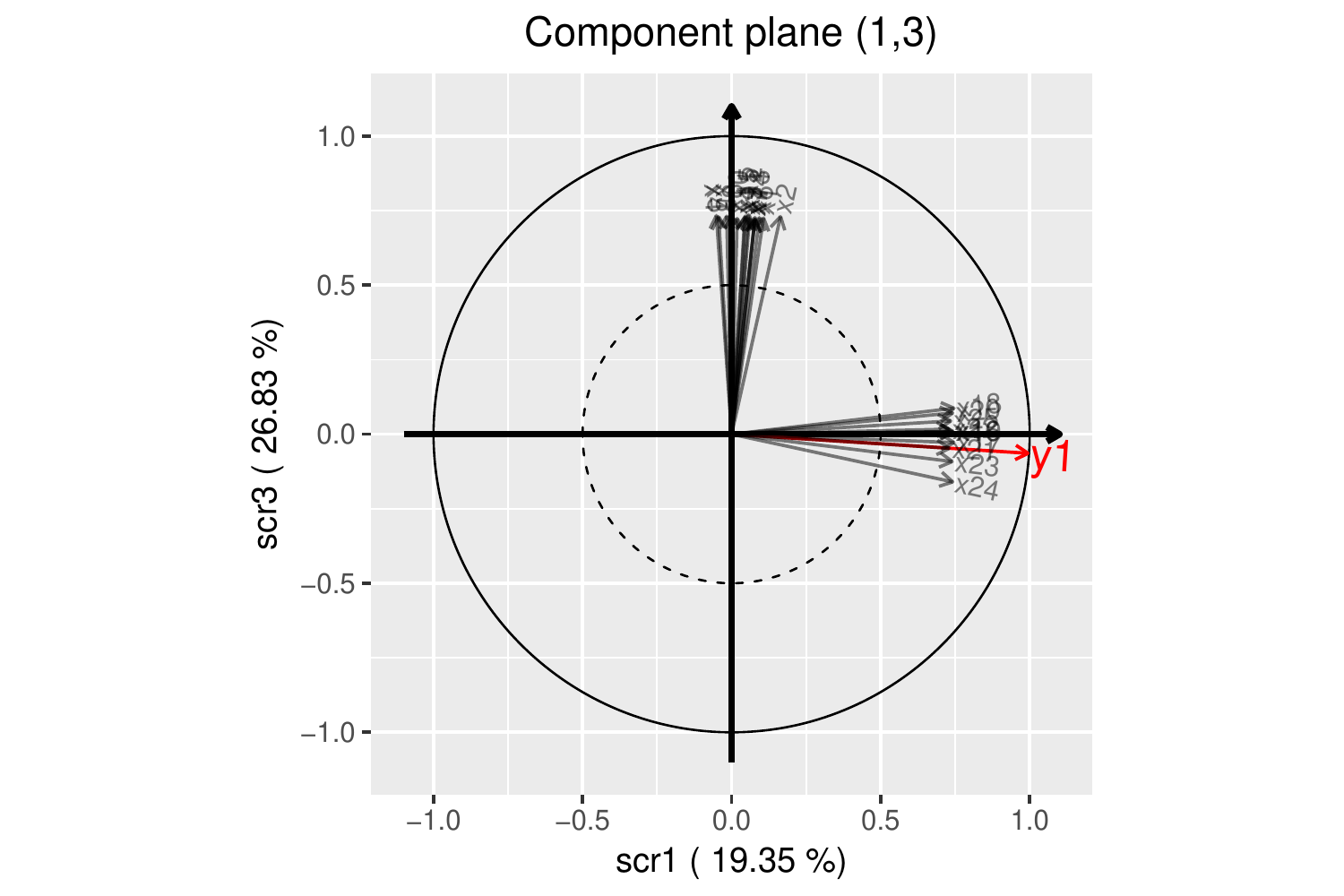}%
\caption{Component planes $(1,2)$ and $(1,3)$ given by mixed--SCGLR on simulated data. 
\textit{The black arrows represent the explanatory variables. The red ones represent the projection of the $\m{X}$--part of the linear predictors associated with $\m{y_1}$ and $\m{y_2}$. The percentage of inertia captured by each component is given in parentheses.}}
\label{Chauvet::componentplanes}
\end{figure}

This section aims at highlighting the power of mixed--SCGLR for model interpretation. 
\autoref{Chauvet::componentplanes} presents an example of the first component planes obtained for 
$\tau=0.5$, with associated optimal parameter values
$s^{\star}=0.7$ and 
$K^{\star}=3$. We still impose $l=4$.
The first two components obtained are the ones which explain the responses.
It clearly appears that $\m{y_1}$ is explained by bundle $\m{X_1}$ and 
$\m{y_2}$ by $\m{X_2}$. 
Interestingly, although bundle $\m{X_0}$ is the one with maximum inertia 
($26.83\%$), it appears only along the third component, for having no explanatory part.

\subsection{Additional simulations involving non--Gaussian outcomes}
\label{Chauvet::simulated_data_binary_poisson}

This section aims at assessing our method in the case of Bernoulli ($\mathcal{B}$) and Poisson ($\mathcal{P}$) distributions of responses. We still consider $N=10$ groups and $R=10$ observations per group. 
We keep design matrices $\m{X}$ and $\m{U}$ defined in \autoref{Chauvet::datageneration}, as well as the values of 
$\m{\beta_1}$, $\m{\beta_2}$, $\sigma_1^2$ and $\sigma_2^2$.
The group variance components are given by 
$\varsigma_1^2 = 0.1 \sigma_1^2$ and $\varsigma_2^2 = \sigma_2^2$ so that
for each $k \in \left\lbrace 1,2 \right\rbrace$, 
$\widetilde{\m{\xi_k}} \widesim{} \mathcal{N}_N
\left(\m{0}, \; \varsigma_k^2 \, \m{I}_N \right)$.
Then given $\widetilde{\m{\xi_1}}$ and $\widetilde{\m{\xi_2}}$, 
we simulate $\m{Y} = \big[\, \m{y_1} \mid \m{y_2} \,\big]$ as
\begin{equation}
\left\lbrace
\begin{aligned}
&\m{y_1} \widesim{} \mathcal{B} \left( 
\m{p} = \text{logit}^{-1} \left[ 
\m{X \theta_1} + \m{U \widetilde{\xi_1}} 
\right] \right) \\
&\m{y_2} \widesim{} \mathcal{P} \left( 
\m{\lambda} = \exp \left[
\m{X \theta_2} + \m{U \widetilde{\xi_2}} \right]
\right),
\end{aligned}
\right.
\label{Chauvet:modele_binary_poisson}
\end{equation}
where $\m{\theta_1} = 0.1 \m{\beta_1}$ and $\m{\theta_2} = \m{\beta_2}$.
Again, for each value of $\tau$, $B=100$ samples are generated according to Model (\ref{Chauvet:modele_binary_poisson}). 
As in \autoref{Chauvet::section_parameter_calibration}, tuning parameters are calibrated so as to minimise the cross--validation error (\ref{Chauvet::crossvalerror}). 
However, since $\m{y_1}$ and $\m{y_2}$ do not have the same range of values, the prediction errors have to be standardised. The cross--validation error for response $\m{y_k}$ in the $b$--th sample is now given by
\begin{equation*}
E^{(b)}_k = 
\dfrac{1}{5} \sum_{j=1}^5 
\sqrt{
\dfrac{1}{20} \sum_{i \in \mathcal{P}_j}
\dfrac{
\left( y_{k,i}^{(b)} - \widehat{y^{(b)}_{k,i (-j)}} \right)^2}
{\widehat{\text{var}}\left( \widehat{y^{(b)}_{k,i (-j)}} \right)}
}.
\end{equation*}

Unlike in \autoref{Chauvet::gaussian_comparison_estimates}, the response--vectors do not come from the same distribution and have different orders of magnitude. The fixed--effect relative errors are thus not comparable.
To compare mixed--SCGLR with GLMM--LASSO and classical GLMM (without regularisation), we thus use the Mean Relative Squared Error (MRSE) defined as
$$\text{MRSE} \left( \m{\theta_k} \right) = 
\dfrac{1}{B} \sum_{b=1}^B \dfrac{\left\lVert \widehat{\m{\theta}}_{\m{k}}^{(b)} - \m{\theta_k} \right\rVert^2}{\left\lVert \m{\theta_k} \right\rVert^2}, \; k \in \left\lbrace 1,2 \right\rbrace,$$
where $\widehat{\m{\theta}}_{\m{k}}^{(b)}$ is the estimate of 
$\m{\theta_k}$ from the $b$--th sample.
MRSE values for the GLMM, mixed--SCGLR and GLMM--LASSO are presented in \autoref{Chauvet::table_simu_binary_poisson}. 
For all methods, estimating a Bernoulli model is obviously a more challenging task than estimating a Poisson model. 
Regardless of the level of redundancy $\tau$,
both mixed--SCGLR and GLMM--LASSO outperform classical GLMM estimation. 
Compared with the Gaussian case (\autoref{Chauvet::gaussian_comparison_estimates}), 
the results deteriorate but (overall) the same behaviours are observed. 
\begin{itemize}[topsep=0cm, itemsep=0cm, parsep=0cm]
\item[$\blacktriangleright$]
For $\tau = 0.1$, fixed--effect estimates provided by mixed--SCGLR 
are less accurate than those provided by GLMM--LASSO.
In this case, GLMM--LASSO has indeed a double advantage. 
First, many $\theta_{k,j}$'s are true zeros. 
Unlike mixed--SCGLR, GLMM--LASSO 
often shrinks their estimates to exactly zero. 
Second, since the level of redundancy is low, GLMM--LASSO also provides accurate coefficient estimates of active variables.
\item[$\blacktriangleright$]
By contrast, for $\tau \geqslant 0.3$, mixed--SCGLR 
takes advantage of redundancies within the explanatory variables.
Thus, mixed-SCGLR outperforms GLMM--LASSO in this case, despite the sparse structure of the $\m{\theta_k}$'s.
\end{itemize}

\begin{table}[!ht]
\caption{Mean Relative Squared Error (MRSE) values obtained with Bernoulli and Poisson responses. } 
\label{Chauvet::table_simu_binary_poisson}
\vspace{3mm}
\begin{tabular}{p{2cm} p{1.8cm} p{1.75cm} p{1.75cm} p{1.75cm} p{1.75cm} p{1.75cm} }
& 
\multicolumn{2}{c}{\centering \textbf{GLMM}} & 
\multicolumn{2}{c}{\multirow{2}{*}{\textbf{GLMM--LASSO}}} & 
\multicolumn{2}{c}{\multirow{2}{*}{\textbf{mixed--SCGLR}}} \\
 & \multicolumn{2}{c}{\centering (no regularisation)} & & & & \\
& \centering  Bernoulli & 
\centering  Poisson & 
\centering  Bernoulli & 
\centering  Poisson & 
\centering  Bernoulli & 
\centering  Poisson
\tabularnewline 
\hline \noalign{\smallskip}
$\tau=0.1$ & 
\centering  \tablenum{316.48} &  \centering  \tablenum{0.54} & 
\centering  \tablenum{8.61}   &  \centering  \tablenum{0.30} & 
\centering  \tablenum{14.71}  &  \centering  \tablenum{0.46}   \tabularnewline  
$\tau=0.3$ & 
\centering  \tablenum{398.78} &  \centering  \tablenum{0.64} & 
\centering  \tablenum{9.23}   &  \centering  \tablenum{0.36} & 
\centering  \tablenum{7.21}   &  \centering  \tablenum{0.21}   \tabularnewline 
$\tau=0.5$ & 
\centering  \tablenum{576.68} &  \centering  \tablenum{0.87} & 
\centering  \tablenum{14.48}  &  \centering  \tablenum{0.44} & 
\centering  \tablenum{2.01}   &  \centering  \tablenum{0.09}   \tabularnewline 
$\tau=0.7$ & 
\centering  \tablenum{886.04} &  \centering  \tablenum{1.28} & 
\centering  \tablenum{17.37}  &  \centering  \tablenum{0.47} & 
\centering  \tablenum{1.50}   &  \centering  \tablenum{0.07}   \tabularnewline 
$\tau=0.9$ & 
\centering  \tablenum{2840.10} &  \centering  \tablenum{3.72} & 
\centering  \tablenum{17.24}   &  \centering  \tablenum{0.59} & 
\centering  \tablenum{1.31}    &  \centering  \tablenum{0.05}   \tabularnewline 
\hline
\end{tabular}
\end{table}

\paragraph{}
Even if the response variables are not Gaussian, 
the power of mixed--SCGLR for model interpretation is preserved.
Graphical diagnoses similar to those provided in \autoref{Chauvet::model_interpret_simulated_data} are available in the Supplementary Material.

\section{An application to forest ecology data}
\label{Chauvet::realDATA}

\subsection{Data description}

The present study is based on the 
\textit{Genus} dataset of the CoForChange project 
(see \url{http://www.coforchange.eu}).
The subsample we consider gives the abundance of $8$ common tree genera on $2615$ Congo Basin land plots. 
These plots are grouped into $22$ forest concessions.
To predict abundances, we have
$56$ environmental variables, plus $2$ explanatory variables which code geology and  anthropogenic interference.
$\m{X}$ consists of all environmental variables which are:
\begin{itemize}[topsep=0cm, itemsep=0cm, parsep=0cm]
\item[$\blacktriangleright$]
$29$ physical factors linked to topography, rainfall or soil moisture,
\item[$\blacktriangleright$]
$25$ photosynthesis activity indicators 
(the Enhanced Vegetation Indices, EVI, 
the Near--InfraRed indices, NIR, 
and the Mid--InfraRed indices, MIR),
\item[$\blacktriangleright$]
$2$ indicators which describe the tree height.
\end{itemize}
Physical factors are many and redundant: monthly rainfalls are highly correlated, and so are photosynthesis activity indicators.
By contrast, geology and anthropogenic interference are weakly correlated and interesting per se. 
These variables are then considered as additional explanatory variables and included in matrix $\m{A}$.

\subsection{Model and parameter calibration}

\begin{figure}[!ht]
\centering
\includegraphics[width=0.9\linewidth, trim={0cm 0 0cm 0},clip]{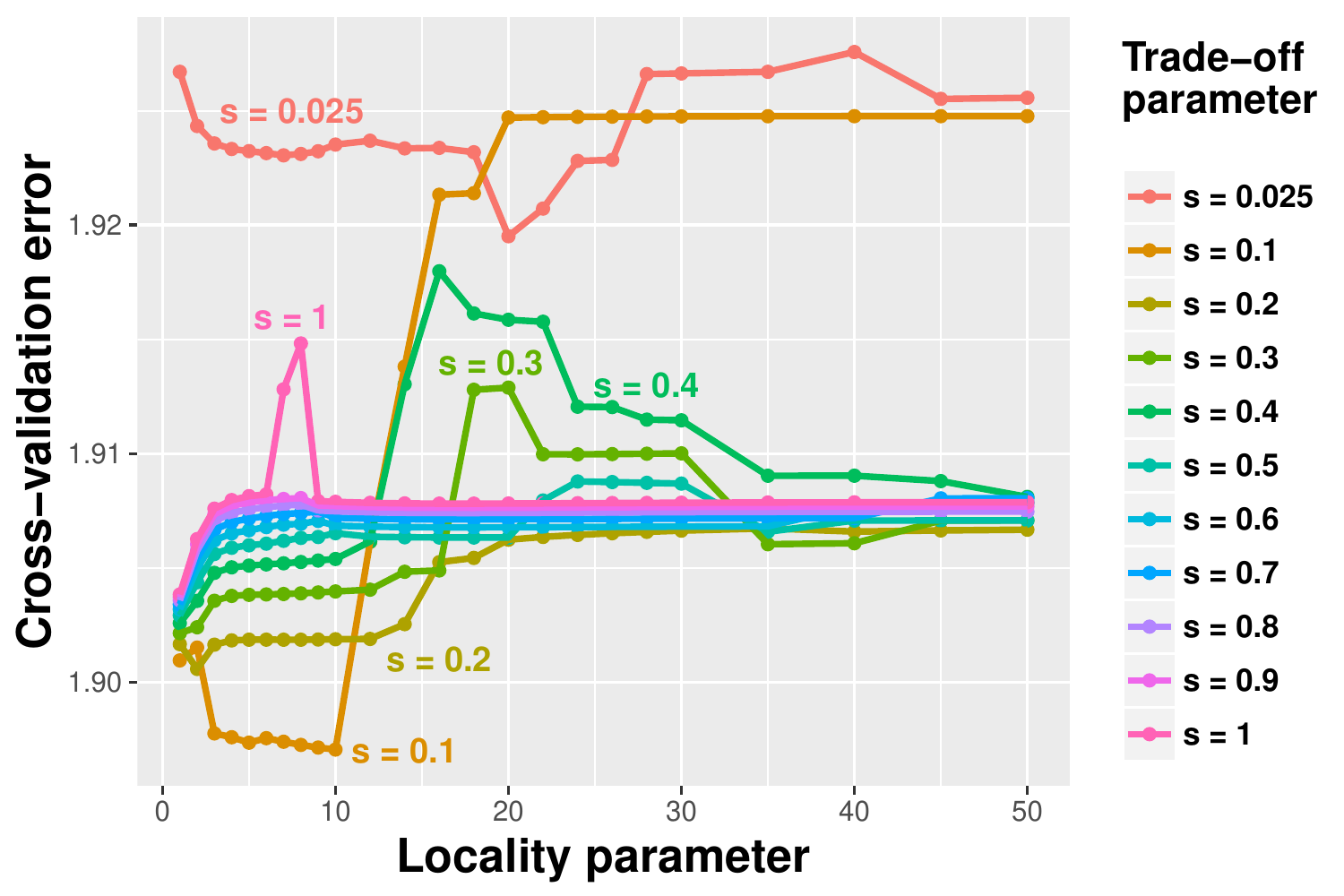}
\caption{Behaviour of the cross--validation error $E$ for trade--off parameter
$s \in \left\lbrace 0.025, 0.1, 0.2, \ldots, 1 \right\rbrace$,
as a function of locality parameter $l \in \left[1,50\right]$.}
\label{Chauvet::5CVerrors}
\end{figure}

Abundances of species given in \textit{Genus} are count data. 
For each $k \in \left\lbrace 1, \ldots, 8 \right\rbrace$, 
we consider a Poisson regression with $\log$ link 
\begin{equation*}
\m{y_k} \widesim{} \mathcal{P} \left( 
\m{\lambda} = 
\exp \left[
\sum_{j=1}^K \left( \m{Xu_{j}} \right) \gamma_{k,j} + \m{A\delta_k} + 
\m{U\xi_k}
\right]
\right),
\end{equation*}
where $\m{\xi_k}$ is the $22$--level random--effect vector 
used to model the dependence between the observations of  
$\m{y_k}$ within concessions.
The first cross--validations we performed 
--- with different fixed values of parameters $s$ and $l$ --- 
indicated that four components were sufficient to capture most of the information in $\m{X}$ needed to model and predict responses. We therefore keep $K^{\star} = 4$.
The optimal values of trade--off and locality parameter 
$s^{\star}$ and $l^{\star}$ are then determined through another cross--validation.
Using the same procedure and notations as in \autoref{Chauvet::section_parameter_calibration},
the data are divided into five parts $\mathcal{P}_{1}, \ldots, \mathcal{P}_{5}$. 
Let $n_j$ be the size of $\mathcal{P}_{j}$.
\begin{equation*}
E = \dfrac{1}{8} \sum_{k=1}^8 E_k,
\end{equation*}
where
\begin{equation}
E_k = \dfrac{1}{5} \sum_{j=1}^5 
\sqrt{
\dfrac{1}{n_j} \sum_{i \in \mathcal{P}_j} 
\dfrac{\left( y_{k,i} - \widehat{y_{k,i(-j)}} \right)^2}
{\widehat{\text{var}}\left( \widehat{y_{k,i(-j)}} \right)}
}. 
\label{Chauvet::CVerrors}
\end{equation}
On \autoref{Chauvet::5CVerrors}, we plot the errors $E$ for parameter pairs 
$(s,l) \in \mathcal{E}_s \times \mathcal{E}_l$, where 
\begin{align*}
\mathcal{E}_s &= \left\lbrace 0.025,0.1,0.2,\ldots,1 \right\rbrace \\
\mathcal{E}_l &= 
\left\lbrace 1,2,\ldots,10,12,14, \ldots, 30,35,40,45,50 \right\rbrace.
\end{align*}
Parameter grid $\mathcal{E}_s \times \mathcal{E}_l$ therefore contains $264$ pair values.
Selecting the best parameter pair from 
$\mathcal{E}_s \times \mathcal{E}_l$ through a $5$--fold cross--validation requires a computation time of about $65$ minutes (parallel computing on $6$ CPU cores, Intel Core i7--6700HQ, 2.6GHz).
It should be noted that there is a risk of non--convergence when the trade--off parameter $s$ is too close to $0$.
Indeed, if we consider no structural information ($s$ exactly equal to $0$) in $\m{X}$, mixed--SCGLR merely performs classical GLMM estimation and does not converge with this data. 
When $s=0.025$, our algorithm converges but leads to fairly unstable estimates and high cross--validation errors because regularisation is then very weak.
By contrast, the components calculated with 
$s \in \left\lbrace 0.5, 0.6, \ldots, 1 \right\rbrace$ are close to principal components. 
The associated errors are therefore stable in most cases, but rather high.
Finally, $s \in \left\lbrace 0.1, \ldots, 0.4 \right\rbrace$ leads to the lowest cross--validation errors, but only for $l \leqslant 10$. Indeed, when $s$ is not too high, mixed--SCGLR may focus on the most predictive structures of $\m{X}$. However, parameter $l$ must not exceed a certain value, in order to avoid being drawn towards too local variable--bundles.
As can be seen, choosing 
$\left( s^{\star}, l^{\star} \right) = \left( 0.1, 10 \right)$ 
minimises the cross--validation error.

\subsection{Prediction quality and interpretation results}

This part evaluates the benefits obtained by taking within--group dependence into account.
The predictions we get with mixed--SCGLR and with initial version of SCGLR are compared with respect to the cross--validation criterion given by (\ref{Chauvet::CVerrors}).
\autoref{Chauvet::improvement_prediction} summarises the $E_{k}$'s for both SCGLR and mixed--SCGLR methods.
Optimal parameter value triplet
$\left( K^{\star}, s^{\star}, l^{\star} \right) = 
\left( 4, 0.1, 10 \right)$ is selected for both methods.
For each $k \in \left\lbrace 1, \ldots, 8 \right\rbrace$, mixed--SCGLR gives a lower cross--validation error than SCGLR: taking into account the within--group dependence has clearly improved prediction performances.

\begin{table}[!ht]
\caption{Cross--validation errors for each response variable.} 
\label{Chauvet::improvement_prediction}
\vspace{3mm}
\centering
\begin{tabular}{lcccccccc} 
 & $E_{\text{cv}}^1$ & $E_{\text{cv}}^2$ & $E_{\text{cv}}^3$ & $E_{\text{cv}}^4$ & $E_{\text{cv}}^5$ & $E_{\text{cv}}^6$ & $E_{\text{cv}}^7$ & $E_{\text{cv}}^8$
\tabularnewline 
\hline \noalign{\smallskip}
SCGLR & 
\tablenum{1.32} & \tablenum{2.46} & \tablenum{3.27} & 
\tablenum{1.43} & \tablenum{2.56} & \tablenum{1.28} & 
\tablenum{1.54} & \tablenum{3.44}
\tabularnewline  
mixed--SCGLR & 
\tablenum{1.24} & \tablenum{1.95} & \tablenum{2.92} & 
\tablenum{1.32} & \tablenum{2.27} & \tablenum{1.15} & 
\tablenum{1.31} & \tablenum{3.01}
\tabularnewline
\hline
\end{tabular}
\end{table}

\begin{figure}[!ht]
\centering
\includegraphics[width=\linewidth, trim={0cm 0cm 0cm 0},clip]{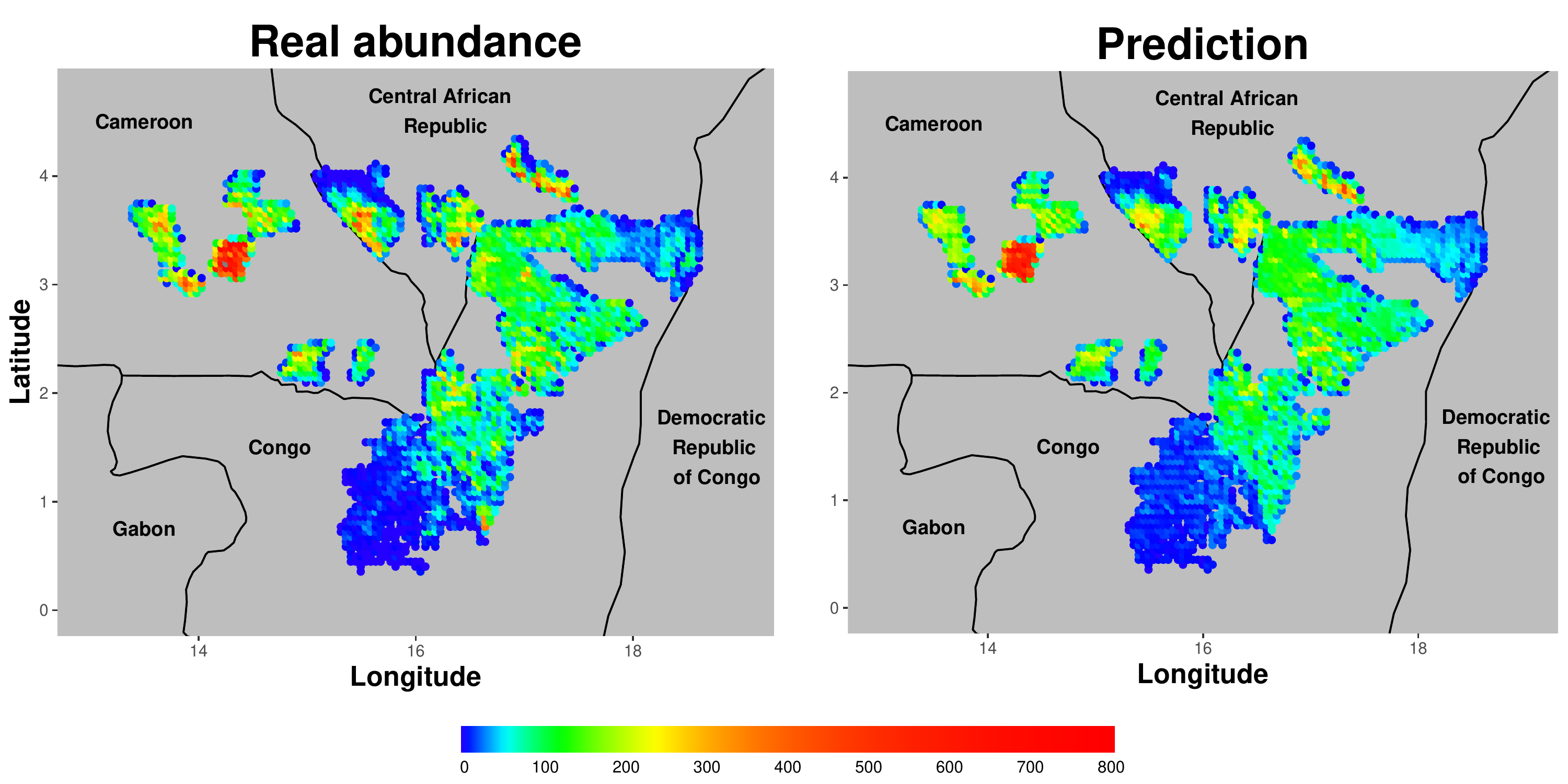}
\caption{Abundance maps issued from mixed--SCGLR. 
\textit{The plots respectively show real abundance (left) and associated conditional predictions (right) of the tree species number $8$. Each point represents a land plot (2615 in total). 
}}
\label{Chauvet::forecastmap}
\end{figure}

\noindent
Moreover, mixed--SCGLR enables to correctly reconstitute observed abundance maps, as illustrated on 
\autoref{Chauvet::forecastmap}.

\begin{figure}[!ht]
\centering
\includegraphics[width=.45\linewidth, trim={2.5cm 0 2.25cm 0},clip]{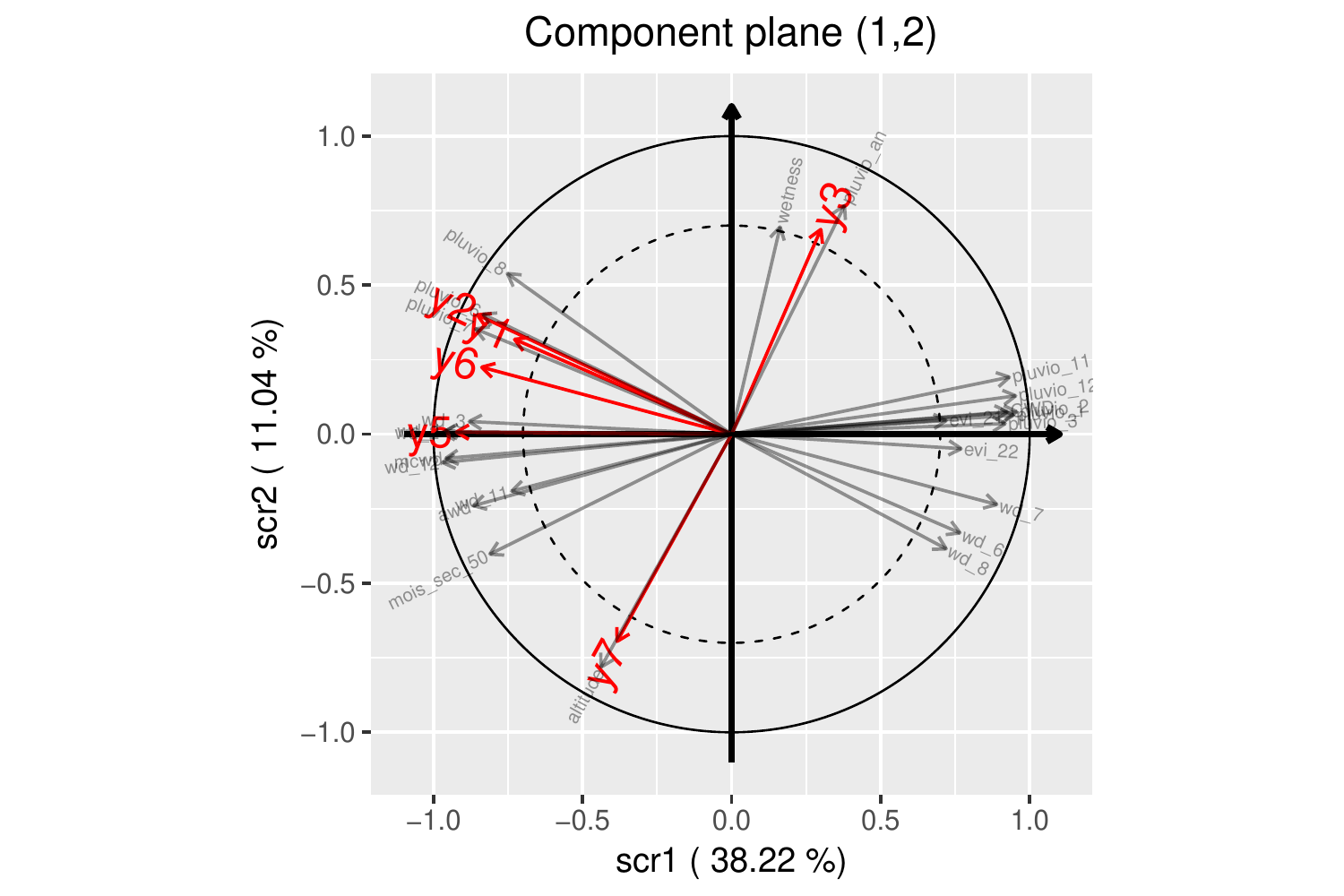}%
\hspace{5mm}
\includegraphics[width=.45\linewidth, trim={2.5cm 0 2.25cm 0},clip]{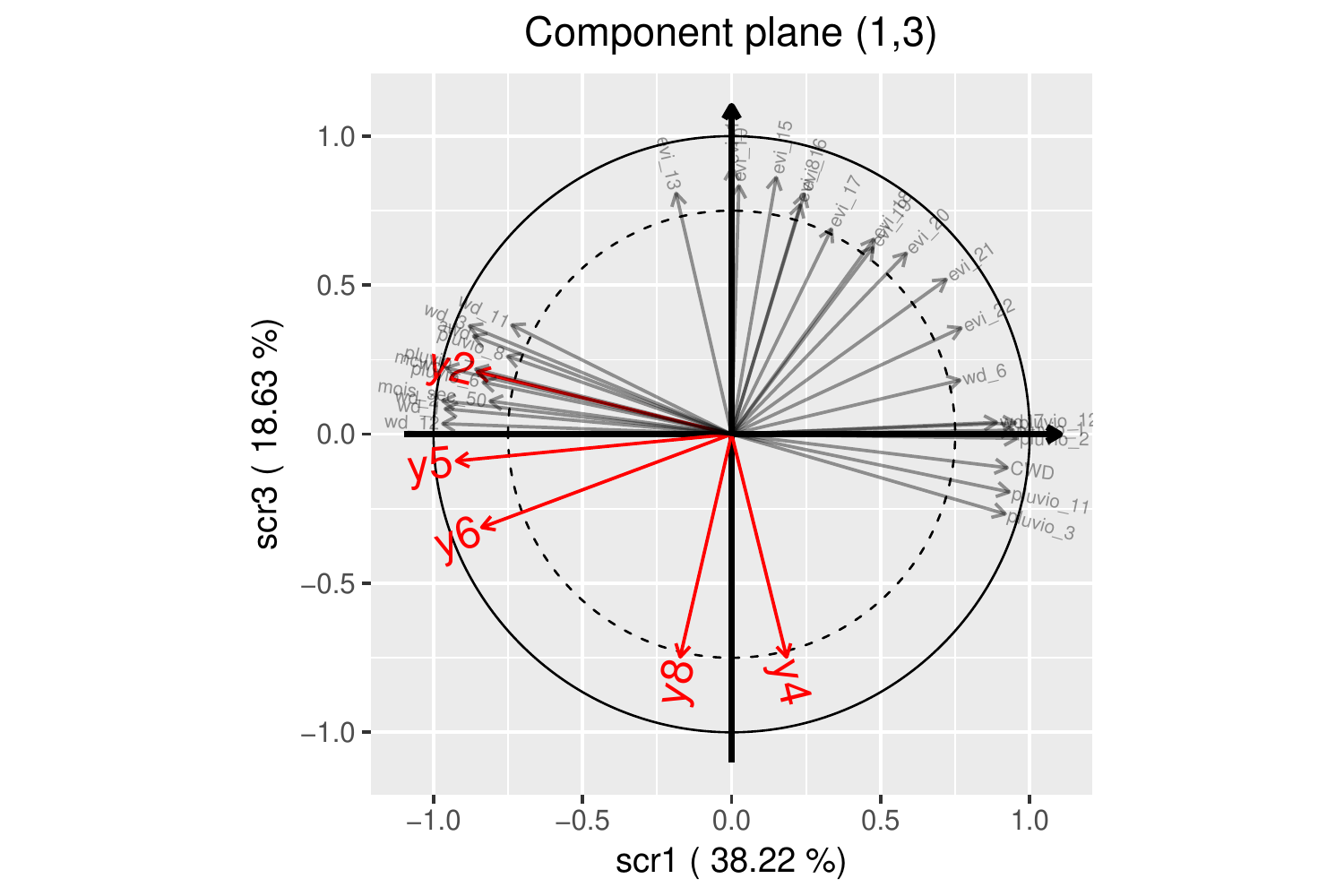}
\caption{Component planes $(1,2)$ and $(1,3)$ output by mixed--SCGLR on dataset \textit{Genus}, with optimal parameter triplet 
$\left( K^{\star}, s^{\star}, l^{\star} \right) = \left( 4, 0.1, 10 \right)$.
\textit{The left--hand side plot displays only variables having cosine greater than $0.7$ with component plane $(1,2)$. 
The right--hand side plots variables having cosine greater than $0.75$ with component plane $(1,3)$.}}
\label{Chauvet::realdata_componentplane}
\end{figure}

As has been seen in 
\autoref{Chauvet::model_interpret_simulated_data}, 
mixed--SCGLR allows an easy interpretation of the model through the decomposition of linear predictors on interpretable components.
\autoref{Chauvet::realdata_componentplane} shows the first two component planes resulting from mixed--SCGLR on real data \textit{Genus}.
Component plane $(1,2)$ reveals two patterns. 
The first one is a global rain--wind pattern driven by the \textit{pluvio}'s and \textit{wd}'s variables which explain the abundances of Species $1,2,5,6$.
The second is a rather local pattern driven by variables \textit{altitude}, \textit{wetness} and annual pluviometry (\textit{pluvio}$\_$\textit{an}) which prove important to model and predict responses $\m{y_3}$ and $\m{y_7}$.
Lastly, Component $3$ reveals a photosynthesis pattern driven by a part of the \textit{Evi}'s, which seems useful to predict $\m{y_4}$ and $\m{y_8}$.

The decomposition of linear predictors on interpretable components
allows to detect the species that tend to share common explanatory dimensions
and those which are more idiosyncratic. 
We can then identify the variable--bundles these dimensions are related to.
The underlying goal is a better understanding of the bio-- and ecosystem diversity with a view to preserve them.
Species 1, 2, 5 and 6 are sensitive to the same rain--wind regime, and Species 4 and 8 are explained by the same photosynthetic pattern. 
On the contrary, Species 3 and 7 are clearly separated. 
Species 7 grows at high altitudes where the atmosphere is rather dry while the abundance of Species 3 is favoured by regular rainfall and high humidity.


\section{Discussion and Conclusions}
\label{sec:conc}

Like Sufficient Dimension Reduction (SDR) methods, 
mixed--SCGLR is based on the construction of a reduction function of dimension less than $p$ which tries to capture all the relevant information that $\m{X}$ contains about $\m{Y}$. 
However, the two approaches do not exactly pursue the same objectives. Indeed, SDR methods look for the ``central subspace'' containing the predictive information irrespective of the structures within $\m{X}$ (e.g. dimensions capturing a large part of $\m{X}$'s variance, or bundles of correlated variables).
Mixed--SCGLR rather aims at basing the explanatory subspace on such structural dimensions so as to both gain interpretability and stabilise prediction. We think that extracting a hierarchy of strong and interpretable dimensions, and decomposing the linear predictor on them, is an essential asset in model--building. The difference in goals entails a difference in means: SDR is based on the sufficiency principle, which is enough to identify a subspace but not to track strong predictive dimensions in it. By contrast, in the wake of PLS regression, 
mixed--SCGLR uses a criterion combining goodness--of--fit and structural relevance of components.

\paragraph*{}
The supervised--component paradigm has proved effective in  situations where regularisation is necessary but where variable selection is inappropriate --- for instance when the true explanatory dimensions are latent and indirectly measured through highly correlated proxies. 
\begin{itemize}[parsep=0cm,itemsep=0cm,topsep=0cm]
\item[$\blacktriangleright$]
When $l=1$, trade--off parameter $s$ allows to continuously tune the attraction of components towards the principal components of explanatory variables. This results in a continuum between classical GLMM estimation ($s=0$ is associated with no regularisation) and principal component generalised linear mixed regression (with $s=1$).
\item[$\blacktriangleright$]
When $l>1$, we take better advantage of local predictive structures in $\m{X}$. The components we build are then usually closer to local gatherings of variables, thus easier to interpret.
\end{itemize}

\paragraph*{}
Mixed--SCGLR is able to identify more or less local predictive structures common to all the $\m{y_k}$'s and performs well on grouped data with Gaussian, Bernoulli, binomial and Poisson outcomes.
Compared to penalty--based approaches as ridge or LASSO, the orthogonal components built by mixed--SCGLR reveal 
the multidimensional explanatory and predictive dimensions, and greatly facilitate the interpretation of the model.

\paragraph*{}
However, a natural question arises as to the accuracy of our methodology under significant deviations from normality. With binary data for instance, variance component estimates are prone to some bias towards zero \citep{mcCulloch1997}. 
That is why other estimation strategies might be considered, especially Monte Carlo integration methods which have the advantage of being based on direct approximations of the likelihood. 
Some examples are the MCMC methods developed by \citet{hadfield2010} in the GLMM framework, and the Monte Carlo Likelihood Approximation (MCLA) proposed by \citet{knudson2016}. Indirect maximisations of the likelihood are also available such as Monte Carlo Expectation--Maximisation (MCEM) and Monte Carlo Newton--Raphson \citep{mcCulloch1997}. 
We think that these methodologies and Schall's could be combined sequentially. 
Indeed we could first take advantage of the linear approximation of the model in order to build the components, and then use MC--based methods to estimate both fixed--effect parameters and variance components. 
This would lead to replacing the current iteration of mixed--SCGLR given by \autoref{Chauvet::PSEUDOCODE} with the following steps (to keep things simple, we take the canonical link):  
\begin{itemize}[parsep=0cm,itemsep=0cm,topsep=0.25cm]
\item[1.] Compute components 
$\m{F} = \big[\, \m{f_1} \mid \ldots \mid \m{f_K} \,\big]$ via the PING algorithm on Schall's linearised models.
\item[2.] For each $k \in \left\lbrace 1, \ldots, q\right\rbrace$, consider the hierarchy
\begin{equation*}
\begin{aligned}
f_{\m{y_k}|\m{\xi_k},\m{\gamma_k},\m{\delta_k}} 
\left( \m{y_k}|\m{\xi_k}, \m{\gamma_k}, \m{\delta_k} \right)
&= \exp\left\lbrace
\m{y}_{\m{k}}^{\mbox{\tiny T}} \m{\eta^{\xi}_k}
- \m{1}^{\mbox{\tiny T}} 
c\left( \m{\eta^{\xi}_k} \right) 
+ \m{1}^{\mbox{\tiny T}} d\left( \m{y_k} \right)
\right\rbrace 
\\
\m{\xi_k}|\m{D_k} &\sim \mathcal{N} \left( \m{0},\m{D_k} \right),
\end{aligned}
\end{equation*}
where $\m{\eta^{\xi}_k} = 
\m{F\gamma_k} + \m{A\delta_k} + \m{U\xi_k}$, and $c$, $d$ are the functions associated with the natural parametrisation of the GLM. For example, for the  
Bernoulli--logistic regression, we have: 
$c(x)=\log(1+e^x)$ and $d(x)=0$.
\item[3.]
Apply MC--based methods such as MCMC, MCLA, MCEM or MCNR to update $\m{\gamma_k}$, $\m{\delta_k}$, $\m{\xi_k}$ and 
$\m{D_k}$, $k \in \left\lbrace 1, \ldots, q\right\rbrace$.
\item[4.]
Update working variables and weight matrices to define the  new Schall's linearised models.
\end{itemize}
Even though such MC--based methods are computationally much more intensive than the ``Joint--Maximisation'' and have intrinsic disadvantages (particularly in the assessment of convergence and in the choice of prior distributions), they could give better results in case of binary data.

\clearpage
\newpage
\bigskip
\begin{center}
{\large\bf SUPPLEMENTARY MATERIAL}
\end{center}
\begin{description}
\item[Additional simulations:] 
The first simulation reproduces that of \autoref{Chauvet::simulated_data_binary_poisson} in the case of binomial and Poisson outcomes. The second simulation explores a different structure of variable--bundles, considers Gaussian, binomial and Poisson outcomes, and presents results concerning variance component estimates. The third one involves high dimensional data.
(pdf file)
\item[Projected Iterated Normed Gradient (PING) algorithm:] 
We give some technical details about the PING algorithm, which maximises, at least locally, any criterion on the unit sphere. (pdf file)
\item[R package \texttt{mixedSCGLR}:] We provide an 
R package to perform mixed--SCGLR, also available at 
\url{https://github.com/SCnext/mixedSCGLR}. 
It contains the dataset \textit{Genus} used in \autoref{Chauvet::realDATA}. 
The package also provides demo codes, in particular for visualising the component planes (mixedSCGLR.tar.gz). 
\item[Code for running simulations:]
We also provide the R codes required to reproduce most of the simulation results (R and Rdata files). 
\end{description}

\bigskip
\begin{center}
{\large\bf ACKNOWLEDGMENTS}
\end{center}
The extended data \textit{Genus} required the arrangement and the inventory of 140.000 developed plots across four countries : Central African Republic, Gabon, Cameroon and Democratic Republic of Congo. The authors thank the members of the CoForTips project for allowing the use of this data.
We are also grateful to the editor, the associate editor and to the referees for their thorough and constructive review of this work.

\cleardoublepage
\newpage

\def\spacingset#1{\renewcommand{\baselinestretch}%
{#1}\small\normalsize} \spacingset{1}

\if0\blind
{\centering
  \title{\bf \LARGE Supplementary file to ``Component--based regularisation of multivariate generalised linear mixed models'': \\
  \vspace{5mm}
 \centering THE PING ALGORITHM} \\
 \vspace{7mm}
  \author{
\textbf{Jocelyn Chauvet}$^{1}$ \quad \and
\textbf{Catherine Trottier}$^{1,2}$ \quad \and
\textbf{Xavier Bry}$^{1}$ \\
\vspace{3mm}
$^{1}$ IMAG, Univ Montpellier, CNRS, Montpellier, France. \\ 
\href{mailto:jocelyn.chauvet@umontpellier.fr}{jocelyn.chauvet@umontpellier.fr} ; \href{mailto:xavier.bry@umontpellier.fr}{xavier.bry@umontpellier.fr} \\
$^{2}$ Univ Paul-Val\'ery Montpellier 3, Montpellier, France.
\vspace{-2.5mm}
\begin{center}
\href{mailto:catherine.trottier@univ-montp3.fr}{catherine.trottier@univ-montp3.fr}
\end{center}
}
  \maketitle
} \fi

\if1\blind
{
  \bigskip
  \bigskip
  \bigskip
  \begin{center}
    {\LARGE\bf Component-based regularisation of multivariate generalised linear mixed models}
\end{center}
  \medskip
} \fi

\bigskip

\spacingset{1.45} 

\spacingset{1.45} 

\vspace{1cm}
The Projected Iterated Normed Gradient (PING) is a basic extension of the iterated power algorithm, for solving any program of the form
\begin{equation}
\left\lbrace
\begin{aligned}
&\text{max} \quad  \mathcal{J}_h \left( \m{u} \right),    \\
&\text{subject to:}  
\quad \m{u}^{\mbox{\tiny T}} \m{M}^{-1} \m{u} = 1 \; \text{and} \; 
\m{\Delta}_{\m{h}}^{\mbox{\tiny T}} \m{u} = \m{0}.
\end{aligned}
\right.
\label{Chauvet::genericPING1}
\end{equation}
Note that putting $\m{v} = \m{M}^{-1/2} \m{u}$, 
$\mathcal{G}_h \left( \m{v} \right) = \mathcal{J}_h \left( \m{M}^{1/2} \m{v} \right)$ and 
$\m{B_h} = \m{M}^{1/2} \m{\Delta_h}$, 
Program (\ref{Chauvet::genericPING1}) is strictly equivalent to Program (\ref{Chauvet::genericPING2}):
\begin{equation}
\left\lbrace
\begin{aligned}
&\text{max} \quad  \mathcal{G}_h \left( \m{v} \right),    \\
&\text{subject to:}  \quad \m{v}^{\mbox{\tiny T}} \m{v} = 1 \; \text{and} \; 
\m{B}_{\m{h}}^{\mbox{\tiny T}} \m{v} = \m{0}.
\end{aligned}
\right.
\label{Chauvet::genericPING2}
\end{equation}
Denoting 
\begin{align*}
\Pi_{\spn \left\lbrace \m{B}_{\m{h}} \right\rbrace^{\bot}} 
&= \m{I} - \m{B_h}
\left( \m{B}_{\m{h}}^{\mbox{\tiny T}} \m{B}_{\m{h}} \right)^{-1} \m{B}_{\m{h}}^{\mbox{\tiny T}} \;\;
\text{  and  } \\
\Gamma_h \left( \m{v} \right) &= \underset{\m{v}}{\nabla} \: 
\mathcal{G}_h \left( \m{v} \right),
\end{align*}
a Lagrange multiplier-based reasoning gives the basic iteration of the PING algorithm:
\begin{equation}
\m{v}^{[t+1]} = 
\dfrac{
\Pi_{\spn \left\lbrace \m{B_h} \right\rbrace^{\bot}} \, 
\Gamma_h \left( \m{v}^{[t]} \right)
}
{\left\lVert 
\Pi_{\spn \left\lbrace \m{B_h} \right\rbrace^{\bot}} \,
\Gamma_h \left( \m{v}^{[t]} \right)
\right\rVert}.
\label{Chauvet::basinPINGiter}
\end{equation}
Although Iteration (\ref{Chauvet::basinPINGiter}) follows a direction of ascent, it does not guarantee that $\mathcal{G}_h$ actually increases on every step.
We therefore propose a 
generic iteration of PING (\autoref{Chauvet::algoPING}) and an alternative one (\autoref{Chauvet::algoPING2}), which both ensure that the criterion increases.

\vspace{10mm}
\begin{center}
\fbox{\begin{minipage}{0.95\linewidth}
\begin{algorithm}[H]
\caption{Generic iteration of the PING algorithm}
\bigskip
\While{convergence of $\m{v}$ non reached}{
\smallskip   
$\m{\kappa}^{[t]} = 
\dfrac{
\Pi_{\spn \left\lbrace \m{B_h} \right\rbrace^{\bot}} \, 
\Gamma_h \left( \m{v}^{[t]} \right)
}
{\left\lVert 
\Pi_{\spn \left\lbrace \m{B_h} \right\rbrace^{\bot}} \,
\Gamma_h \left( \m{v}^{[t]} \right)
\right\rVert}$ \\
\bigskip
A unidimensional Newton--Raphson maximisation procedure is used to find  \\
\vspace{-1mm}
the maximum of $\mathcal{G}_h \left( \m{v} \right)$ on the arc 
$\left( \m{v}^{[t]}, \m{\kappa}^{[t]} \right)$ and take it as 
$\m{v}^{[t+1]}$. \\
\smallskip
$t \longleftarrow t+1$
}
\label{Chauvet::algoPING}
\end{algorithm}
\end{minipage}}
\end{center}

\vspace{10mm}
\begin{center}
\fbox{\begin{minipage}{0.95\linewidth}
\begin{algorithm}[H]
\caption{Alternative generic iteration of the PING algorithm}
\bigskip
\While{convergence of $\m{v}$ non reached}{
\smallskip
$
\m{m} \longleftarrow  
\dfrac{
\Pi_{\spn \left\lbrace \m{B_h} \right\rbrace^{\bot}} \, 
\Gamma_h \left( \m{v}^{[t]} \right)
}
{\left\lVert 
\Pi_{\spn \left\lbrace \m{B_h} \right\rbrace^{\bot}} \,
\Gamma_h \left( \m{v}^{[t]} \right)
\right\rVert}
$ \\
\medskip
\While{$\mathcal{G}_h \left( \m{m} \right) < \mathcal{G}_h \left( \m{v}^{[t]} \right)$}{
$\m{m} \longleftarrow 
\dfrac{\m{v}^{[t]} + \m{m}}
{\left\lVert \m{v}^{[t]} + \m{m} \right\rVert}$
}
$\m{v}^{[t+1]} = \m{m}$ \\
\smallskip
$t \longleftarrow t+1$
}
\medskip
\label{Chauvet::algoPING2}
\end{algorithm}
\end{minipage}}
\end{center}

\bigskip
\paragraph*{First rank component.}
Component $\m{f_{1}} = \m{X} \m{u_{1}}$ is obtained by solving
\begin{equation*}
\left\lbrace
\begin{aligned}
&\text{max} \quad  
s \log \left[\phi\left(\m{u}\right)\right] + 
(1-s) \log \left[\psi_{\m{A}}\left(\m{u}\right)\right]    \\
&\text{subject to:}  \quad \m{u}^{\mbox{\tiny T}} \m{M}^{-1} \m{u} = 1.
\end{aligned}
\right.
\end{equation*}
This corresponds to Program (\ref{Chauvet::genericPING1}) with $h=0$, where
\begin{itemize}[topsep=0cm, itemsep=0cm, parsep=0cm]
\item[$\blacktriangleright$]
$\mathcal{J}_0 \left( \m{u} \right) = 
s \log \left[\phi\left(\m{u}\right)\right] + 
(1-s) \log \left[\psi_{\m{A_0}}\left(\m{u}\right)\right]$, 
\item[$\blacktriangleright$]
$\m{A_0}=\m{A}$ (the matrix of additional explanatory variables), and 
\item[$\blacktriangleright$]
$\m{\Delta_0}=\m{0}$.
\end{itemize}
In this particular case, we have $\m{B_0}=\m{M}^{1/2}\m{\Delta_0}=\m{0}$, and so: 
$$\Pi_{\spn \left\lbrace \m{B_0} \right\rbrace^{\bot}} = \m{I}.$$

\paragraph*{Higher rank components.}
Let ${\m{F_h}} = \big[\, \m{f_1} \mid \ldots \mid \m{f_h} \,\big]$ be the matrix of the first $h$ components and $\m{A_h} = \big[\, \m{F_h} \mid \m{A} \,\big]$.
Let $\m{P}$ denote the weight matrix reflecting the a priori relative importance of observations 
($\m{P}=\frac{1}{n} \m{I}_n$ if all observations are of equal importance).
Component $\m{f_{h+1}} = \m{X} \m{u_{h+1}}$ is obtained by solving
\begin{equation*}
\left\lbrace
\begin{aligned}
&\text{max} \quad  
s \log \left[\phi\left(\m{u}\right)\right] + 
(1-s) \log \left[\psi_{\m{A_h}}\left(\m{u}\right)\right]    \\
&\text{subject to:}  \quad \m{u}^{\mbox{\tiny T}} \m{M}^{-1} \m{u} = 1 \; \text{and} \; 
\m{F}_{\m{h}}^{\mbox{\tiny T}} \m{P} \m{Xu} = \m{0}.
\end{aligned}
\right.
\end{equation*}
This corresponds to Program (\ref{Chauvet::genericPING1}), where
\begin{itemize}[topsep=0cm, itemsep=0cm, parsep=0cm]
\item[$\blacktriangleright$]
$\mathcal{J}_h \left( \m{u} \right) = 
s \log \left[\phi\left(\m{u}\right)\right] + 
(1-s) \log \left[\psi_{\m{A_h}}\left(\m{u}\right)\right]$, 
\item[$\blacktriangleright$]
$\m{A_h}=\big[\, \m{F_h} \mid \m{A} \,\big]$, and 
\item[$\blacktriangleright$]
$\m{\Delta_h}=\m{X}^{\mbox{\tiny T}} \m{P} \m{F_h}$.
\end{itemize}

\paragraph*{Initialisation.}
To quickly find $\m{f_1}$, algorithm PING is initialised with the first PLS component of the responses on $\m{X}$.
In like manner, for $h \geqslant 2$, 
PING is initialised with the first PLS component of the responses on $\m{X}$ deflated on components 
$\m{F_{h-1}}= \big[\, \m{f_1} \mid \ldots \mid \m{f_{h-1}} \,\big]$.

\cleardoublepage
\newpage

\def\spacingset#1{\renewcommand{\baselinestretch}%
{#1}\small\normalsize} \spacingset{1}

\if0\blind
{\centering
  \title{\bf \LARGE Supplementary file to ``Component--based regularisation of multivariate generalised linear mixed models'': \\
  \vspace{5mm}
 \centering ADDITIONAL SIMULATIONS} \\
 \vspace{7mm}
  \author{
\textbf{Jocelyn Chauvet}$^{1}$ \quad \and
\textbf{Catherine Trottier}$^{1,2}$ \quad \and
\textbf{Xavier Bry}$^{1}$ \\
\vspace{3mm}
$^{1}$ IMAG, Univ Montpellier, CNRS, Montpellier, France. \\ 
\href{mailto:jocelyn.chauvet@umontpellier.fr}{jocelyn.chauvet@umontpellier.fr} ; \href{mailto:xavier.bry@umontpellier.fr}{xavier.bry@umontpellier.fr} \\
$^{2}$ Univ Paul-Val\'ery Montpellier 3, Montpellier, France.
\vspace{-2.5mm}
\begin{center}
\href{mailto:catherine.trottier@univ-montp3.fr}{catherine.trottier@univ-montp3.fr}
\end{center}
}
  \maketitle
} \fi

\if1\blind
{
  \bigskip
  \bigskip
  \bigskip
  \begin{center}
    {\LARGE\bf Component-based regularisation of multivariate generalised linear mixed models}
\end{center}
  \medskip
} \fi

\bigskip

\spacingset{1.45} 

\section{Comparative results with binomial and Poisson outcomes}

In this section, we simply extend the simulation scheme presented in \autoref{Chauvet::simulated_data_binary_poisson} to binomial and Poisson outcomes. We maintain design matrices $\m{X}$ and $\m{U}$ as defined in \autoref{Chauvet::datageneration}. 
Fixed--effect parameters $\m{\theta_k}$'s and random--effect vectors $\m{\widetilde{\xi_k}}$'s are defined in \autoref{Chauvet::simulated_data_binary_poisson}. Given $\m{\widetilde{\xi_1}}$ and $\m{\widetilde{\xi_2}}$, we then simulate
$\m{Y} = \big[\, \m{y_1} \mid \m{y_2} \,\big]$ as
\begin{equation*}
\left\lbrace
\begin{aligned}
&\m{y_1} \widesim{} \mathcal{B}\text{in} 
\left( \textbf{trials} = 50 \, \mathbf{1}_n, \,
\m{p} = \text{logit}^{-1} \left[ 
\m{X \theta_1} + \m{U \widetilde{\xi_1}} 
\right] \right) \\
&\m{y_2} \widesim{} \mathcal{P} \left( 
\m{\lambda} = \exp \left[
\m{X \theta_2} + \m{U \widetilde{\xi_2}} \right]
\right).
\end{aligned}
\right.
\label{Chauvet:modele_binomial_poisson}
\end{equation*}
\autoref{Chauvet::table_simu_binomial_poisson} gives the Mean Relative Squared Error (MRSE) values for 
$\m{\theta_1}$ and $\m{\theta_2}$ obtained on 100 samples for each value of 
$\tau$.

\paragraph*{}
For the Poisson distribution, the results in \autoref{Chauvet::table_simu_binomial_poisson} are essentially identical to those in the article: mixed--SCGLR outperforms GLMM--LASSO
\citep[R package \texttt{glmmLasso}]{packageGLMMLASSO} except for $\tau=0.1$.
As for the binomial distribution, the regularisation provided by mixed--SCGLR improves the results obtained without regularisation 
\citep[R package \texttt{lme4}]{packagelme4}, regardless of the level of redundancy within the explanatory variables.
Unsurprisingly, the errors are much smaller than in the binary case.

\begin{table}[!ht]
\caption{Mean Relative Squared Error (MRSE) values obtained with binomial and Poisson distributions. The R package \texttt{glmmLasso} does not handle binomial outcomes but only Bernoulli ones, which precludes comparison in this case.} 
\vspace{3mm}
\begin{tabular}{l p{2cm} p{2cm} p{3.5cm} p{2cm} p{2cm}}
& 
\multicolumn{2}{c}{\centering \textbf{GLMM}} & 
\multirow{2}{*}{\textbf{GLMM--LASSO}} & 
\multicolumn{2}{c}{\multirow{2}{*}{\textbf{mixed--SCGLR}}} \\
 & \multicolumn{2}{c}{\centering (no regularisation)} & & &  \\
& 
\centering  Binomial & 
\centering  Poisson & 
\centering  Poisson & 
\centering  Binomial & 
\centering  Poisson
\tabularnewline 
\hline \noalign{\smallskip}
$\tau=0.1$ & 
\centering  \tablenum{2.31} & \centering \tablenum{0.50} & 
\centering \tablenum{0.31} & 
\centering \tablenum{0.51} & \centering \tablenum{0.45}   
\tabularnewline  
$\tau=0.3$ & 
\centering  \tablenum{3.07} & \centering \tablenum{0.60} & 
\centering \tablenum{0.33} & 
\centering \tablenum{0.28} & \centering \tablenum{0.18}  
\tabularnewline 
$\tau=0.5$ & 
\centering  \tablenum{3.93} & \centering \tablenum{0.75} & 
\centering \tablenum{0.39} & 
\centering \tablenum{0.15} & \centering \tablenum{0.09}   
\tabularnewline 
$\tau=0.7$ & 
\centering  \tablenum{6.50} & \centering \tablenum{1.07} & 
\centering \tablenum{0.40} & 
\centering \tablenum{0.10} & \centering \tablenum{0.07}   
\tabularnewline 
$\tau=0.9$ & 
\centering  \tablenum{19.29} & \centering \tablenum{2.71} & 
\centering \tablenum{0.42} & 
\centering \tablenum{0.07} & \centering \tablenum{0.05}   
\tabularnewline 
\hline
\end{tabular}
\label{Chauvet::table_simu_binomial_poisson}
\end{table}

\begin{figure}[!ht]
\centering
\includegraphics[width=.45\linewidth, trim={2.5cm 0 2.25cm 0},clip]{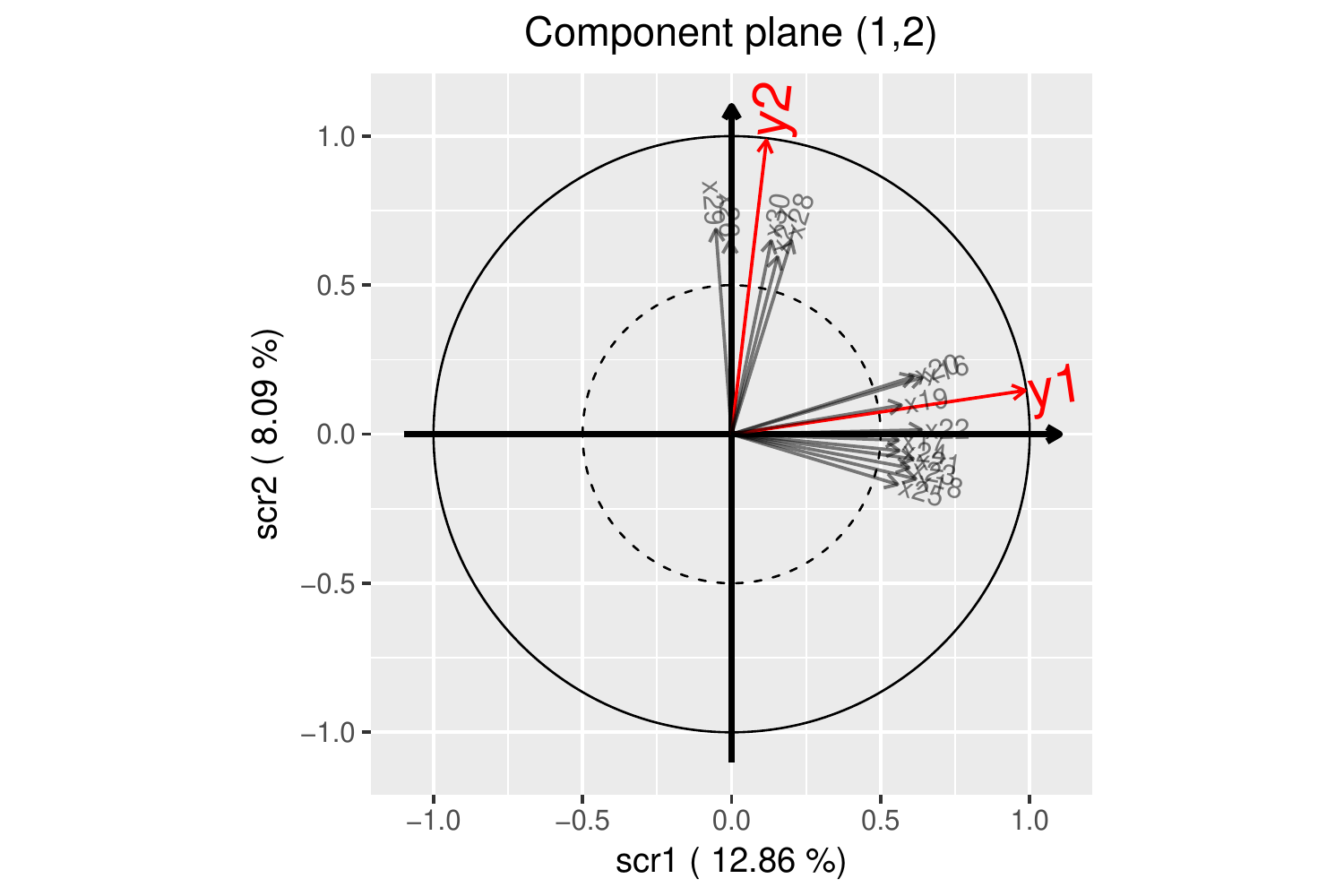}
\hspace{5mm}
\includegraphics[width=.45\linewidth, trim={2.5cm 0 2.25cm 0},clip]{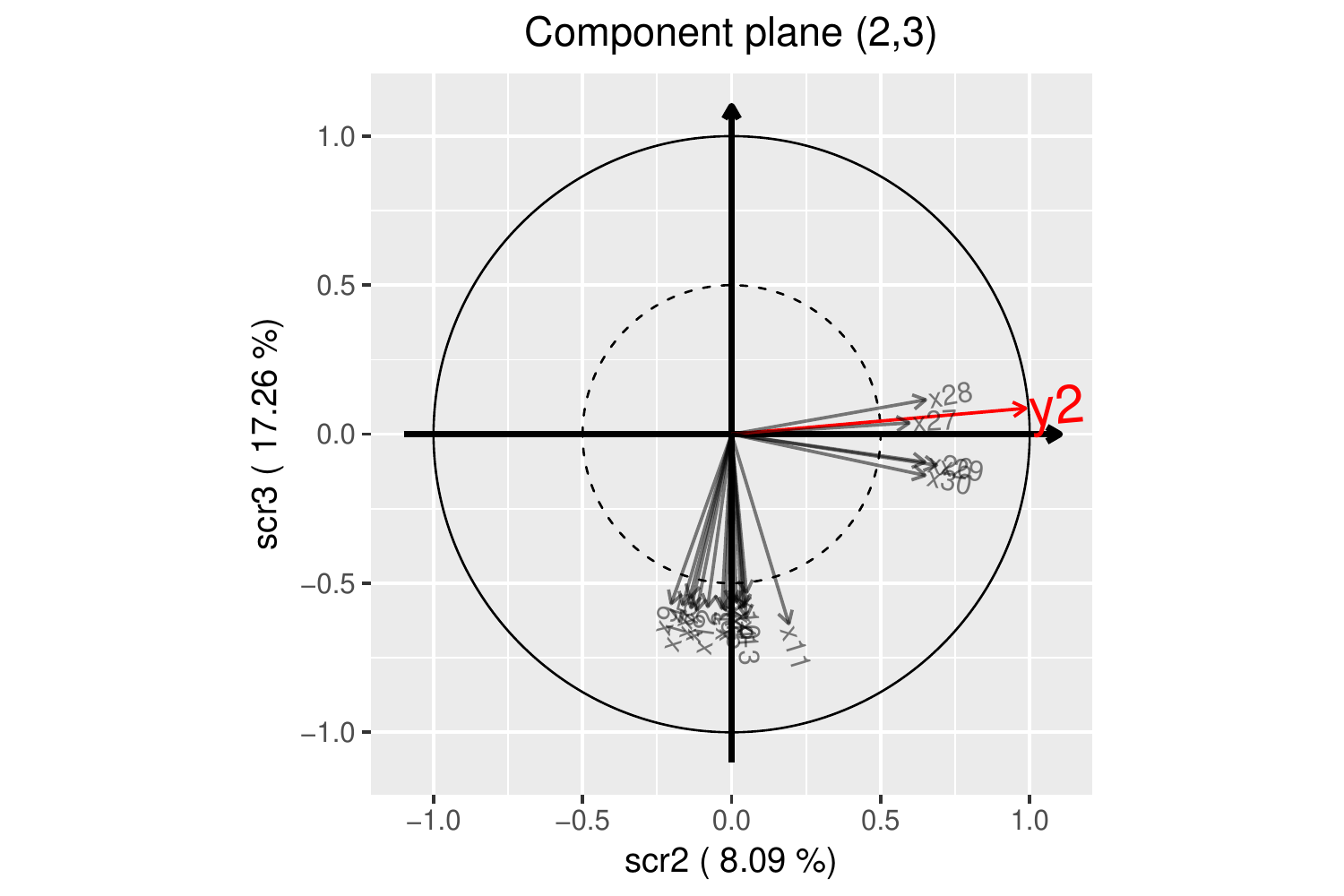}
\caption{Example of component planes given by mixed--SCGLR in the binomial/Poisson case for 
$\tau=0.3$, with parameter triplet
$\left( K, s, l \right) = \left( 3, 0.5, 2 \right)$.}
\label{Chauvet::binomial_poisson_componentplane}
\end{figure}

\begin{figure}[!ht]
\centering
\includegraphics[width=.45\linewidth, trim={2.5cm 0 2.25cm 0},clip]{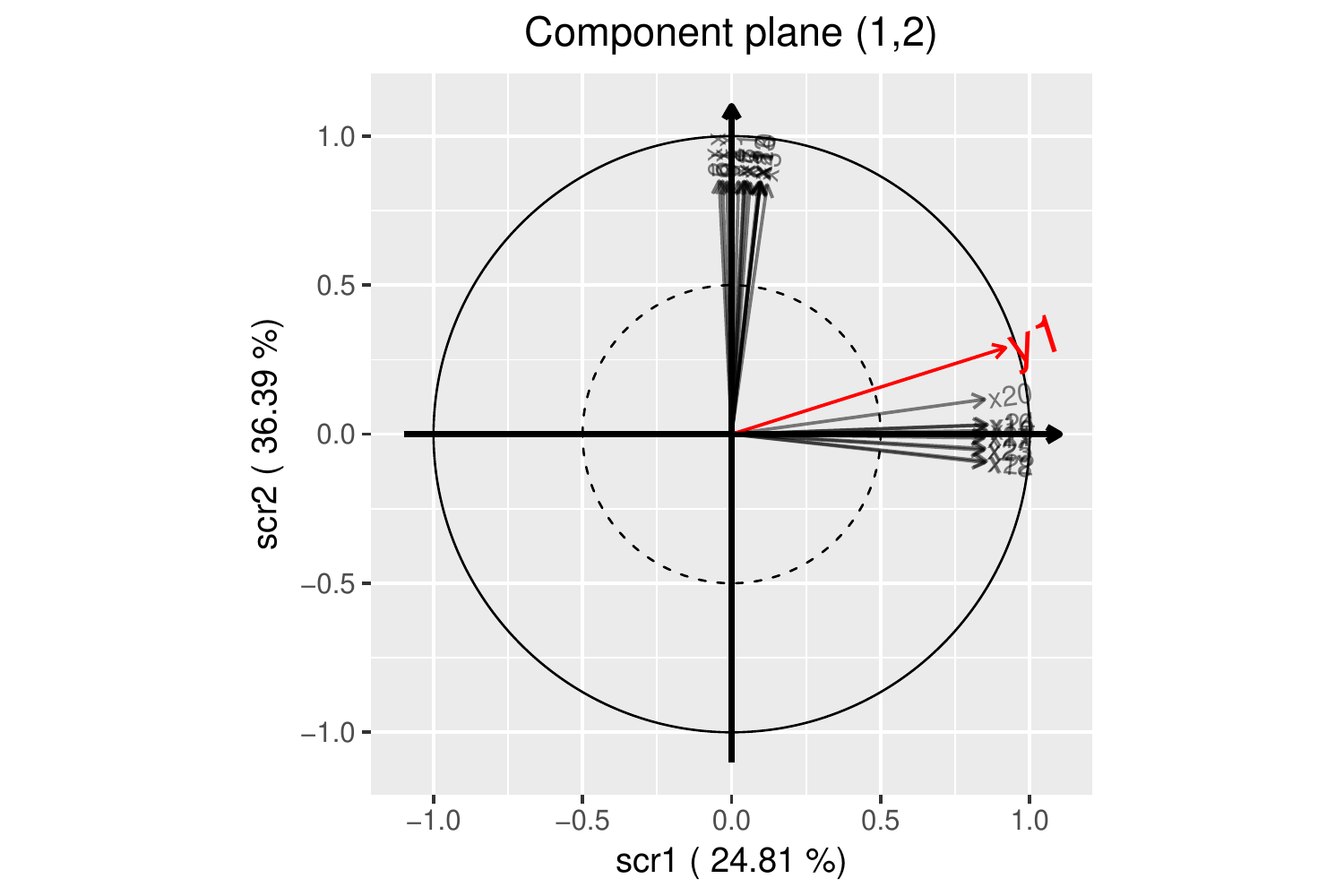}
\hspace{5mm}
\includegraphics[width=.45\linewidth, trim={2.5cm 0 2.25cm 0},clip]{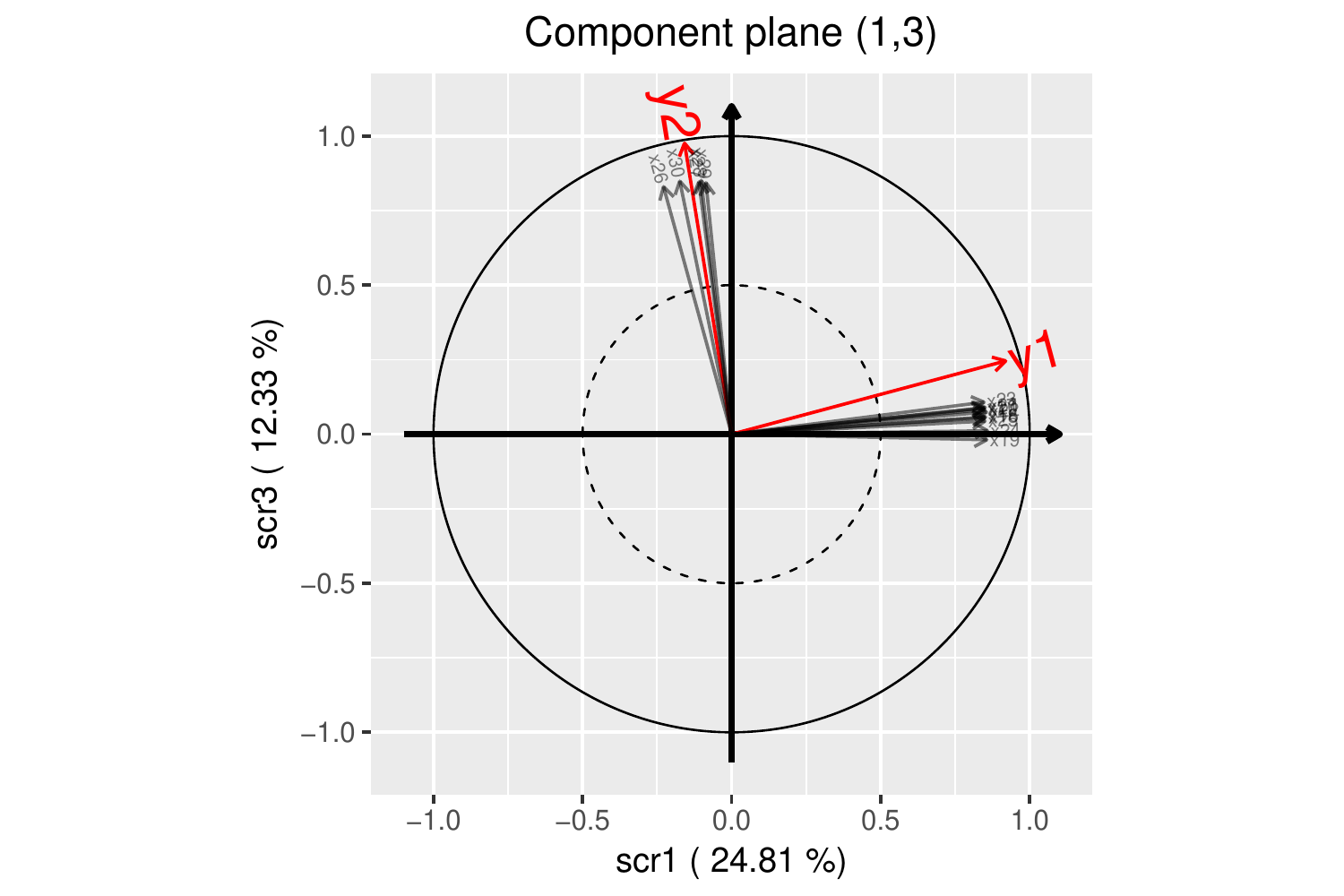}
\caption{Example of component planes given by mixed--SCGLR in the Bernoulli/Poisson case for $\tau=0.7$, with parameter triplet
$\left( K, s, l \right) = \left( 3, 0.9, 4 \right)$.}
\label{Chauvet::bernoulli_poisson_componentplane}
\end{figure}

The power of mixed--SCGLR in terms of model interpretation remains the same for non--Gaussian outcomes. \autoref{Chauvet::binomial_poisson_componentplane} (respectively \autoref{Chauvet::bernoulli_poisson_componentplane}) presents an example of the first component planes output by mixed--SCGLR in the binomial/Poisson (respectively Bernoulli/Poisson) case.
As for Gaussian outcomes, the component planes reveal that $\m{y_1}$ is explained by bundle 
$\m{X_1}$ and $\m{y_2}$ by $\m{X_2}$.
In the binomial/Poisson case with $\tau=0.3$ (\autoref{Chauvet::binomial_poisson_componentplane}), predictive bundles 
$\m{X_1}$ and $\m{X_2}$ are captured respectively by the first and the second components. The third component aligns on nuisance bundle $\m{X_0}$, despite its high inertia.
\autoref{Chauvet::bernoulli_poisson_componentplane} illustrates what may happen when the level of redundancy is very high ($\tau=0.7$ here). Since the explanatory variables are highly correlated, mixed--SCGLR regularisation requires that the structural relevance be given a heavy weight with respect to the goodness--of--fit, which leads to a trade--off parameter $s$ close to 1 ($s=0.9$ here). 
Having the greatest structural strength, the nuisance bundle is captured by the second component despite its lack of explanatory power.
This is sometimes the price to be paid for the trade--off.
In our example, the second explanatory bundle is captured by the third component, so that the predictive dimensions are accurately represented in component plane $(1,3)$.

\section{A new structure for the bundles --- Gaussian, Poisson and binomial distributions}
This simulation study tests mixed-SCGLR on a slightly more complex bundle structure. 
Results concerning variance component estimates are also presented.

\paragraph*{}
We consider a fixed--effect design matrix $\m{X}_{n \times p}$ partitioned into 3 blocks $\m{\mathcal{G}_1}$, $\m{\mathcal{G}_2}$ and $\m{\mathcal{G}_3}$. 
Block $\m{\mathcal{G}_1}$ contains 10 predictive explanatory variables structured about a latent variable 
$\m{\varphi_1} \widesim{} \mathcal{N}_n 
\left( \m{0}, \sigma_{\text{LV}}^2 \m{I}_n \right)$. 
Thus for each
$j \in \left\lbrace 1, \ldots, 10 \right\rbrace $,
$\m{x_j} = \m{\varphi_1} + \m{\varepsilon_j}$, where 
$\m{\varepsilon_j} \widesim{} \mathcal{N}_n 
\left( \m{0}, \sigma_{\text{noise}}^2 \m{I}_n \right)$ such that 
$\sigma_{\text{LV}}^2 + \sigma_{\text{noise}}^2 = 1$.
The correlation within $\m{\mathcal{G}_1}$ is tuned by signal to noise (StN) ratio 
$\sigma_{\text{LV}}^2 / \sigma_{\text{noise}}^2$
(chosen in $\left\lbrace \frac{1}{3},1,3 \right\rbrace$ in practice).
$\m{\mathcal{G}_2}$ contains a single predictive variable 
$\m{\varphi_2} = \m{x_{11}} \widesim{} \mathcal{N}_n \left( \m{0}, \m{I}_n \right)$. 
$\m{\mathcal{G}_3}$ contains 20 unstructured noise variables: for each 
$j \in \left\lbrace  12, \ldots, 31 \right\rbrace$,
$\m{x_j} \widesim{} \mathcal{N}_n \left( \m{0}, \m{I}_n \right)$.
For each $k \in \left\lbrace 1,2,3 \right\rbrace$, random--effect vectors are simulated as  
$\m{\xi_k} \widesim{\text{ind.}} \mathcal{N}_N \left( \m{0}, \sigma_k^2 \, \m{I}_N \right)$.
Given $\m{\xi_1}$, $\m{\xi_2}$, $\m{\xi_3}$, 
we simulate 3 responses having different distributions, 
$\m{Y} = \big[\, \m{y_1} \mid \m{y_2} \mid \m{y_3} \,\big]$, as
\begin{equation*}
\left\lbrace
\begin{aligned}
&\m{y_1} \widesim{} \mathcal{N}_n \Big(  
\m{\mu} = \alpha_1 \m{\varphi_1} + \m{U \xi_1}, \,
\m{\Sigma} = \m{I}_n
\Big) \\
&\m{y_2} \widesim{} \mathcal{P} \Big( 
\m{\lambda} = \exp \Big[
\alpha_2 \m{\varphi_2} + \m{U \xi_2} \Big]
\Big) \\
&\m{y_3} \widesim{} \mathcal{B}\text{in} \Big( 
\textbf{trials} = 25 \, \mathbf{1}_n, \,
\m{p} = \text{logit}^{-1} \Big[ 
\alpha_3 \left( \m{\varphi_1}+\m{\varphi_2} \right) 
+ \m{U \xi_3}
\Big] \Big).
\end{aligned}
\right.
\label{Chauvet:sim_responses_newbundles}
\end{equation*}
In our simulations, we set 
$\alpha_1 = \sigma_1^2 = 2$,
$\alpha_2 = \sigma_2^2= 1$ and 
$\alpha_3 =\sigma_3^2 = 0.5$.

\paragraph*{}
We consider in turn $N=10$ and $N=50$ groups, and $R=10$ observations per group 
($n=100$ and $n=500$ observations in total). 
$B=100$ samples are generated for each pair of values 
$\left( N,\text{StN} \right)$.
The main goal of the study is to assess the ability of mixed--SCGLR to track down both latent variable $\m{\varphi_1}$ and predictive variable
$\m{\varphi_2}$. For $j=1$ and $2$, we then define
\begin{equation*}
\text{cor}_j = \dfrac{1}{B} \sum_{b=1}^B
\left|
\text{cor} \left( \m{\varphi_j}, \m{f}_{\m{j}}^{(b)} \right)
\right|,
\end{equation*}
where $\m{f}_{\m{j}}^{(b)}$ is the component most correlated with $\m{\varphi_j}$
issued from mixed--SCGLR in the $b$--th sample.
Consistency of fixed--effect estimates is assessed through criteria 
$\text{err}_1$, $\text{err}_2$ and $\text{err}_3$ defined by
\begin{equation*}
\begin{aligned}
\text{err}_j &= \dfrac{1}{B} \sum_{b=1}^B
\dfrac{\left\lVert \alpha_j \m{\varphi_j} - 
\m{X} \widehat{\m{\beta}}_{\m{j}}^{(b)} \right\rVert^2}{\left\lVert \alpha_j \m{\varphi_j} \right\rVert^2}, \;
j \in \left\lbrace 1,2 \right\rbrace
\\
\text{err}_3 &= \dfrac{1}{B} \sum_{b=1}^B
\dfrac{\left\lVert \alpha_3 \left( \m{\varphi_1}+\m{\varphi_2} \right) - \m{X} \widehat{\m{\beta}}_{\m{3}}^{(b)} \right\rVert^2}{\left\lVert \alpha_3 \left( \m{\varphi_1}+\m{\varphi_2} \right) \right\rVert^2},
\end{aligned}
\end{equation*}
where $\widehat{\m{\beta}}_{\m{j}}^{(b)}$ is the fixed--effect estimate related to response $\m{y_j}$ associated with sample $b$.

\paragraph*{}
\autoref{Chauvet::table_new_bundles} summarises the values of the afore--defined criteria and presents biases and standard errors of variance components estimates.
For a given value of $N$, $\text{cor}_1$ increases towards $1$ with ratio $\sigma_{\text{LV}}^2 / \sigma_{\text{noise}}^2$:
the tighter the block $\m{\mathcal{G}_1}$ is structured about its latent variable, the better mixed--SCGLR can reconstruct it. The associated criterion $\text{err}_1$ then naturally decreases towards $0$.
On the other hand, $\text{cor}_2$ and $\text{err}_2$ are very stable, which proves that mixed--SCGLR is able to detect an isolated predictive variable among a large number of irrelevant others.
As $\text{err}_3$ depends on how accurately mixed--SCGLR recovers $\m{\varphi_1}$ and $\m{\varphi_2}$, it slightly decreases when the StN ratio increases.
Both variance component biases and standard errors seem rather stable regardless of the value of StN.
Finally, when $N$ increases, all the $\text{cor}_j$'s increase towards 1 and all the $\text{err}_j$'s decrease towards 0. 
As far as variance component estimates are concerned, the biases are getting slightly closer to 0 and the standard errors decrease significantly.

\begin{table}[!ht]
\caption{Summary of $\text{cor}_j$ and $\text{err}_j$ values, and presentation of biases and standard errors of estimated variance components.} 
\vspace{3mm}
\centering
\begin{tabular}{l p{2cm} p{1.2cm} p{2cm} p{2cm} p{1.2cm} p{1.5cm} }
& 
\multicolumn{3}{c}{\centering $N=10, R=10 \; (n=100)$} & 
\multicolumn{3}{c}{\centering $N=50, R=10 \; (n=500)$} \\
$\sigma_{\text{LV}}^2 / \sigma_{\text{noise}}^2$ & 
\centering  $\frac{1}{3}$ & 
\centering  $1$ & 
\centering  $3$ & 
\centering  $\frac{1}{3}$ & 
\centering  $1$ & 
\centering  $3$
\tabularnewline 
\hline \noalign{\smallskip}
$\text{cor}_1$ & 
\centering  \tablenum{0.71} & \centering \tablenum{0.91} & 
\centering \tablenum{0.96} & \centering \tablenum{0.75} & 
\centering \tablenum{0.92} & \centering \tablenum{0.96}   
\tabularnewline  
$\text{cor}_2$ & 
\centering  \tablenum{0.93}  &  \centering \tablenum{0.94} & 
\centering  \tablenum{0.94}  &  \centering \tablenum{0.97} & 
\centering  \tablenum{0.98}  &  \centering \tablenum{0.98}  
\tabularnewline 
\hline \noalign{\smallskip}
$\text{err}_1$ & 
\centering  \tablenum{0.47}  &  \centering \tablenum{0.15} & 
\centering  \tablenum{0.06}  &  \centering \tablenum{0.38} & 
\centering  \tablenum{0.13}  &  \centering \tablenum{0.05}   
\tabularnewline 
$\text{err}_2$ & 
\centering  \tablenum{0.12}  &  \centering \tablenum{0.12} & 
\centering  \tablenum{0.12}  &  \centering \tablenum{0.05} & 
\centering  \tablenum{0.04}  &  \centering \tablenum{0.04}   
\tabularnewline 
$\text{err}_3$ & 
\centering  \tablenum{0.19}  &  \centering \tablenum{0.14} & 
\centering  \tablenum{0.11}  &  \centering \tablenum{0.11} & 
\centering  \tablenum{0.07}  &  \centering \tablenum{0.04}   
\tabularnewline 
\hline \noalign{\smallskip}
$\text{bias} \left( \widehat{\sigma_1^2} \right)$ & 
\centering  \tablenum{-0.02} &  \centering  \tablenum{-0.01} & 
\centering  \tablenum{0.00}  &  \centering  \tablenum{0.02} & 
\centering  \tablenum{0.00}  &  \centering  \tablenum{-0.02}  
\tabularnewline 
$\text{sd} \left( \widehat{\sigma_1^2} \right)$ & 
\centering  \tablenum{1.04}  &  \centering   \tablenum{1.05} & 
\centering  \tablenum{1.06}  &  \centering   \tablenum{0.41} & 
\centering  \tablenum{0.40}   &  \centering  \tablenum{0.39}   
\tabularnewline 
$\text{bias} \left( \widehat{\sigma_2^2} \right)$ & 
\centering  \tablenum{-0.11}  &  \centering  \tablenum{-0.08}  & 
\centering  \tablenum{-0.06}  &  \centering  \tablenum{-0.06}  & 
\centering  \tablenum{-0.06}  &  \centering  \tablenum{-0.06}    
\tabularnewline 
$\text{sd} \left( \widehat{\sigma_2^2} \right)$ & 
\centering  \tablenum{0.50}  &  \centering  \tablenum{0.51} & 
\centering  \tablenum{0.52}  &  \centering  \tablenum{0.21} & 
\centering  \tablenum{0.21}  &  \centering  \tablenum{0.21}   
\tabularnewline 
$\text{bias} \left( \widehat{\sigma_3^2} \right)$ & 
\centering  \tablenum{-0.03}  &  \centering  \tablenum{-0.04}  & 
\centering  \tablenum{-0.04}  &  \centering  \tablenum{-0.02}  & 
\centering  \tablenum{-0.02}  &  \centering  \tablenum{-0.02}   
\tabularnewline 
$\text{sd} \left( \widehat{\sigma_3^2} \right)$ & 
\centering  \tablenum{0.22}  & \centering  \tablenum{0.21} & 
\centering  \tablenum{0.21}  & \centering  \tablenum{0.11} & 
\centering  \tablenum{0.11}  & \centering  \tablenum{0.11} 
\tabularnewline 
\hline
\vspace{5mm}
\end{tabular}
\label{Chauvet::table_new_bundles}
\end{table}

\begin{figure}[!ht]
\centering
\includegraphics[width=.49\linewidth, trim={2.5cm 0 2.25cm 0},clip]{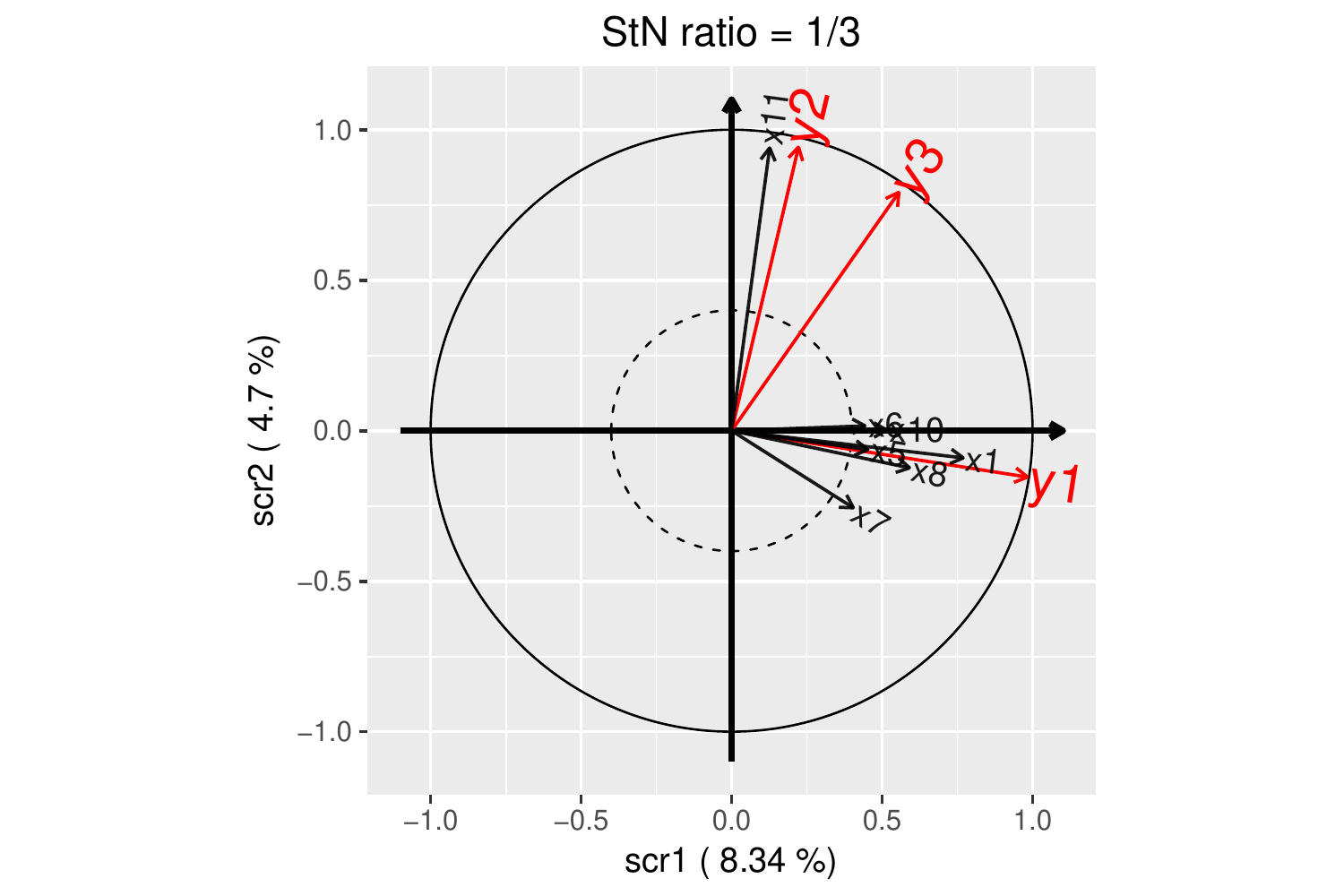}
\includegraphics[width=.49\linewidth, trim={2.5cm 0 2.25cm 0},clip]{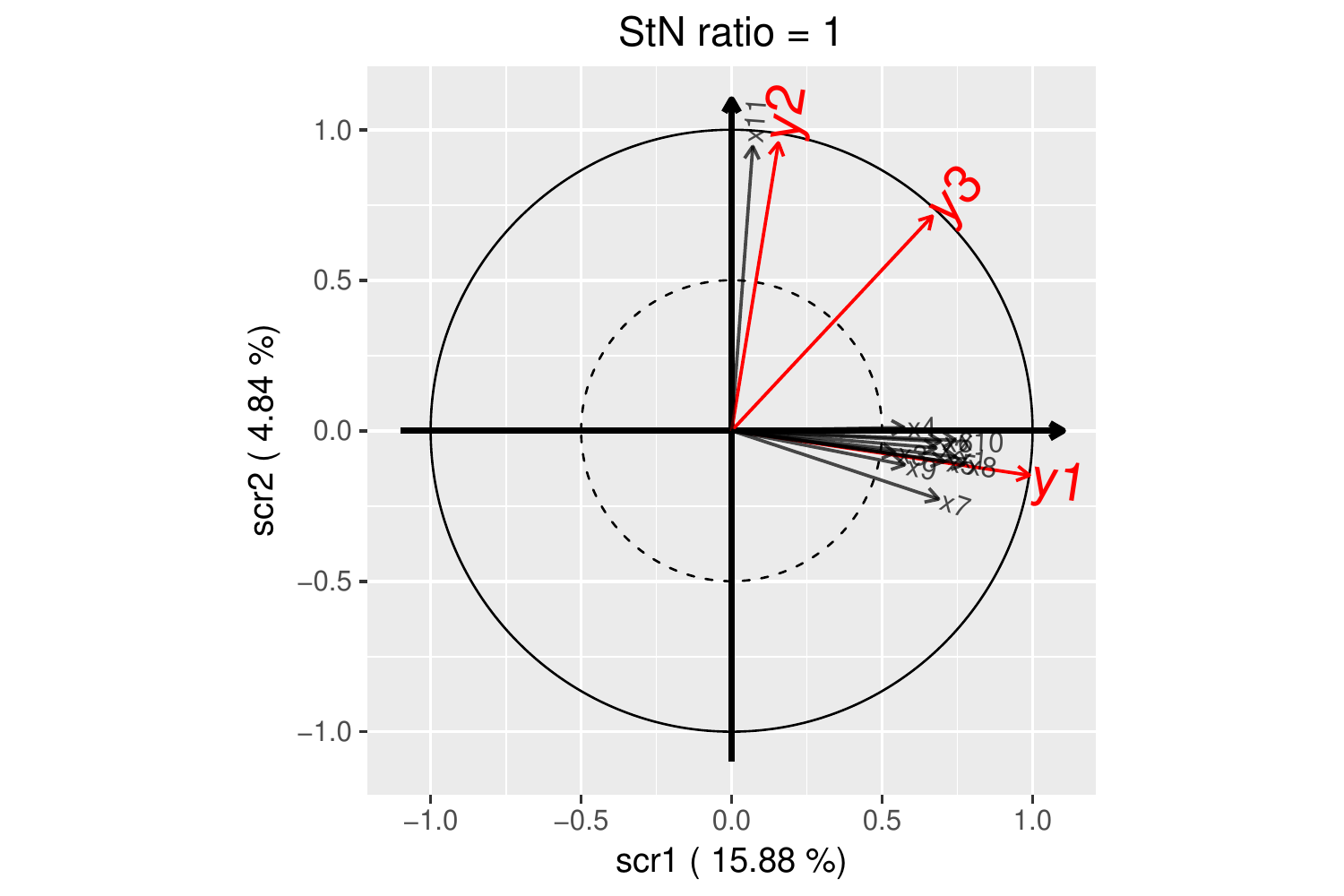}  \\
\includegraphics[width=.49\linewidth, trim={2.5cm 0 2.25cm 0},clip]{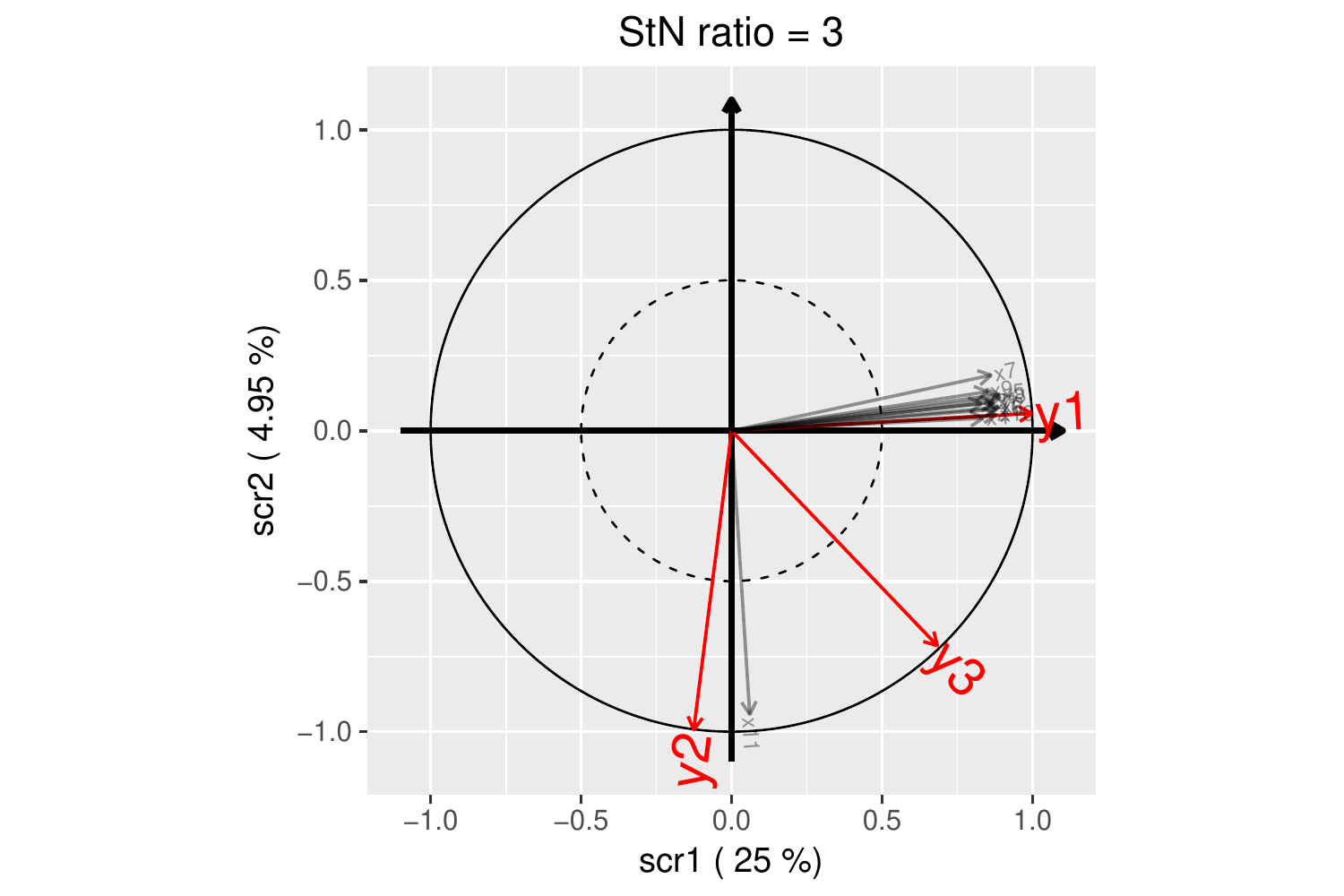} 
\caption{
Examples of the first two--component planes given by mixed--SCGLR when 
$\sigma_{\text{LV}}^2 / \sigma_{\text{noise}}^2 = 1/3$ (top left),  
$\sigma_{\text{LV}}^2 / \sigma_{\text{noise}}^2 = 1$ (top right), and
$\sigma_{\text{LV}}^2 / \sigma_{\text{noise}}^2 = 3$ (bottom). 
When $\text{StN ratio} = 1/3$ (resp. $\text{StN ratio} \in \left\lbrace 1,3 \right\rbrace$), only the variables having cosine greater than $0.4$ (resp. $0.5$) with component plane $(1,2)$ are represented.
}
\label{Chauvet::SNTratio131}
\end{figure}

\paragraph*{}
Model interpretation is revealed by \autoref{Chauvet::SNTratio131} 
in the case of $N=10$ groups and $R=10$ observations per group.
The first component aligns with block $\m{\mathcal{G}_1}$ which alone explains response 
$\m{y_1}$. The second aligns with $\m{\mathcal{G}_2}$ (containing single explanatory variable 
$\m{x_{11}}$) which alone explains $\m{y_2}$.
Finally, note that the projection of the $\m{X}$--part of the linear predictor related to 
$\m{y_3}$ is well represented on component plane $(1,2)$. This indicates that $\m{y_3}$ is explained jointly by $\m{\mathcal{G}_1}$ and $\m{\mathcal{G}_2}$.

\clearpage
\section{High dimensional data}

\subsection{Key idea}
To cope with high dimensional data, the key idea is to replace the fixed--effect design matrix, $\m{X}$, with the matrix $\m{C}$ of its principal components associated with non--zero eigenvalues. 
More precisely, $\lambda_j$ being the eigenvalue associated with the 
$j$--th eigenvector $\m{v_j}$, the last eigenvector we consider, $\m{v_r}$, is such that 
$$
\dfrac{\lambda_r}{\sum_{j=1}^r \lambda_j} > \frac{1}{p},
$$
where $p$ is the number of columns of matrix $\m{X}$.
The matrix of the corresponding unit-eigenvectors is denoted
$\m{V} = \big[\, \m{v_1} \mid \ldots \mid \m{v_r} \,\big]$, and $\m{C}=\m{XV}$. 
The component $\m{f}$ is then sought as a combination of the principal components:
$\m{f}=\m{Cu}=\m{X\widetilde{u}}$, where
$\m{\widetilde{u}}=\m{Vu}$.  
Mixed--SCGLR then solves
\begin{equation*}
\left\lbrace
\begin{aligned}
&\text{max } \qquad\;\; 
s \log \left[\phi \left( \m{u} \right)\right] 
+ (1-s) \log \left[\psi_{\m{A}} \left( \m{u} \right)\right] \\
&\text{subject to } \; \m{u}^{\mbox{\tiny T}} \m{C}^{\mbox{\tiny T}} \m{PC}  \m{u}=1,
\end{aligned}
\right.
\end{equation*}
where the goodness--of--fit measure, $\psi_{\m{A}}$, is given by 
\begin{align*}
\psi_{\m{A}} \left( \m{u} \right) 
&= 
\sum_{k=1}^q \left\lVert \m{z_k}^{\vphantom{t}} \right\rVert_{\m{W^{\xi}_k}}^2  
\cos_{\m{W^{\xi}_k}}^2 
\Big( \m{z_k},\, \spn \left\lbrace \m{Cu},\m{A} \right\rbrace \Big)
\nonumber \\
&=
\sum_{k=1}^q \left\lVert \m{z_k}^{\vphantom{t}} \right\rVert_{\m{W^{\xi}_k}}^2  
\cos_{\m{W^{\xi}_k}}^2 
\left( \m{z_k}, \,
\Pi_{\spn \left\lbrace \m{Cu},\m{A} \right\rbrace}^{\m{W^{\xi}_k}} \m{z_k} 
\right),
\end{align*}
and the structural relevance by
\begin{equation*}
\phi \left( \m{u} \right)
= \left[    
\sum_{j=1}^p \omega_j 
\left(
\spr{\m{Cu}}{\m{x^j}}{\m{P}}^2
\right)^l
\right]^{\frac{1}{l}}
= \left[   
\sum_{j=1}^p \omega_j 
\Big( 
\m{u}^{\mbox{\tiny T}} \, 
\m{C}^{\mbox{\tiny T}} \m{P} \m{x_j} \m{x}_{\m{j}}^{\mbox{\tiny T}} \m{P C} \, 
\m{u} 
\Big)^{l}      
\right]^{\frac{1}{l}}.
\end{equation*}

\noindent
This idea is tested on simulated data where the number of explanatory
variables $p$ exceeds the number of observations $n$.

\subsection{Data generation}
To generate grouped data,
we consider $N=10$ groups and $R=10$ observations per group 
(i.e. a total of $n=100$ observations). 
The random effects' design matrix is then
$\m{U} = \m{I}_N \otimes \mathbf{1}_R$. 
Explanatory variables consist of
four independent bundles  
$\m{X_j}, j \in \left\lbrace 1,2,3,4 \right\rbrace$, such as 
$\m{X} = 
\big[\, \m{X_0} \mid \m{X_1} \mid \m{X_2} \mid \m{X_3} \,\big]$.
Each explanatory variable is normally simulated with mean $0$ and variance $1$. 
Parameter $\tau \in \left\lbrace 0.3, 0.5, 0.7 \right\rbrace$
tunes the level of redundancy within each bundle: 
the correlation matrix of bundle $\m{X_j}$ is
\begin{equation*}
\text{cor} \left(\m{X_j} \right) = \tau \mathbf{1}_{p_j} 
\mathbf{1}_{p_j}^{\mbox{\tiny T}} + (1-\tau) \m{I}_{p_j},
\end{equation*}
where $p_j$ is the number of variables in $\m{X_j}$.
For each $k \in \left\lbrace 1,2,3,4 \right\rbrace$, random--effect vectors are simulated as  
$\m{\xi_k} \widesim{\text{ind.}} \mathcal{N}_N \left( \m{0}, \sigma_k^2 \, \m{I}_N \right)$.
Given $\m{\xi_1}$, $\m{\xi_2}$, $\m{\xi_3}$, $\m{\xi_4}$, 
we simulate 4 responses having different distributions, 
$\m{Y} = \big[\, \m{y_1} \mid \m{y_2} \mid \m{y_3}  \mid \m{y_4} \,\big]$, as
\begin{equation}
\left\lbrace
\begin{aligned}
&\m{y_1} \widesim{} \mathcal{N}_n \Big(  
\m{\mu} = \m{X\beta_1} + \m{U \xi_1}, \,
\m{\Sigma} = \m{I}_n
\Big) \\
&\m{y_3} \widesim{} \mathcal{B} \Big( 
\m{p} = \text{logit}^{-1} \Big[ 
\m{X\beta_2} + \m{U \xi_2}
\Big] \Big) \\
&\m{y_3} \widesim{} \mathcal{B}\text{in} \Big( 
\textbf{trials} = 30 \, \mathbf{1}_n, \,
\m{p} = \text{logit}^{-1} \Big[ 
\m{X\beta_3} + \m{U \xi_3}
\Big] \Big) \\
&\m{y_4} \widesim{} \mathcal{P} \Big( 
\m{\lambda} = \exp \Big[
\m{X\beta_4} + \m{U \xi_4} \Big]
\Big). \\
\end{aligned}
\right.
\label{Chauvet:modele_simulateddata_add}
\end{equation}
Response $\m{y_1}$ is predicted only by bundle $\m{X_1}$, 
$\m{y_2}$ only by bundle $\m{X_2}$, $\m{y_3}$ only by bundle $\m{X_3}$, 
$\m{y_4}$ by both bundles $\m{X_2}$ and $\m{X_3}$,
and bundle $\m{X_0}$ plays no explanatory role.
Our choice for the fixed--effect parameters is
\begin{align*}
\m{\beta_1} &= (\; 
\underbrace{0, \ldots\ldots\ldots\ldots\ldots\ldots, 0}_{p_0 \: \text{times}},
\underbrace{0.1, \ldots\ldots\ldots\ldots\, 0.1}_{p_1 \: \text{times}}, 
\underbrace{0, \ldots\ldots\ldots\ldots, 0}_{p_2 \: \text{times}},
\underbrace{0, \ldots\ldots\ldots, 0}_{p_3 \: \text{times}}
\;)^{\mbox{\tiny T}}, \\
\m{\beta_2} &= (\; 
\underbrace{0, \ldots\ldots\ldots\ldots\ldots\ldots, 0}_{p_0 \: \text{times}},
\underbrace{0, \ldots\ldots\ldots\ldots\ldots\, 0}_{p_1 \: \text{times}}, 
\underbrace{0.1, \ldots\ldots\ldots, 0.1}_{p_2 \: \text{times}},
\underbrace{0, \ldots\ldots\ldots, 0}_{p_3 \: \text{times}}
\;)^{\mbox{\tiny T}}, \\
\m{\beta_3} &= (\; 
\underbrace{0, \ldots\ldots\ldots\ldots\ldots\ldots, 0}_{p_0 \: \text{times}},
\underbrace{0, \ldots\ldots\ldots\ldots\ldots\, 0}_{p_1 \: \text{times}}, 
\underbrace{0, \ldots\ldots\ldots\ldots, 0}_{p_2 \: \text{times}},
\underbrace{0.05, \ldots.., 0.05}_{p_3 \: \text{times}}
\;)^{\mbox{\tiny T}}, \\
\m{\beta_4} &= (\; 
\underbrace{0, \ldots\ldots\ldots\ldots\ldots\ldots, 0}_{p_0 \: \text{times}},
\underbrace{0.025, \ldots\ldots..\, 0.025}_{p_1 \: \text{times}}, 
\underbrace{0.025, \ldots.., 0.025}_{p_2 \: \text{times}},
\underbrace{0, \ldots\ldots\ldots, 0}_{p_3 \: \text{times}}
\;)^{\mbox{\tiny T}}
\end{align*}
We consider in turn $p=150$ ($p_0=60$, $p_1=45$, $p_2=30$, $p_3=15$) and $p=200$ ($p_0=80$, $p_1=60$, $p_2=40$, $p_3=20$) explanatory variables.
Variance components are set to 
$\sigma_1^2 = \sigma_2^2 = \sigma_3^2 = 0.1$, and $\sigma_4^2 = 0.05$. 
For each value of $p$ and for each value of $\tau$, $B=20$ samples are generated according to Model (\ref{Chauvet:modele_simulateddata_add}).

\subsection{Results}
\autoref{Chauvet::table_hdim150} and \autoref{Chauvet::table_hdim200}
present the results for respectively $150$ and $200$ explanatory variables.
They give the Mean Relative Squared Error (MRSE) values for 
$\m{\beta_k}, k \in \left\lbrace 1,\ldots,4\right\rbrace$, as well as biases and standard errors of estimated variance components, obtained on 20 samples for each value of $\tau$.

\begin{table}[!ht]
\caption{Mean Relative Squared Error (MRSE) values for fixed--effect estimates, and biases and standard errors for estimated variance components, obtained with $100$ observations and $150$ explanatory variables.} 
\vspace{3mm}
\centering
\begin{tabular}{lcccccccc}
 & $\m{\beta_1}$ & $\m{\beta_2}$ & $\m{\beta_3}$ & $\m{\beta_4}$ & 
 $\sigma_1^2$ & $\sigma_2^2$ & $\sigma_3^2$ & $\sigma_4^2$ \\ 
 \hline \noalign{\smallskip}
\multirow{2}{*}{$\tau=0.3$} & \multirow{2}{*}{\tablenum{0.06}} & 
\multirow{2}{*}{\tablenum{0.26}} & \multirow{2}{*}{\tablenum{0.19}} & 
\multirow{2}{*}{\tablenum{0.13}} &  
\tablenum{-0.01}  & \tablenum{-0.03}  & \tablenum{-0.02}  & \tablenum{0.02}  \\
&&&&& \tablenum{0.09} & \tablenum{0.09} & \tablenum{0.03} & \tablenum{0.06} \\
\multirow{2}{*}{$\tau=0.5$} & \multirow{2}{*}{\tablenum{0.03}} & 
\multirow{2}{*}{\tablenum{0.20}} & \multirow{2}{*}{\tablenum{0.10}} & 
\multirow{2}{*}{\tablenum{0.07}} &  
\tablenum{0.01}  & \tablenum{-0.03}  & \tablenum{0.00}  & \tablenum{0.01}  \\
&&&&& \tablenum{0.11} & \tablenum{0.08} & \tablenum{0.07} & \tablenum{0.07} \\
\multirow{2}{*}{$\tau=0.7$} & \multirow{2}{*}{\tablenum{0.01}} & 
\multirow{2}{*}{\tablenum{0.10}} & \multirow{2}{*}{\tablenum{0.05}} & 
\multirow{2}{*}{\tablenum{0.04}} &  
\tablenum{0.01}  & \tablenum{-0.05}  & \tablenum{0.01}  & \tablenum{0.02}  \\
&&&&& \tablenum{0.07} & \tablenum{0.09} & \tablenum{0.10} & \tablenum{0.07} 
\tabularnewline 
\hline
\end{tabular}
\label{Chauvet::table_hdim150}
\end{table}

\begin{table}[!ht]
\caption{Mean Relative Squared Error (MRSE) values for fixed--effect estimates, and biases and standard errors for estimated variance components, obtained with $100$ observations and $200$ explanatory variables.} 
\vspace{3mm}
\centering
\begin{tabular}{lcccccccc}
 & $\m{\beta_1}$ & $\m{\beta_2}$ & $\m{\beta_3}$ & $\m{\beta_4}$ & 
 $\sigma_1^2$ & $\sigma_2^2$ & $\sigma_3^2$ & $\sigma_4^2$ \\ 
 \hline \noalign{\smallskip}
\multirow{2}{*}{$\tau=0.3$} & \multirow{2}{*}{\tablenum{0.06}} & 
\multirow{2}{*}{\tablenum{0.15}} & \multirow{2}{*}{\tablenum{0.18}} & 
\multirow{2}{*}{\tablenum{0.10}} &  
\tablenum{-0.04}  & \tablenum{-0.05}  & \tablenum{0.01}  & \tablenum{-0.02}  \\
&&&&& \tablenum{0.04} & \tablenum{0.09} & \tablenum{0.05} & \tablenum{0.05} \\
\multirow{2}{*}{$\tau=0.5$} & \multirow{2}{*}{\tablenum{0.03}} & 
\multirow{2}{*}{\tablenum{0.17}} & \multirow{2}{*}{\tablenum{0.09}} & 
\multirow{2}{*}{\tablenum{0.05}} &  
\tablenum{-0.05}  & \tablenum{0.00}  & \tablenum{-0.02}  & \tablenum{-0.01}  \\
&&&&& \tablenum{0.06} & \tablenum{0.19} & \tablenum{0.04} & \tablenum{0.04} \\
\multirow{2}{*}{$\tau=0.7$} & \multirow{2}{*}{\tablenum{0.01}} & 
\multirow{2}{*}{\tablenum{0.15}} & \multirow{2}{*}{\tablenum{0.04}} & 
\multirow{2}{*}{\tablenum{0.03}} &  
\tablenum{0.03}  & \tablenum{0.00}  & \tablenum{-0.01}  & \tablenum{-0.02}  \\
&&&&& \tablenum{0.08} & \tablenum{0.14} & \tablenum{0.05} & \tablenum{0.05} 
\tabularnewline 
\hline
\end{tabular}
\label{Chauvet::table_hdim200}
\end{table}

\begin{figure}[!ht]
\centering
\includegraphics[width=.49\linewidth, trim={2.5cm 0 2.25cm 0},clip]{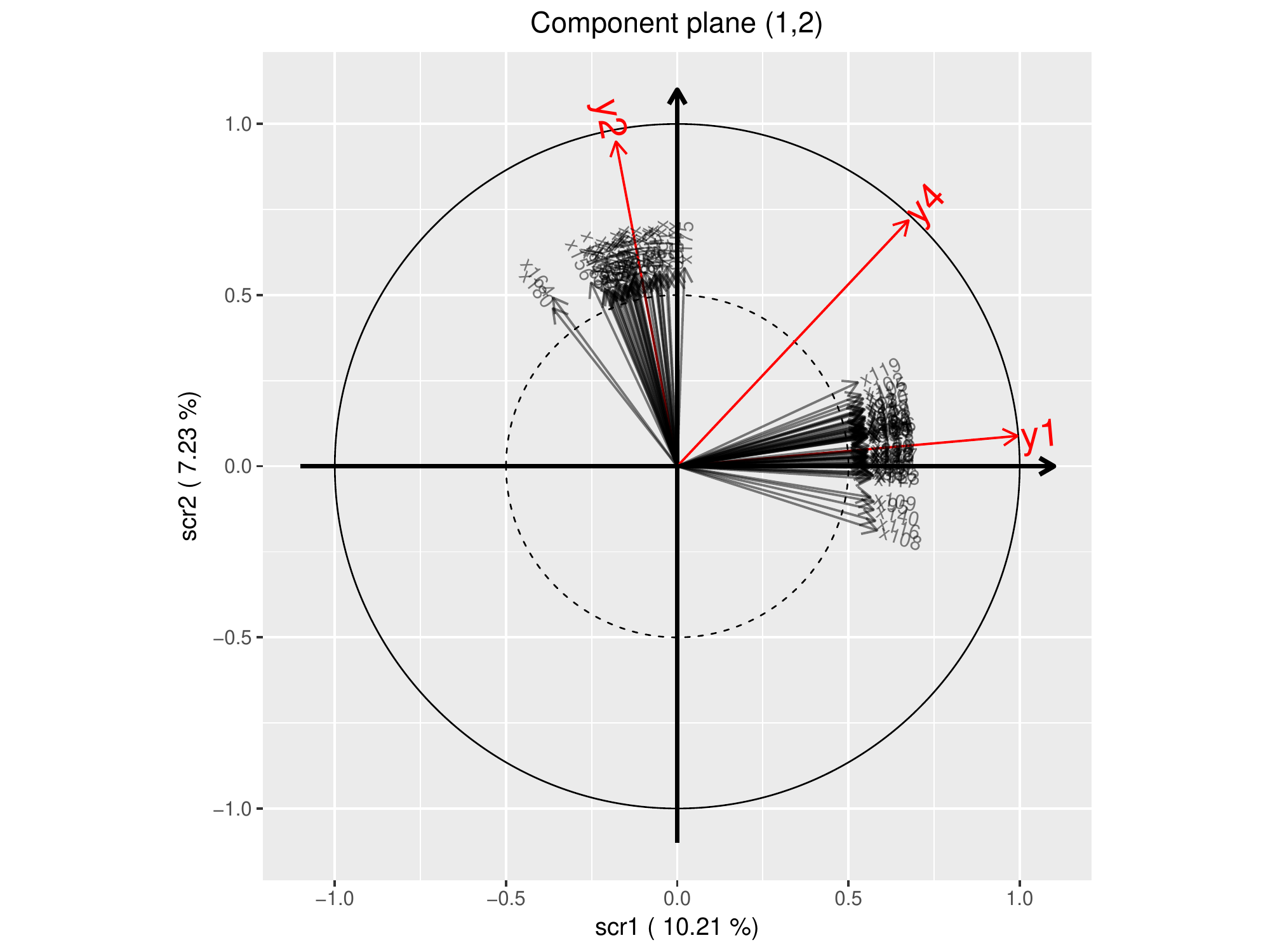}
\includegraphics[width=.49\linewidth, trim={2.5cm 0 2.25cm 0},clip]{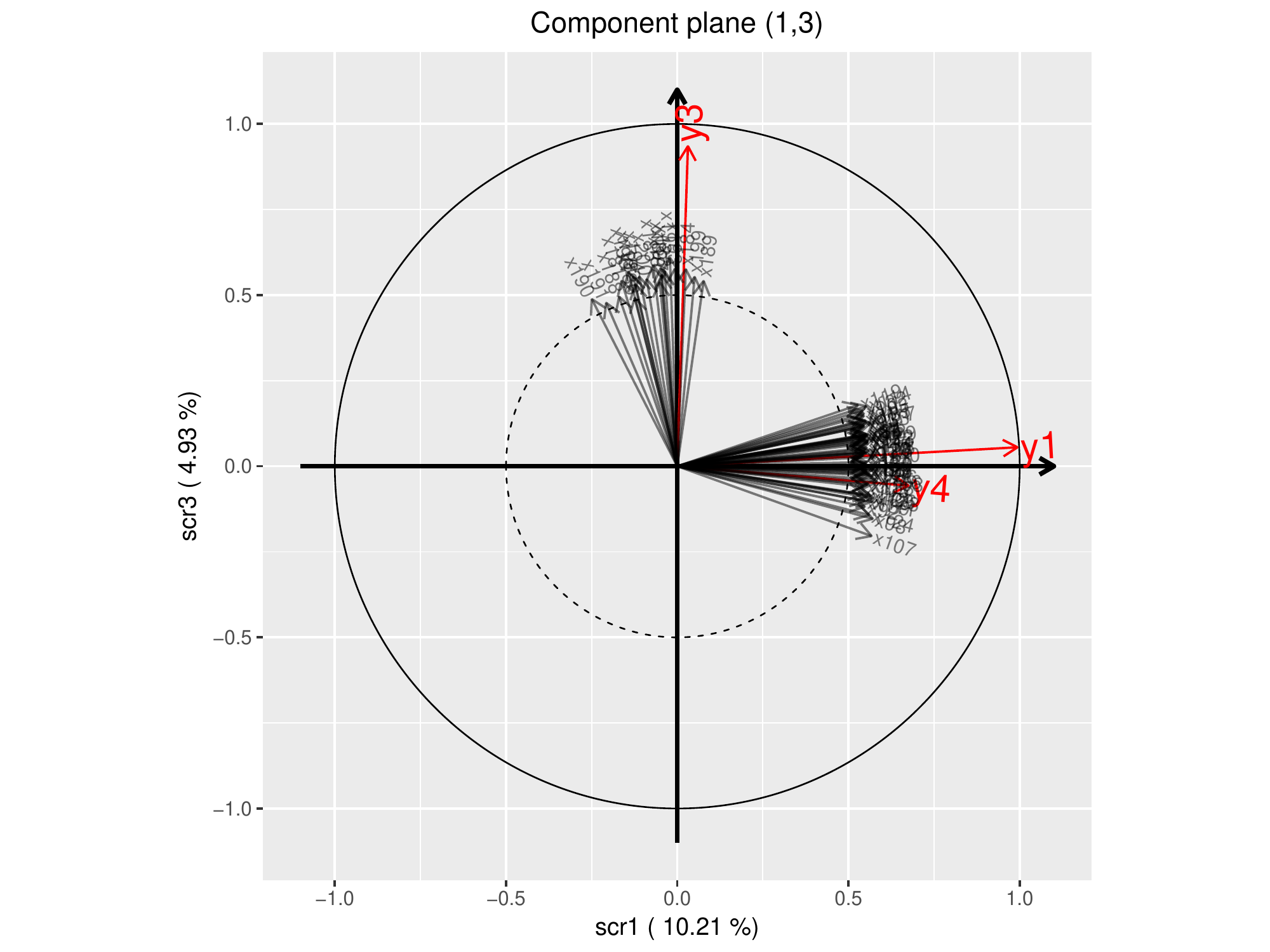}  \\
\includegraphics[width=.49\linewidth, trim={2.5cm 0 2.25cm 0},clip]{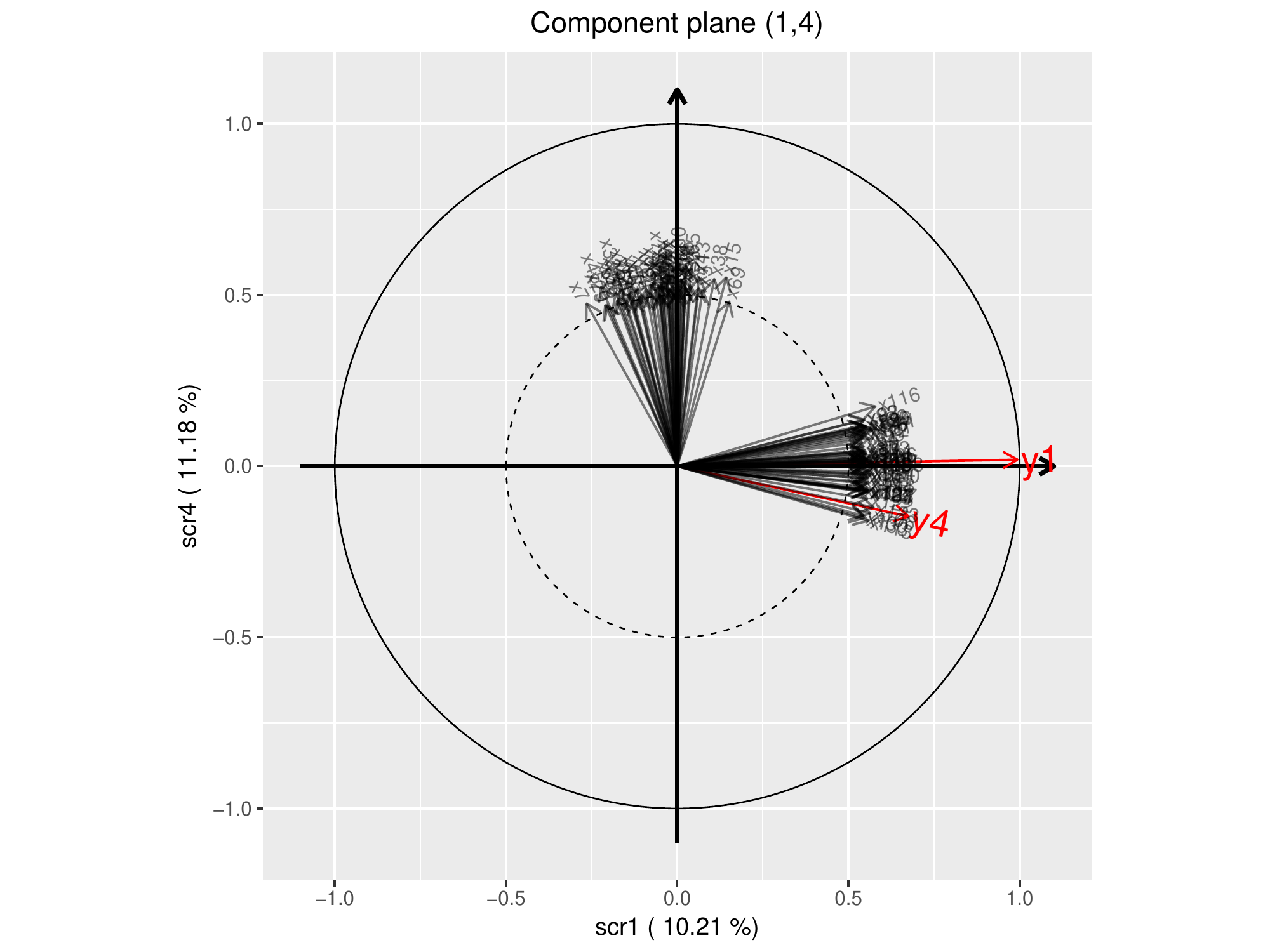} 
\caption{
Component planes $(1,2)$, $(1,3)$ and $(1,4)$ given by mixed--SCGLR
for $100$ observations, $150$ explanatory variables and $\tau=0.3$. 
The tuning parameter triplet $(K,s,l)$ is set to $(4,0.5,4)$.
}
\label{Chauvet::fig150var}
\end{figure}

\begin{figure}[!ht]
\centering
\includegraphics[width=.49\linewidth, trim={2.5cm 0 2.25cm 0},clip]{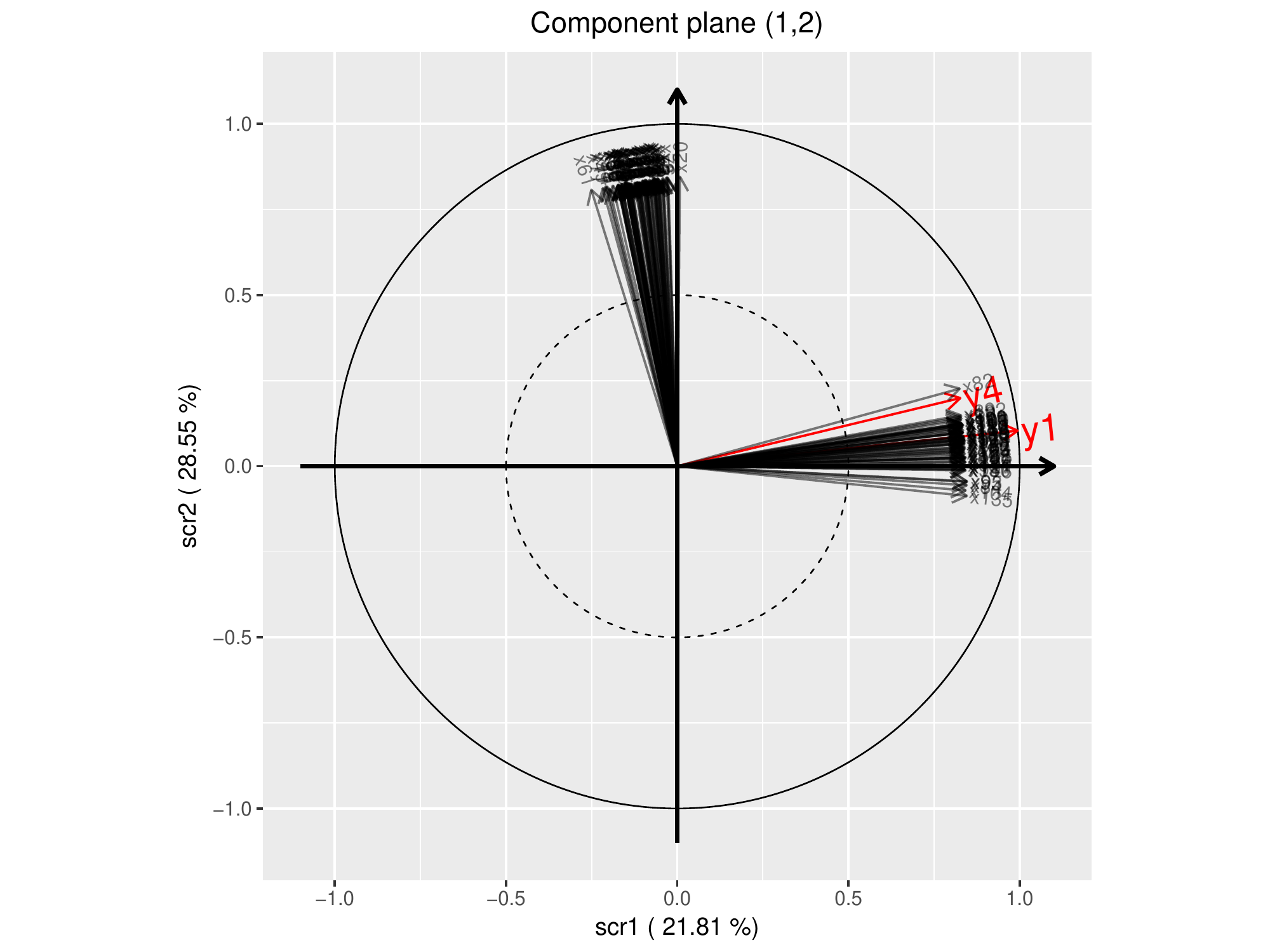}
\includegraphics[width=.49\linewidth, trim={2.5cm 0 2.25cm 0},clip]{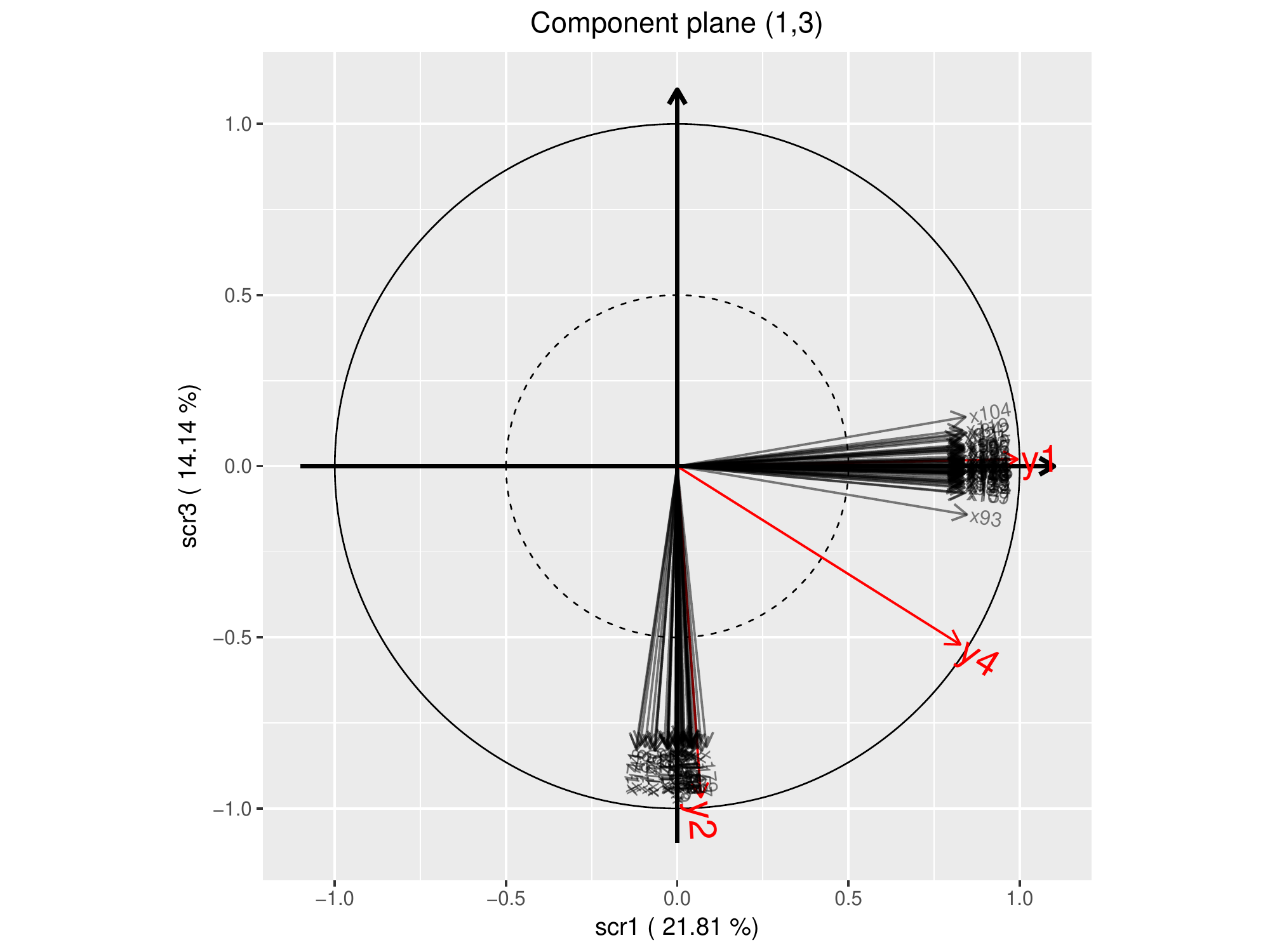}  \\
\includegraphics[width=.49\linewidth, trim={2.5cm 0 2.25cm 0},clip]{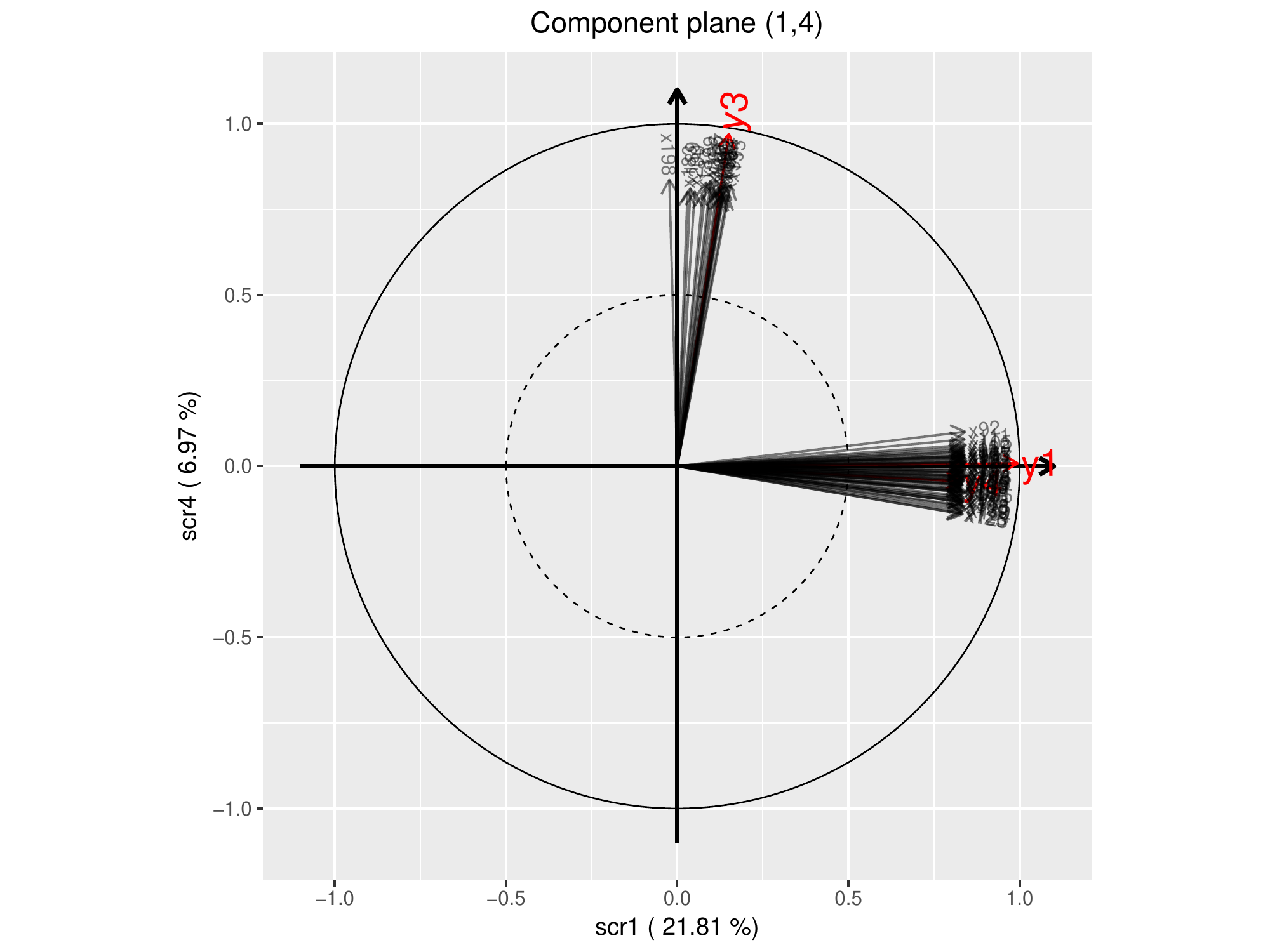} 
\caption{
Component planes $(1,2)$, $(1,3)$ and $(1,4)$ given by mixed--SCGLR
for $100$ observations, $200$ explanatory variables, and $\tau=0.9$.
The tuning parameter triplet $(K,s,l)$ is set to $(4,0.9,4)$.
}
\label{Chauvet::fig200var}
\end{figure}

Some component planes are given on
\autoref{Chauvet::fig150var} ($150$ explanatory variables) and
\autoref{Chauvet::fig200var} ($200$ explanatory variables).

\clearpage
\bibliographystyle{apalike}
\bibliography{chauvetBIBLIO}

\end{document}